\documentclass{amsart}

\usepackage[T1]{fontenc}

\usepackage[dvipsnames]{xcolor}
\usepackage[colorlinks=true,hyperindex, linkcolor=magenta, pagebackref=false, citecolor=cyan,pdfpagelabels]{hyperref}
\hypersetup{linktocpage}

\hypersetup{colorlinks, linkcolor=BurntOrange, urlcolor=BurntOrange, citecolor=NavyBlue}

\usepackage{amssymb}
\usepackage{amsmath}
\usepackage[matrix,arrow,curve,cmtip]{xy}
\entrymodifiers={+!!<0pt,\fontdimen22\textfont2>} 
\usepackage{units}
\usepackage{enumitem}
\setenumerate[1]{leftmargin=2em}
\usepackage{eucal}
\usepackage[nameinlink,capitalise,noabbrev]{cleveref} 

\numberwithin{equation}{section}

\theoremstyle{plain}   

\newtheorem{theorem}[equation]{Theorem}
\newtheorem{corollary}[equation]{Corollary}     
\newtheorem{lemma}[equation]{Lemma}         
\newtheorem{proposition}[equation]{Proposition}
\newtheorem{addendum}[equation]{Addendum}

\theoremstyle{definition}
\newtheorem{definition}[equation]{Definition}

\theoremstyle{remark}
\newtheorem{remark}[equation]{Remark}
\newtheorem{example}[equation]{Example}
\newtheorem{warning}[equation]{Warning}


\begin{document}

\title{Dirac geometry I: Commutative algebra}

\author{Lars Hesselholt}
\address{Nagoya University, Japan, and University of Copenhagen, Denmark}
\email{larsh@math.nagoya-u.ac.jp}

\author{Piotr Pstr\k{a}gowski}
\address{Harvard University, USA, and Institute for Advanced Study, USA}
\email{pstragowski.piotr@gmail.com}

\thanks{Both authors were partially supported by the Danish National Research Foundation through the Copenhagen Center for Geometry and Topology (DNRF151). The first author also received support from JSPS Grant-in-Aid for Scientific Research number 21K03161.}

\maketitle

\setcounter{tocdepth}{2}
\tableofcontents

\section{Introduction}

The purpose of this paper and its sequel is to develop the geometry built from the commutative algebras that naturally appear as the homology of differential graded algebras and, more generally, as the homotopy of algebras in spectra. The commutative algebras in question are those in the symmetric monoidal category of \emph{graded} abelian groups, and, being commutative, they form the affine building blocks of a geometry, as commutative rings form the affine building blocks of algebraic geometry. We name this geometry Dirac geometry, because the grading exhibits the hallmarks of spin in that it is a remnant of the internal structure encoded by \emph{anima}, it distinguishes symmetric and anti-symmetric behavior, and the coherent cohomology of Dirac schemes and Dirac stacks, which we develop in the sequel, admits half-integer Serre twists.
\footnote{\,We remark that Dirac geometry is related to but distinct from supergeometry in that the latter is built from commutative algebras in the symmetric monoidal category of $\mathbb{Z}/2$-graded abelian groups instead of $\mathbb{Z}$-graded abelian groups. The distinction is most relevant in the situation of arithmetic geometry, where Bott periodicity fails. We also note that our usage of the term Dirac geometry is unrelated to the usage in differential geometry as in \cite{courant1990dirac}.}
Thus, informally, Dirac geometry constitutes a ``square root'' of $\mathbb{G}_m$-equivariant algebraic geometry.

As mentioned, the bulding blocks of Dirac geometry are Dirac rings by which we mean anticommutative $\mathbb{Z}$-graded rings. In this paper, we develop the commutative algebra of Dirac rings. We assign to a Dirac ring a Dirac prime spectrum and define a Dirac scheme to be a Dirac ringed space, locally isomorphic to the Dirac prime spectrum of a Dirac ring. Many definitions and results from ordinary geometry extend mutatis mutandis to Dirac geometry. These include the definitions of most geometric properties of maps between Dirac schemes as well as results such as Nakayama's lemma and Hilbert's basis theorem. But there are also some notable differences. The underlying space of the Dirac prime spectrum of a Dirac ring of finite type over a field may not be Jacobson, so Hilbert's Nullstellensatz fails, and more importantly, smooth maps between Dirac schemes may not be flat! However, as a consequence of our main theorem, every \'{e}tale map between Dirac schemes is flat. This theorem involves the key new geometric property of \emph{evenness} of maps in Dirac geometry: A map between Dirac rings is even if it is obtained by extension of scalars from its underlying map of sub-Dirac rings consisting of the elements of even degree. Now, our main theorem states that every \'{e}tale map of Dirac schemes is even. To prove this result, we show that Zariski's main theorem holds in Dirac geometry. The proof consists of an elaborate series of reductions to the even case, which, in turn, follows from Zariski's main theorem in ordinary geometry. Finally, we (re)define a map of commutative algebras in spectra $A \to B$ to be \'{e}tale if the induced map of Dirac rings $\pi_*(A) \to \pi_*(B)$ is \'{e}tale and show that, as a consequence of our main theorem, the functor that to $A \to B$ assigns $\pi_*(A) \to \pi_*(B)$ restricts to an equivalence from the $\infty$-category of \'{e}tale $A$-algebras to the $1$-category of \'{e}tale $\pi_*(A)$-algebras. This generalizes a result of Lurie~\cite[Corollary~7.5.0.6]{lurieha}, who uses a more restrictive notion of \'{e}taleness. For instance, the canonical map
$$\xymatrix{
{ (KU_p)^{h\operatorname{Aut}(\mu_p)} } \ar[r]^-{\theta} &
{ KU_p } \cr
}$$
with $\operatorname{Aut}(\mu_p) \subset \operatorname{Aut}(\mu_{p^{\infty}}) = \mathbb{Z}_p^{\times}$ acting on $KU_{p}$ via Adams operations is \'{e}tale in our sense, but not in Lurie's.

In more detail, a Dirac ring is a commutative algebra in the symmetric monoidal category of $\mathbb{Z}$-graded abelian groups with the symmetry isomorphism
$$x \otimes y \mapsto (-1)^{\deg(x)\deg(y)} y \otimes x,$$
and a map of Dirac rings is a ring homomorphism that preserves the grading. We will use homological grading throughout, and we write
$$\textstyle{ \operatorname{spin}(x) = -\frac{1}{2}\deg(x) }$$
for a homogeneous element $x$. A graded ideal in a Dirac ring is a maximal graded ideal if it is maximal among proper graded ideals, and a Dirac ring is local if it has a unique maximal graded ideal. A local Dirac ring is a Dirac field if its maximal graded ideal is the zero ideal. It is easy to see that a Dirac field either is a field $k$ or a Laurent polynomial algebra $k[x^{\pm1}]$ over a field $k$ generated by a homogeneous element $x$ of nonzero spin $s$ and that $s$ is an integer, unless $k$ is of characteristic $2$. 

\begin{proposition}[{\hyperref[lem:nakayama]{Nakayama's lemma}}]
\label{proposition:nakayama's_lemma_in_introduction}
Let $A$ be a local Dirac ring with maximal graded ideal $\mathfrak{m} \subset A$ and residue Dirac field $k = A/\mathfrak{m}$. If $M$ is a finitely generated graded $A$-module such that $M \otimes_Ak = 0$, then $M$ is zero.
\end{proposition}

For graded commutative rings, an analogous statement is stated as an exercise in the book by Bruns and Herzog~\cite[Exercise~1.5.24]{bruns1998cohen}. The usual proof of Nakayama's lemma uses determinants, which are unavailable in the Dirac setting. However, as indicated in~\cite{matsumura}, one can instead use induction on the number of generators.

In Dirac geometry, finite generation and finite presentation will always refer to families of homogeneous elements. We define a Dirac ring $A$ to be noetherian if every graded ideal $I \subset A$ is finitely generated. This is equivalent to the requirement that the Grothendieck abelian category $\operatorname{Mod}_R(\operatorname{Ab})$ of graded $A$-modules is locally noetherian in the sense of \cite[Definition~C.6.8.5]{luriesagtemp}. Moreover, in this situation, the noetherian objects of $\operatorname{Mod}_R(\operatorname{Ab})$ in the sense of~\cite[Definition~C.6.8.1]{luriesag} are precisely the finitely generated graded $A$-modules.

\begin{proposition}[{\hyperref[proposition:hilbert_basis]{Hilbert's basis theorem}}]
\label{proposition:hilbert's_basis_theorem_in_introduction}
If a Dirac ring $A$ is noetherian, then so is every finitely generated $A$-algebra $B$.
\end{proposition}

We follow the usual proof, but some care must be taken in order to work with homogeneous elements. It turns out that the expression of homogeneous elements of the free Dirac $A$-algebra $B = A[x]$ on a homogeneous generator as homogeneous polynomials is not unique, but this is irrelevant for the argument.

Let $A$ be a Dirac ring. We define a graded $A$-module $M$ to be flat if the functor from the abelian category of graded $A$-modules to itself that to $N$ assigns $M \otimes_AN$ is exact. We prove that, with this definition, the equivalent characterizations of flatness due to Lazard are valid in the Dirac context as well. This also appears in unpublished work of Davies~\cite{davies_lazard_theorem}.

\begin{proposition}[{\hyperref[thm:lazard]{Lazard's theorem}}]
\label{proposition:lazard_in_introduction}
Let $A$ be a Dirac ring, and let $M$ be a graded $A$-module. The following are equivalent.
\begin{enumerate}
\item[{\rm (1)}]The graded $A$-module $M$ is flat.
\item[{\rm (2)}]Given maps of graded $A$-modules
$$\xymatrix{
{ F' } \ar[r]^-{a} &
{ F } \ar[r]^-{x} &
{ M } \cr
}$$
with $F$ and $F'$ finitely generated free and $xa = 0$, there exists a factorization
$$\begin{xy}
(0,12)*+{ {}\phantom{'}F }="11";
(19.5,12)*+{ F'' }="12";
(10,0)*+{ M }="21";
{ \ar@<-.2ex>^-{b} "12";"11";};
{ \ar_-(.35){x} "21";"11";};
{ \ar^-(.35){y} "21";"12";};
\end{xy}$$
with $F''$ finitely generated free and $ba = 0$.
\item[{\rm (3)}]The graded $A$-module $M$ is a filtered colimit of finitely generated free graded $A$-modules.
\end{enumerate}
\end{proposition}

The characterization~(2), known as the equational criterion for flatness, will be of particular importance to us in the sequel~\cite{hesselholtpstragowski2}. We use it here to show that if a graded $A$-module is evenly generated and flat, then it is evenly presented. We note that the notion of flatness employed by Lurie in~\cite[Definition~7.2.2.10]{lurieha} is more restrictive in that, in terms of the characterization~(3), only free modules generated by finite families of homogeneous elements of spin $0$ are allowed.

We include a proof of faithfully flat descent for graded modules. We write $\operatorname{Ab}$ for the symmetric monoidal category of graded abelian groups with the Koszul sign in the symmetry isomorphism so that the category of Dirac rings is the category $\operatorname{CAlg}(\operatorname{Ab})$ of commutative algebras in $\operatorname{Ab}$. The assignment to a Dirac ring $A$ of its category $\operatorname{Mod}_A(\operatorname{Ab})$ of graded modules is a functor
$$\xymatrix@C=10mm{
{ \operatorname{CAlg}(\operatorname{Ab}) } \ar[r]^-{\phantom{,}\operatorname{Mod}\phantom{,}} &
{ \operatorname{LPr} } \cr
}$$
to the $\infty$-category of presentable $\infty$-categories and left adjoint functors. It takes values in the full subcategory spanned by the presentable $1$-categories. Given a map of Dirac rings $\phi \colon A \to B$, we obtain the augmented cosimplicial diagram
$$\xymatrix@C=12mm{
{ \Delta_+ } \ar[r]^-{B^{\otimes_{\!A}[-]}} &
{ \operatorname{CAlg}(\operatorname{Ab}) } \cr
}$$
that to $[n] = \{0,1,\dots,n\}$ assigns the $(n+1)$-fold tensor product $B^{\otimes_{\!A}[n]}$. Now, Grothendieck's faithfully flat descent is the following statement.

\begin{proposition}[{\hyperref[prop:flat_descent_for_abelian_categories_of_modules]{Faithfully flat descent for modules}}]
\label{prop:flat_descent_for_abelian_categories_of_modules_intro}
If $\phi \colon A \to B$ is a faithfully flat map of Dirac rings, then the canonical maps
$$\xymatrix{
{ \operatorname{Mod}_A(\operatorname{Ab}) } \ar[r] &
{ \varprojlim_{[n] \in \Delta} \operatorname{Mod}_{B^{\otimes_{\!A}[n]}}(\operatorname{Ab}) } \ar[r] &
{ \varprojlim_{[n] \in \Delta_{\leq 2}} \operatorname{Mod}_{B^{\otimes_{\!A}[n]}}(\operatorname{Ab}) } \cr
}$$
are equivalences of categories.
\end{proposition}

The right-hand category in \cref{prop:flat_descent_for_abelian_categories_of_modules_intro} is the category of graded $B$-modules with descent data along the map $\phi \colon A \to B$, and the result is commonly proved by showing that the composite functor is an equivalence. The argument that the left-hand functor is an equivalence, however, is much clearer, and the fact that the right-hand functor is an equivalence is an instance of the folklore fact that the limit of a cosimplicial diagram of $n$-categories agrees with the limit of its restriction along the inclusion $\Delta_{\leq n+1} \to \Delta$. We include a proof of this fact in the appendix.

To discuss further results, we associate to a Dirac ring $A$ its Dirac prime spectrum given by a pair of a Dirac ringed space
$$X = (|X|,\mathcal{O}_X)$$
and a map of Dirac rings $\epsilon_A \colon A \to \mathcal{O}_X(|X|)$. The topological space $|X|$ is the set of graded prime ideals $\mathfrak{p} \subset A$ with the Zariski topology. A basis for the topology is given by the distinguished open subsets $|X_f| \subset |X|$ with $f \in A$ homogeneous consisting of the graded prime ideals $\mathfrak{p} \subset A$ such that $f \notin \mathfrak{p}$. The space $|X|$, which we call the Dirac--Zariski space of $A$, is well-known and appears in Quillen's work on equivariant cohomology~\cite{quillen7}, and is an important tool in tensor-triangular geometry \cite{balmer2010spectra, benson2008local}. But to our knowledge, the sheaf of Dirac rings $\mathcal{O}_X$ has not been considered before. To define it, we follow Grothendieck's definition of the structure sheaf in ordinary geometry.

\begin{theorem}[{\hyperref[thm:structuresheaf]{Grothendieck's theorem}}]
\label{thm:structuresheaf_in_introduction}
If $A$ is a Dirac ring, and if $|X|$ is its Dirac--Zariski space, then, up to unique isomorphism, there is a unique pair $(\mathcal{O}_X,\epsilon_A \colon A \to \mathcal{O}_X(|X|)$ of a sheaf of Dirac rings on $|X|$ and map of Dirac rings such that for every homogeneous element $f \in A$, the composite map
$$\xymatrix{
{ A } \ar[r]^-{\epsilon_A} &
{ \mathcal{O}_X(|X|) } \ar[r] &
{ \mathcal{O}_X(|X_f|) } \cr
}$$
is a localization with respect to $S = \{1,f,f^2,\dots\} \subset A$.
\end{theorem}

As in ordinary geometry, the Dirac prime spectrum $X = (|X|,\mathcal{O}_X)$ is a locally Dirac ringed space in the sense that for every $x \in |X|$, the stalk $\mathcal{O}_{X,x}$ is a local Dirac ring. The category of locally Dirac ringed spaces is the non-full subcategory of the category of Dirac ringed spaces, whose objects are the locally Dirac ringed spaces, and whose maps are the maps of ringed spaces
$$\xymatrix@C=14mm{
{ Y = (|Y|,\mathcal{O}_Y) } \ar[r]^-{f = (p,\phi)} &
{ X = (|X|,\mathcal{O}_X) } \cr
}$$
such that for every $y \in |Y|$ with image $x = p(y) \in |X|$, the induced map
$$\xymatrix@C=10mm{
{ \mathcal{O}_{X,x} } \ar[r]^-{\phi_x} &
{ \mathcal{O}_{Y,y} } \cr
}$$
is local in the sense that $\phi_x(\mathfrak{m}_x) \subset \mathfrak{m}_y$. An advantage of the Dirac--Zariski space is that it is often quite small, yet captures the essential properties of the Dirac ring. For example, faithfulness of flat maps of Dirac rings is equivalent to surjectivity of the induced map of Dirac--Zariski spaces. 

\begin{definition}
\label{definition:dirac_scheme}
A locally Dirac ringed space $X = (|X|,\mathcal{O}_X)$ is a Dirac scheme if, locally on $|X|$, it is isomorphic to the Dirac prime spectrum of a Dirac ring. The category of Dirac schemes is the full subcategory of the category of locally Dirac ringed spaces spanned by the Dirac schemes.
\end{definition}

The geometric properties of maps of Dirac schemes of being an open and closed immersion, quasi-compact and quasi-separated, locally of finite type and locally finitely presented, and of being \'{e}tale and smooth are defined analogously to the corresponding geometric properties of a map of schemes. We define a map
$$\xymatrix@C=14mm{
{ Y = (|Y|,\mathcal{O}_Y) } \ar[r]^-{f = (p,\phi)} &
{ X = (|X|,\mathcal{O}_X) } \cr
}$$
to be flat if $\phi_x \colon \mathcal{O}_{X,x} \to \mathcal{O}_{Y,y}$ is flat for every $y \in |Y|$ with image $x = p(y) \in |X|$.

The anti-symmetric behavior of homogeneous elements of half-integer spin means that the underlying ring of a Dirac ring $A$ typically is not commutative. However, the sub-Dirac ring $A^{\operatorname{ev}} \subset A$ spanned by the homogeneous elements of integer spin is a $\mathbb{Z}$-graded commutative ring. Moreover, the map of Dirac--Zariski spaces
$$\xymatrix{
{ |X| } \ar[r] &
{ |X^{\operatorname{ev}}| } \cr
}$$
induced by $A^{\operatorname{ev}} \to A$ is a homeomorphism, and we identify the common space with the space of orbits for the $\mathbb{G}_m$-action on the Zariski space $|Y|$ of the underlying commutative ring of $A^{\operatorname{ev}}$ defined by the $\mathbb{Z}$-grading. This space is often quite small. For instance, if $k$ is a field and $A = k[t]$ the free Dirac $k$-algebra on generator $t$ of nonzero spin, then its Dirac--Zariski space is the Sierpinski space 
$$|\mathbb{A}_k^1(t)| = \{\eta,s\}.$$
Thus, Hilbert's Nullstellensatz fails in Dirac geometry: The Dirac--Zariski space of a finitely generated Dirac algebra over a field need not be Jacobson.\footnote{\,A topological space is Jacobson if every closed subset is equal to the closure of its intersection with the subset of closed points.} Nevertheless, the Dirac--Zariski space $|X|$ carries the essence of the topological information encoded by the larger Zariski space $|Y|$. In general, there is a natural map
$$\xymatrix@C=10mm{
{ X = (|X|,\mathcal{O}_X) } \ar[r]^-{\eta_X} &
{ X^{\operatorname{ev}} = (|X|,\mathcal{O}_X^{\operatorname{ev}}) } \cr
}$$
given by the canonical inclusion $\mathcal{O}_X^{\operatorname{ev}} \to \mathcal{O}_X$, and we define:

\begin{definition}
\label{definition_even_map}
A map $f \colon Y \to X$ of Dirac schemes is even, if the diagram
$$\xymatrix@C=10mm{
{ Y } \ar[r]^{\eta_Y} \ar[d]^-(.45){f\phantom{\operatorname{ev}}} &
{ Y^{\operatorname{ev}} } \ar@<-.7ex>[d]^-(.45){f^{\operatorname{ev}}} \cr
{ X } \ar[r]^-{\eta_X} &
{ X^{\operatorname{ev}} } \cr
}$$
is cartesian.
\end{definition}

To be even is a geometric property of maps of Dirac schemes in the sense that it preserved under base-change along arbitrary maps of Dirac schemes. It is also local on the source and target and closed under composition.

To prove statements about maps of Dirac schemes, we employ the following strategy using evenness: Suppose that $P$ and $Q$ are geometric properties of maps of Dirac schemes and that we wish to prove that $P$ implies $Q$. If $P$ and $Q$ are local on both source and target, then it suffices to do so for maps between Dirac rings. Suppose that there are geometric properties $P_0$ and $Q_0$ of maps of schemes such that a map of even Dirac rings has property $P$ and $Q$ if and only if the map of underlying rings has property $P_0$ and $Q_0$, respectively. In this situation, if we know that $P_0$ implies $Q_0$, then we can conclude that $P$ implies $Q$ for all even maps between Dirac schemes. Thus, to prove that $P$ implies $Q$ for general maps between Dirac schemes, we may try to bootstrap our way from the case of even maps between Dirac schemes.

We use this strategy to prove that \'{e}tale maps between Dirac rings are flat. The fact from ordinary geometry that \'{e}tale maps between schemes are flat implies that even \'{e}tale maps between Dirac schemes are flat. Hence, it remains to bootstrap our way to general \'{e}tale maps between Dirac schemes. In fact, we will show that every \'{e}tale map between Dirac schemes is even. The key to proving this fact is to prove that Zariski's main theorem holds in Dirac geometry. The proof of this, in turn, uses the same strategy, albeit in a more refined way. We define a map of Dirac rings $\phi \colon A \to B$ to be finite if it exhibits $B$ as a finitely generated graded $A$-module and to be quasi-finite if it is of finite type and if the fiber
$$\xymatrix{
{ k(x) } \ar[r] &
{ B \otimes_Ak(x) } \cr
}$$
is a finite map of Dirac rings for every $x \in |X| = \lvert\operatorname{Spec}(A)\rvert$.

\begin{theorem}[{\hyperref[corollary:finite_open_immersion_factorization_of_qf_morphism]{Zariski's main theorem}}]
\label{thm:zariski's_main_theorem_in_introduction}
If $\phi \colon A \to B$ is a quasi-finite map of Dirac rings, then the induced map of Dirac prime spectra  $f \colon Y \to X$ admits a factorization as the composition
$$\xymatrix{
{ Y } \ar[r]^-{j_0} &
{ Y_0 } \ar[r]^-{p_0} &
{ X } \cr
}$$
of an open immersion $j_0$ and a finite map $p_0$.
\end{theorem}

Now, to show that \'{e}tale maps are even, we proceed as follows. It suffices to show that every \'{e}tale map $\phi \colon A \to B$ between Dirac rings is even. The usual argument shows that it is quasi-finite and that for every $x \in |X|$, the base-change
$$\xymatrix@C=10mm{
{ k(x) } \ar[r]^-{\phi(x)} &
{ B \otimes_Ak(x) } \cr
}$$
is the product of a finite family of \'{e}tale extensions of Dirac fields. Moreover, the classification of Dirac fields above immediately shows that every \'{e}tale extension of Dirac fields is even. We wish to conclude that the composite map
$$\xymatrix{
{ A } \ar[r]^-{\phi} &
{ B } \ar[r]^-{\gamma} &
{ B_{\mathfrak{q}} } \cr
}$$
is evenly generated for every graded prime ideal $\mathfrak{q} \subset B$, and it is here that we use Zariski's main theorem to bring ourselves in a situation, where Nakayama's lemma applies. This proves our main result:

\begin{theorem}[{\hyperref[theorem:equivalence_of_categories_between_etale_extensions_of_a_and_even_subring]{Evenness of \'{e}tale maps}}]
\label{thm:etale_is_even_in_introduction}
Every \'{e}tale map between Dirac schemes is even, and hence, flat.
\end{theorem}

It follows that the flat and \'{e}tale topologies on the category $\operatorname{Aff}$ of affine Dirac schemes are comparable. The former is finer than the latter, and both are finitary. One may also prove that, \'{e}tale locally, every smooth map of Dirac schemes admits a section, so that the \'{e}tale and smooth topologies on the category $\operatorname{Aff}$ agree. 

We note that smooth maps of Dirac schemes generally are not flat. The problem is that affine space is not flat. For example, the unique map $\mathbb{Z} \to \mathbb{Z}[t]$ to the free Dirac ring on a generator $t$ of half-integer spin is not flat, because
$$\mathbb{Z}[t] \simeq \mathbb{Z} \cdot 1 \oplus \mathbb{Z} \cdot t \oplus \mathbb{Z}/2 \cdot t^2 \oplus \dots \oplus \mathbb{Z}/2 \cdot t^n \oplus \dots$$
as a graded abelian group.

As an application of \cref{thm:etale_is_even_in_introduction}, we prove the following generalization of a theorem of Lurie~\cite[Theorem~7.5.0.6]{lurieha}. If $R$ is an $\mathbf{E}_{k+1}$-algebra in spectra with $1 \leq k \leq \infty$, then $\operatorname{Mod}_R(\operatorname{Sp})$ can be promoted to an $\mathbf{E}_k$-monoidal category, so we may consider the $\infty$-category $\operatorname{Alg}_{\mathbf{E}_{k}}(\operatorname{Mod}_{R}(\operatorname{Sp}))$ of $\mathbf{E}_k$-algebras.

We (re)define an $\mathbf{E}_{k}$-algebra $A$ in $\operatorname{Mod}_R(\operatorname{Sp})$ to be \'{e}tale if $\pi_*(A)$ is a Dirac ring, which is automatic for $k \geq 2$, and if the map $\pi_*(R) \to \pi_*(A)$ induced by the unit map is an \'{e}tale map of Dirac rings. We say ``redefine'' because Lurie also requires that $\pi_0(A) \otimes_{\pi_0(R)}\pi_n(R) \to \pi_n(A)$ be an isomorphism for all integers $n$, and it is in this sense that the following result generalizes Lurie's, as well as some of earlier work of Baker--Richter and Rognes in the context of Galois extensions \cite{baker2007realizability, rognes2008galois}.

\begin{theorem}[{\hyperref[theorem:dirac_etale_extensions_of_ek_rings]{\'{E}tale rigidity}}]
\label{theorem:etale_E_k-algebras_in_introduction}
If $R$ is an $\mathbf{E}_{k+1}$-algebra in spectra with $1 \leq k \leq \infty$, then taking homotopy groups
$$\xymatrix@C=10mm{
{ \operatorname{Alg}_{\mathbf{E}_k}(\operatorname{Mod}_R(\operatorname{Sp}))^{\text{{\rm \'{e}t}}} } \ar[r]^-{\pi_*} &
{ \operatorname{CAlg}(\operatorname{Mod}_{\pi_*(R)}(\operatorname{Ab}))^{\text{{\rm\'{e}t}}} } \cr
}$$
is an equivalence from the full subcategory of $\operatorname{Alg}_{\mathbf{E}_k}(\operatorname{Mod}_R(\operatorname{Sp}))$ spanned by the \'{e}tale algebras to the full subcategory of $\operatorname{CAlg}(\operatorname{Mod}_{\pi_*(R)}(\operatorname{Ab}))$ spanned by the \'{e}tale algebras.
\end{theorem}

Finally, we mention another instance of the square root phenomenon encoded by Dirac geometry. The second author has introduced the synthetic deformation of stable homotopy theory with respect to any homology theory~\cite{patchkoria2021adams, pstrkagowski2018synthetic}. In the case of complex cobordism, the resulting symmetric monoidal stable $\infty$-category of synthetic spectra contains ($p$-complete) cellular motivic spectra over $\mathbb{C}$ as a symmetric monoidal subcategory \cite{pstrkagowski2018synthetic}. We note that the inclusion is proper, and that, notably, the projective line acquires a tensor-square root
$$\mathbb{P}^1 \simeq S^{2,1} \simeq (S^{1,\frac{1}{2}})^{\otimes 2}$$
in the larger $\infty$-category of synthetic spectra.\footnote{\,Here, we use the motivic grading convention. In terms of the notation used in \cite{pstrkagowski2018synthetic}, the square root of $\mathbb{P}^{1} \simeq \nu(S^{2})$ would be denoted by $\nu(S^{1})$.}

\noindent\textbf{Notation.} We implicitly use the language of $\infty$-categories of Joyal and Lurie, for which the standard reference is \cite{luriehtt}, but with the exception of \cref{subsection:etale_maps_of_ring_spectra} most of the paper requires only ordinary 1-category theory and commutative algebra. We use the term \emph{anima} for what is called a space in~\cite{luriehtt}, and we write $\mathcal{S}$ for the $\infty$-category of anima. The inclusion $\operatorname{Set} \simeq \tau_{\leq 0}\mathcal{S} \to \mathcal{S}$ preserves limits, but not colimits. In fact, the $1$-category of sets is the free cocomplete $1$-category on a single generator, whereas the $\infty$-category of anima the free cocomplete $\infty$-category on a single generator.

We define the $\infty$-category of $\mathbb{Z}$-graded objects in an $\infty$-category $\mathcal{C}$ to be the functor $\infty$-category $\operatorname{Fun}(\mathbb{Z},\mathcal{C})$, where we consider $\mathbb{Z}$ as a category with identity maps only. We will write $A = (A_{k})_{k \in \mathbb{Z}}$ to indicate a $\mathbb{Z}$-graded object in $\mathcal{C}$.

\begin{remark}
\label{remark:graded_objects_in_an_infinity_category}
If the $\infty$-category $\mathcal{C}$ admits $\mathbb{Z}$-indexed colimits, then the functor $\operatorname{Fun}(\mathbb{Z}^{\triangleright},\mathcal{C}) \to \operatorname{Fun}(\mathbb{Z},\mathcal{C})$ given by restriction along $\mathbb{Z} \to \mathbb{Z}^{\triangleright}$ induces an equivalence from the full subcategory of the domain spanned by the colimit diagrams. Therefore, it is immaterial whether or not we include the data of a colimit in the definition of a $\mathbb{Z}$-graded object, as does Bourbaki~\cite[Chapter~II, \S11, Definition~1]{bourbaki1}.
\end{remark}\begin{warning}
\label{warning:category_of_graded_abelian_groups_notation}
To avoid cluttering notation, we simply write $\operatorname{Ab}$ for the symmetric monoidal category of $\mathbb{Z}$-graded abelian groups with the Koszul sign in the symmetry isomorphism. This does not lead to conflicts here, since we will have no occasion to consider the symmetric monoidal category of abelian groups.
\end{warning}

\noindent\textbf{Acknowledgements.} We would like to express our gratitude to Bhargav Bhatt, Sanath Devalapurkar, Mel Hochster, Jacob Lurie, Denis Nardin, Ryszard Nest, Maxime Ramzi, and Dylan Wilson for enlightening conversations related to this work. The second author would also like to thank Kazuhiro Fujiwara and Nagoya University and Takeshi Saito and the University of Tokyo for their support and hospitality during his visit to Japan, where this paper was completed. 
\section{Dirac schemes}
\label{sec:diracschemes}

In this section , we define the category of Dirac schemes, emulating the definition of the category of schemes due to Grothendieck. So we first associate to a Dirac ring $A$ its Dirac spectrum, which is a locally Dirac ringed space
$$\operatorname{Spec}(A) = (|X|,\mathcal{O}_X)$$
with underlying $|X|$ given by the set of graded prime ideals $\mathfrak{p} \subset A$ with the Zariski topology. The structure sheaf $\mathcal{O}_X$ is such that its stalk $\mathcal{O}_{X,x}$ at a point $x \in |X|$ corresponding to a graded prime ideal $\mathfrak{p} \subset A$ is the local Dirac ring given by the localization $A_{\mathfrak{p}}$ of $A$ with respect to the set of homogeneous elements $f \in A \smallsetminus \mathfrak{p}$. We next define a Dirac scheme to be a locally Dirac ringed space that locally is isomorphic to $\operatorname{Spec}(A)$ for some Dirac ring $A$.

\subsection{Dirac rings and ideals}
\label{subsec:diracrings}

Let $\operatorname{Sp}$ be the symmetric monoidal $\infty$-category of spectra. It is the initial presentably symmetric monoidal stable $\infty$-category. We define the $\infty$-category of $\mathbb{Z}$-graded spectra to be the $\infty$-category $\operatorname{Fun}(\mathbb{Z},\operatorname{Sp})$, where we view $\mathbb{Z}$ as category with identity maps only. It promotes via Day convolution to a presentably symmetric monoidal stable $\infty$-category.

We recall that the stable $\infty$-category $\operatorname{Fun}(\mathbb{Z},\operatorname{Sp})$ admits the Beilinson $t$-structure, where $X = (X_{k})_{k \in \mathbb{Z}}$ is defined to be connective if $X_k \in \operatorname{Sp}_{\geq k}$ for all $k \in \mathbb{Z}$. The connective part $\operatorname{Fun}(\mathbb{Z},\operatorname{Sp})_{\geq 0}$ is closed under the tensor product and contains the monoidal unit, and hence, the heart $\operatorname{Fun}(\mathbb{Z},\operatorname{Sp})^{\heartsuit}$ can be promoted canonically to a a symmetric monoidal category.

Let us denote by $\operatorname{Ab}$ the category of $\mathbb{Z}$-graded abelian groups. The functor 
$$\xymatrix{
{ \operatorname{Fun}(\mathbb{Z}, \operatorname{Sp})^{\heartsuit} } \ar[r] &
{ \operatorname{Ab} } \cr
}$$
given by $(X_{k})_{k \in \mathbb{Z}} \mapsto (\pi_k(X_{k}))_{k \in \mathbb{Z}}$ is an equivalence of categories, through which $\operatorname{Ab}$ inherits a symmetric monoidal structure from that of graded spectra. Explicitly, the monoidal product is given by the usual formula
$$\textstyle{ (A \otimes B)_{k} \simeq \bigoplus_{i+j = k} A_{i} \otimes B_{j} }$$
with the symmetry isomorphism is given by
$$a \otimes b \mapsto (-1)^{\deg(a) \deg(b)} b \otimes a.$$

\begin{definition}\label{definition:dirac_ring}
A Dirac ring is a commutative algebra in the symmetric monoidal category $\operatorname{Ab}$ of $\mathbb{Z}$-graded abelian groups. The category of Dirac rings is the category $\operatorname{CAlg}(\operatorname{Ab})$ of commutative algebras in $\operatorname{Ab}$.
\end{definition}

More concretely, a Dirac ring is a $\mathbb{Z}$-graded ring in which the multiplication satisfies the ``anticommutative'' law that
$$a \cdot b = (-1)^{\deg(a) \deg(b)} b \cdot a.$$
A Dirac ring $A = (A_k)_{k \in \mathbb{Z}}$ has an underlying ring $\smash{\widetilde{A}} = \bigoplus_{k \in \mathbb{Z}}A_k$. We stress, however, that, due to the Koszul sign, the underlying ring is not a commutative ring, unless $A$ is concentrated in even degrees or $2 = 0$ in $A$.

\begin{example}
\label{example:homotopy_groups_of_einfty_ring_form_a_dirac_ring}
If $E$ is a commutative algebra in spectra, then the homotopy groups $\pi_*(E)$ form a Dirac ring. To see this, note that we can identify $\pi_*(E)$ with the image of $E$ under the composite 
$$\xymatrix{
{ \operatorname{Sp} } \ar[r] &
{ \operatorname{Fun}(\mathbb{Z},\operatorname{Sp}) } \ar[r] &
{ \operatorname{Fun}(\mathbb{Z},\operatorname{Sp})^{\heartsuit} } \ar[r] &
{ \operatorname{Ab}, } \cr
}$$
where the first map forms the constant graded spectrum, the second takes zeroth homology with respect to the Beilinson $t$-structure, and the last is the equivalence defined above. Each map promotes to a lax symmetric monoidal functor (the last one even to a symmetric monoidal functor, by construction), which, in turn induces a functor between the associated $\infty$-categories of commutative algebras.

Concretely, the sign rule in the homotopy groups comes from the fact that given two spheres $S^i$ and $S^j$ in spectra, the symmetry map
$$\xymatrix{
{ S^{i+j} \simeq S^i \otimes S^j } \ar[r] &
{ S^j \otimes S^i \simeq S^{j+i} } \cr
}$$
is homotopic to $(-1)^{ij}$ times the identity map of $S^{i+j}$. 
\end{example}

Motivated by this example, we will consider the grading on a Dirac ring $A$ to be homological. Given a homogeneous element $x$, we will refer to the half-integer given by the half of the cohomological degree as the spin, so that
$$\textstyle{ \operatorname{spin}(x) = - \frac{1}{2} \deg(x). }$$
Associated with $A \in \operatorname{CAlg}(\operatorname{Ab})$, we have the category $\operatorname{Mod}_A(\operatorname{Ab})$ of $A$-modules in the symmetric monoidal category $\operatorname{Ab}$. An $A$-module $M = (M_k)_{k \in \mathbb{Z}}$ is the same as a graded module over the graded ring $A$.

\begin{definition}
\label{definition:Serre_twist}
Let $A$ be a Dirac ring, and let $M = (M_k)_{k \in \mathbb{Z}}$ be an $A$-module. For every $s \in \frac{1}{2}\mathbb{Z}$, the $A$-module $M(s)$ defined by
$$M(s)_k = M_{k+2s}$$
is called the spin $s$ Serre twist of $M$.
\end{definition}

The graded $A$-modules $M(s)$ given by the Serre twists are in general pairwise non-isomorphic: There is no process in Dirac geometry that can undo spin. We also note that, for integer spin, \cref{definition:Serre_twist} recovers the usual definition of Serre twist of graded modules over graded commutative rings.

\begin{example}
\label{example:Serre_twist}
If $E \in \operatorname{CAlg}(\operatorname{Sp})$  and $M \in \operatorname{Mod}_E(\operatorname{Sp})$, then, as in \cref{example:homotopy_groups_of_einfty_ring_form_a_dirac_ring}, the graded abelian group $\pi_*(M)$ acquires a canonical structure of module over the Dirac ring $\pi_*(E)$. Moreover, in terms of Serre twist, we have
$$\pi_*(M)(s) \simeq \pi_*(\Sigma^{-2s}M).$$
\end{example}

We now observe that, although the underlying rings of Dirac rings generally are non-commutative, there is a good theory of ideals. 

\begin{lemma}
\label{lemma:left_and_right_ideals_coincide}
Let $A$ be a Dirac ring. If $I \subset A$ is a sub-graded abelian group, then the following conditions are equivalent: 
\begin{enumerate}
    \item[{\rm (1)}] $I \subset A$ is a left sub-$A$-module.
    \item[{\rm (2)}] $I \subset A$ is a right sub-$A$-module.
\end{enumerate}
\end{lemma}

\begin{proof}
If $I$ satisfies~(1), then for all homogeneous elements $x \in I$ and $a \in A$, we have $a \cdot x \in I$. But, in this case, we also have
$$x \cdot a = (-1)^{\mathrm{deg}(x) \mathrm{deg}(a)} a \cdot x \in I,$$
since $I$ is a sub-graded abelian group, and hence, is closed under forming additive inverses. So $I$ satisfies~(2). The opposite implication is proved analogously.
\end{proof}

\begin{definition}
\label{definition:graded_ideal}
A sub-graded abelian group $I \subset A$ of a Dirac ring $A$ is a graded ideal if it satisfies the equivalent conditions of \cref{lemma:left_and_right_ideals_coincide}.
\end{definition}

As in ordinary commutative algebra, two particular important classes of graded ideals are the graded prime ideals and the maximal graded ideals.

\begin{definition}
\label{def:localdiracring}
Let $A$ be a Dirac ring, and let $\mathfrak{a} \subset A$ be a graded ideal.
\begin{enumerate}
\item[(1)]The ideal $\mathfrak{a} \subset A$ is a maximal graded ideal if it is a proper graded ideal and if it is maximal with this property.
\item[(2)]The ideal $\mathfrak{a} \subset A$ is a graded prime ideal if it is a proper ideal and if for all homogeneous elements $f,g \in A$, $fg \in \mathfrak{a}$ implies that $f \in \mathfrak{a}$ or $g \in \mathfrak{a}$ or both.
\end{enumerate}
The Dirac ring $A$ is local if it has a unique maximal graded ideal.
\end{definition}

\begin{warning}
We note that a maximal graded ideal of a Dirac ring is not necessarily a maximal ideal of the underlying ring. 
\end{warning}

\begin{lemma}\label{lem:primeideal}
Let $A$ be a Dirac ring, and let $\mathfrak{a} \subset A$ be a graded ideal.
\begin{enumerate}
\item[{\rm (1)}]The graded ideal $\mathfrak{a} \subset A$ is a graded prime ideal if and only if the underlying ideal in the underlying ring $\smash{\widetilde{\mathfrak{a}}} \subset \smash{\widetilde{A}}$ is a prime ideal.
\item[{\rm (2)}]The ideal $\mathfrak{a} \subset A$ is radical if and only if for every homogeneous element $f \in A$ and integer $n \geq 0$, $f^n \in \mathfrak{a}$ implies that $f \in \mathfrak{a}$.
\end{enumerate}
\end{lemma}

\begin{proof}
Both statements are proved by using the fact that an element $f \in A$ belongs to $\mathfrak{a}$ if and only if all of its homogeneous parts $f_d \in A$ belong to $\mathfrak{a}$.
\end{proof}

\begin{corollary}
\label{cor:primeideal}A maximal graded ideal in a Dirac ring is a prime ideal.
\end{corollary}

\begin{proof}
Let $A$ be a Dirac ring, and let $\mathfrak{m} \subset A$ be a maximal graded ideal. To prove that $\mathfrak{m}$ is a prime ideal, it suffices by \cref{lem:primeideal} to show that if $f,g \in A$ are homogeneous elements with $fg \in \mathfrak{m}$, then either $f \in \mathfrak{m}$ or $g \in \mathfrak{m}$ or both. So let $f,g \in A$ be homogeneous elements with $fg \in \mathfrak{m}$. If $f \in \mathfrak{m}$, then we are done, and if $f \notin \mathfrak{m}$, then the maximality of $\mathfrak{m}$ implies that $A = (f) + \mathfrak{m}$. But then we can write $1 = af+r$ with $r \in \mathfrak{m}$, so $g = afg +rg \in \mathfrak{m}$.
\end{proof}

\begin{lemma}
\label{lem:zorn}
Let $A$ be a Dirac ring. If $\mathfrak{a} \subset A$ is a proper graded ideal, then there exists a maximal graded ideal $\mathfrak{m} \subset A$ such that $\mathfrak{a} \subset \mathfrak{m}$.
\end{lemma}

\begin{proof}
We wish to show that the set $S$ of proper graded ideals $\mathfrak{b} \subset A$ with $\mathfrak{a} \subset \mathfrak{b}$, partially ordered under inclusion, has a maximal element. By Zorn's lemma, it suffices to show that $S$ is non-empty and that every totally ordered subset $T \subset S$ admits an upper bound. First, the set $S$ is non-empty, because $\mathfrak{a} \in S$. Second, if $T \subset S$ is a totally ordered subset, then $\mathfrak{c} = \bigcup_{\mathfrak{b} \in T}\mathfrak{b} \subset A$ is a graded ideal, and it remains to prove that $\mathfrak{c} \neq A$. If not, then $1 \in \mathfrak{c}$. But then $1 \in \mathfrak{b}$ for some $\mathfrak{b} \in T$, which is a contradiction. This completes the proof.
\end{proof}

\begin{definition}
\label{def:diracfield}A Dirac ring $k$ is a Dirac field if any of the following equivalent conditions are satisfied:
\begin{enumerate}
\item The Dirac ring $k$ is not the zero Dirac ring, and every nonzero homogeneous element of $k$ is invertible.
\item The zero ideal $(0) \subset k$ is prime and it is the only graded prime ideal of $k$.
\end{enumerate}
\end{definition}

If $\mathfrak{m} \subset A$ is a maximal graded ideal, then $k(\mathfrak{m}) = A/\mathfrak{m}$ is a Dirac field. A Dirac field is not necessarily a field. We have the following structure theorem.

\begin{proposition}
\label{prop:classification_of_dirac_fields}
Let $k$ be a Dirac field. The subring $k_0 \subset k$ of homogeneous elements of degree $0$ is a field, and if $k_0 \neq k$, then there exists $t \in k$ homogeneous with $\deg(t) \neq 0$ such that $k = k_0[t^{\pm1}]$. If $\operatorname{char}(k_0) \neq 2$, then $\deg(t)$ is even.
\end{proposition}

\begin{proof}We claim that every nonzero homogeneous element $a \in k$ is a unit. If not, then $(a) \subset k$ is a proper graded ideal, so by \cref{lem:zorn}, there exists a maximal graded ideal $(a) \subset \mathfrak{m} \subset k$, and by \cref{cor:primeideal}, $\mathfrak{m} \subset k$ is a nonzero prime ideal. This proves the claim. In particular, the subring $k_0 \subset k$ is a field. So we suppose that $k_0 \neq k$ and let $t \in k$ be a nonzero homogeneous element with $\deg(t) > 0$ and minimal. We have $k_0[t^{\pm1}] \subset k$ and wish to prove that $k_0[t^{\pm1}] = k$. If not, then there exists a homogeneous element $u \in k$ such that $u \notin k_0[t^{\pm1}]$. We write
$$\deg(u) = d \deg(t) + r$$
with $d$ and $0 \leq r < \deg(t)$ integers. If $r = 0$, then $t^{-d}u \in k_0$, which contradicts the fact that $u \notin k_0[t^{\pm1}]$, and if $r \neq 0$, then $0 < \deg(t^{-d}u) < \deg(t)$, which contradicts the fact that $t$ was chosen with $\deg(t) > 0$ minimal. So $k = k_0[t^{\pm1}]$ as desired. Finally, if $\operatorname{char}(k_0) \neq 2$, then homogeneous elements of odd degree square to zero, so they are not units, and therefore, we necessarily have that $\deg(t)$ is even.
\end{proof}

\begin{remark}
We observe that, by \cref{prop:classification_of_dirac_fields}, the underlying ring of every Dirac field is a commutative ring. In general, however, it is not a field.
\end{remark}

\subsection{Localization}

As remarked, the underlying ring of a Dirac ring is typically not commutative. However, the Ore condition is satisfied, and therefore, localizations are well-behaved and admit a calculus of fractions. We discuss this is some detail following Quillen's account in~\cite[Appendix~Q]{friedlandermazur}.

\begin{definition}\label{def:localization}
A localization of a Dirac ring $A$ with respect to a multiplicative subset $S \subset A$ of homogeneous elements is a homomorphism of Dirac rings
$$\xymatrix{
{ A } \ar[r]^-{\gamma} &
{ S^{-1}A } \cr
}$$
that takes every $s \in S$ to a unit $\gamma(s) \in S^{-1}A$ and which is initial with respect to this property.
\end{definition} 
A localization exists and is unique, up to unique isomorphism under $A$, by general arguments, and we now show that it can be constructed by left fractions. So let $S\,\backslash S$ be the category with objects the elements of $S$, with morphisms from $s_1$ to $s_2$ the elements $t \in S$ such that $ts_1 = s_2$, and with composition of morphisms given by multiplication in $S$. We will use that $S\,\backslash S$ is a filtered category. This means that:
\begin{enumerate}
\item[(i)]For all $s_1,s_2 \in S$, there exists $t_1,t_2 \in S$ such that $t_1s_1 = t_2s_2$.
\item[(ii)]For all $s,s_1,s_2 \in S$ such that $s_1s = s_2s$, there exists $t \in S$ such that $ts_1 = ts_2$.
\end{enumerate}
There is a functor from $S\,\backslash S$ to the category graded right $A$-modules that takes the object $s$ to the Serre twist $A(-\operatorname{spin}(s))$ and the morphism $t \colon s_1 \to s_2$ to the map $l_t \colon A(-\operatorname{spin}(s_1)) \to A(-\operatorname{spin}(s_2))$ given by left multiplication by $t$. We define
$$\textstyle{ B = \varinjlim_{s \in S\,\backslash S} A(-\operatorname{spin}(s)) }$$
to be the colimit of this functor. The general description of a filtered colimit of sets~\cite[Expos\'{e}~I, Lemme~2.8.1]{SGA4I} identifies $B$ with the set of fractions $s^{-1}a$, where $s_1^{-1}a_1 = s_2^{-1}a_2$ if and only if there exists $t_1,t_2 \in S$ with $t_1a_1 = t_2a_2$ and $t_1s_1 = t_2s_2$. If $s \in S$ and $a \in A$ are homogeneous of degrees $d$ and $e$, then $s^{-1}a$ is homogeneous of degree $e-d$. We will use that for all $s \in S$, the map $r_s \colon B \to B(-\operatorname{spin}(s))$ given by right multiplication by $s$ is an isomorphism. This means that:
\begin{enumerate}
\item[(iii)]Given $a \in A$ and $s \in S$ with $as = 0$, there exists $t \in S$ with $ta = 0$.
\item[(iv)]Given $a \in A$ and $s \in S$, there exists $b \in A$ and $t \in S$ such that $ta = bs$.
\end{enumerate}
The assumptions~(i)--(iv) are trivially verified.

\begin{proposition}\label{prop:calculusoffractions}Let $A$ be a Dirac ring and let $\gamma \colon A \to S^{-1}A$ be the localization with respect to a multiplicative subset $S \subset A$ that consists of homogeneous elements. In this situation, the map
$$\xymatrix{
{ B = \varinjlim_{s \in S\,\backslash S} A(-\operatorname{spin}(s)) } \ar[r]^-{u} &
{ S^{-1}A } \cr
}$$
that to $s^{-1}a$ assigns $\gamma(s)^{-1}\gamma(a)$ is an isomorphism of graded right $A$-modules. In particular, the localization $\gamma \colon A \to S^{-1}A$ is a filtered colimit of free $A$-modules. 
\end{proposition}

\begin{proof}The internal mapping space $\operatorname{End}(B)$ in the category of $\mathbb{Z}$-graded abelian groups has a canonical structure of $\mathbb{Z}$-graded ring, and moreover, there is a map of $\mathbb{Z}$-graded rings $r \colon A^{\operatorname{op}}\to \operatorname{End}(B)$ that to a homogeneous element $a \in A$ assigns the map $r_a \colon B \to B(-\operatorname{spin}(a))$ given by right multiplication by $a$. Since~(iii)--(iv) hold, the map $r$ extends uniquely to a map of $\mathbb{Z}$-graded rings $\tilde{r} \colon S^{-1}A^{\operatorname{op}} \to \operatorname{End}(B)$. Thus, the structure of graded right $A$-module on $B$ extends uniquely to a structure of graded right $S^{-1}A$-module. Moreover, the map $u$ in question is $S^{-1}A$-linear, since every homogeneous element in $S^{-1}A$ is a finite products of elements of the form $\gamma(s)^{-1}$ and $\gamma(a)$ with $s \in S$ and $a \in A$ homogeneous. Finally, since $B$ is generated as a right $S^{-1}A$-module by the fraction $(1)^{-1}1$, and since $u$ maps this generator to the identity element in $S^{-1}A$, we conclude that $u$ is an isomorphism.
\end{proof}

\begin{remark}
As a consequence of \cref{prop:calculusoffractions}, the localization $S^{-1}A$ is a flat $A$-module. We discuss flatness in more depth in \cref{subsection:commutative_algebra_flatness}.
\end{remark}

\subsection{The Dirac--Zariski space}

We wish to interpret every Dirac ring $A$ as the ring of functions on a topological space $|X| = \lvert\operatorname{Spec}(A)\rvert$. We call $|X|$ the Dirac--Zariski space of $A$, and as already mentioned, its points are defined to be the graded prime ideals $\mathfrak{p} \subset A$. Before we proceed, we introduce some familiar notation which emphasizes the interpretation of $A$ as a Dirac ring of functions on $|X|$. If $x \in |X|$ corresponds to the graded prime ideal $\mathfrak{p} \subset A$, then we write
$$\mathcal{O}_X(|X|) \to \mathcal{O}_{X,x}$$
for the localization $A \to A_{\mathfrak{p}}$ with respect to the multiplicative subset $S \subset A$ that consists of all homogeneous elements $f \in A$ that are not contained in $\mathfrak{p}$. We call $\mathcal{O}_{X,x}$ the local Dirac ring at $x \in |X|$, and we write
$$\mathfrak{m}_x \subset \mathcal{O}_{X,x}$$
for its unique maximal graded ideal. It is the graded ideal $\mathfrak{p}A_{\mathfrak{p}} \subset A_{\mathfrak{p}}$. The quotient
$$k(x) = \mathcal{O}_{X,x}/\mathfrak{m}_x$$
is a Dirac field, which we call the residue Dirac field at $x \in |X|$. We denote the image of $f \in \mathcal{O}_X(|X|)$ by the composite map
$$\xymatrix{
{ \mathcal{O}_X(|X|) } \ar[r] &
{ \mathcal{O}_{X,x} } \ar[r] &
{ k(x) } \cr
}$$
by $f(x)$ and call it the value of $f$ at $x$. Finally, we write
$$\xymatrix{
{ \mathcal{O}_X(|X|) } \ar[r] &
{ \mathcal{O}_Z(|Z|) } \cr
}$$
for the projection $A \to A/\mathfrak{p}$. Here $Z \subset X$ is the corresponding closed sub-Dirac scheme, which we define below. The closed subspace $|Z| = V(\mathfrak{p}) \subset |X|$ is the closure of $\{x\} \subset |X|$. These maps constitute a commutative diagram
$$\xymatrix{
{ \mathcal{O}_X(|X|) } \ar[r] \ar[d] &
{ \mathcal{O}_{X,x} } \ar[d] \cr
{ \mathcal{O}_Z(|Z|) } \ar[r] &
{ k(x) } \cr
}$$
with the lower horizontal map the localization at $(0) = (\mathfrak{p}) \subset A/\mathfrak{p}$.

The topology of the space $|X| = \lvert\operatorname{Spec}(A)\rvert$ is the Zariski topology, for which a basis consists of the family of distinguished open subsets
$$|X_f| = \{x \in |X| \mid f(x) \neq 0 \} \subset |X|$$
with $f \in \mathcal{O}_X(|X|)$ homogeneous.

\begin{definition}\label{def:zariskispace}Let $A$ be a Dirac ring. Its Dirac--Zariski space is the topological space $\lvert\operatorname{Spec}(A)\rvert$ given by the set of graded prime ideals equipped with the Zariski topology.
\end{definition}

We enumerate a number of properties of the Dirac--Zariski space, which it shares with the Zariski space of a commutative ring.

\begin{proposition}\label{prop:specquasicompact}Let $A$ be a Dirac ring, let $|X| = \lvert\operatorname{Spec}(A)\rvert$, and let $(f_i)_{i \in I}$ be a family of homogeneous elements of $A$. The following are equivalent:
\begin{enumerate}
\item[{\rm (1)}]The family $(|X_{f_i}|)_{i \in I}$ of distinguished open subsets of $|X|$ covers $|X|$.
\item[{\rm (2)}]The family $(f_i)_{i \in I}$ of homogeneous elements of $A$ generates the unit ideal.    
\end{enumerate}
In particular, the topological space $|X| = \lvert\operatorname{Spec}(A)\rvert$ is quasi-compact.
\end{proposition}

\begin{proof}By definition,~(1) holds if and only if for all $x \in |X|$, there exists some $i \in I$ such that $f_i(x) \neq 0$. This, in turn, is equivalent to the statement that for every graded prime ideal $\mathfrak{p} \subset A$, there exists $i \in I$ such that $f_i \notin \mathfrak{p}$, and this holds if and only if no graded prime ideal contains the (graded) ideal generated by $(f_i)_{i \in I}$. But \cref{lem:zorn} and \cref{cor:primeideal} show that the latter statement is equivalent to~(2).

To prove that $|X|$ is quasi-compact, suppose that a family $(U_i)_{i \in I}$ of open subsets of $|X|$ covers $|X|$. For every $x \in |X|$, we first choose $i \in I$ such that $x \in U_i$ and then choose a homogeneous element $f_x \in A$ such that $x \in |X_{f_x}| \subset U_i$. It will suffice to show that a finite subfamily of the family $(|X_{f_x}|)_{x \in |X|}$ of distinguished open subsets cover $|X|$. Now, since $(|X_{f_x}|)_{x \in |X|}$ covers $|X|$, the family $(f_x)_{x \in |X|}$ generates the unit ideal, so we can write $1 = a_1f_{x_1} + \dots + a_nf_{x_n}$ for some $x_1,\dots,x_n \in |X|$. But then the subfamily $(f_{x_1},\dots,f_{x_n})$ of $(f_x)_{x \in |X|}$ generates the unit ideal, and therefore, the subfamily $(|X_{f_{x_1}}|, \dots,|X_{f_{x_n}}|)$ of $(|X_{f_x}|)_{x \in |X|}$ covers $|X|$.
\end{proof}

We will also associate to a map of Dirac rings $\phi \colon A \to B$ a map of locally Dirac ringed spaces $p \colon (|Y|,\mathcal{O}_Y) \to (|X|,\mathcal{O}_X)$ in the opposite direction. We define the map of underlying sets to be the map that to a point $y \in |Y|$ corresponding to the graded prime ideal $\mathfrak{q} \subset B$ assigns the point $x \in |X|$ that corresponds to the graded prime ideal $\mathfrak{p} = \phi^{-1}(\mathfrak{q}) \subset A$. If $f \in \mathcal{O}_X(|X|)$ is homogeneous, then
$$p^{-1}(|X_f|) = |Y_{\phi(f)}|,$$
so in particular the map $p \colon |Y| \to |X|$ is continuous.

\begin{proposition}\label{prop:embedding}The map of $p \colon |Y| \to |X|$ of Dirac--Zariski spaces induced by a map $\phi \colon A \to B$ of Dirac rings is an embedding in the following cases:
\begin{enumerate}
\item[{\rm (1)}]If $\mathfrak{a} \subset A$ is a graded ideal and $\phi \colon A \to A/\mathfrak{a}$ the canonical projection, then $p \colon |Y| \to |X|$ is an embedding with closed image
$$\textstyle{ V(\mathfrak{a}) = \{x \in |X| \mid \text{$f(x) = 0$ for all $f \in \mathfrak{a}$ homogeneous}\}. }$$
\item[{\rm (2)}]If $S \subset A$ is a multiplicative subset consisting of homogeneous elements and $\phi \colon A \to S^{-1}A$ the localization, then $p \colon |Y| \to |X|$ is an embedding with image
$$\textstyle{ p(|Y|) = \{x \in |X| \mid \text{$f(x) \neq 0$ for all $f \in S$}\}, }$$
which in general it is neither open nor closed.
\end{enumerate}
\end{proposition}

\begin{proof}This is proved as in~\cite[Tag~00E3, Tag~00E5]{stacks-project}.
\end{proof}

\begin{example}\label{ex:embedding}Let $A$ be a Dirac ring, and let $\phi_f \colon A \to A_f$ be the localization with respect to the multiplicative subset $S \subset A$ consisting of all powers of a homogeneous element $f \in A$. In this situtation, \cref{prop:embedding}~(2) shows that the induced map of Dirac--Zariski spaces $j \colon |Y| \to |X|$ is an embedding with image the distinguised open subset $|X_f| \subset |X|$.
\end{example}

\begin{lemma}\label{lemma:radical_ideals_as_intersection_of_primes}
Let $A$ be a Dirac ring. The intersection of a family of graded prime ideals is a graded radical ideal, and every graded radical ideal is the intersection of the family of graded prime ideals that contain it.
\end{lemma}

\begin{proof}The first statement follows from \cref{lem:primeideal} and from the fact that every intersection of prime ideals is a radical ideal. To prove the second statement, we must show that if $f \in A$ is homogeneous and $f \notin \mathfrak{a}$, then there exists a graded prime ideal $\mathfrak{p} \subset A$ such that $\mathfrak{a} \subset \mathfrak{p}$ and $f \notin \mathfrak{p}$. So we let $\mathfrak{q} \subset A$ be a graded ideal that is maximal among the homogeneous ideals that contain $\mathfrak{a}$ and that do not contain any power of $f$. We claim $\mathfrak{q}$ is a prime ideal, which will complete the proof. If not, then we conclude from \cref{lem:primeideal} that there exists homogeneous elements $g,h \in A$ with $g,h \notin \mathfrak{q}$ but with $g \cdot h \in \mathfrak{q}$. By the maximality of $\mathfrak{q}$, both $\mathfrak{q} + (g)$ and $\mathfrak{q} + (h)$ contain some power of $f$. Hence, for some $m,n \geq 0$, we can write
$$\begin{aligned}
f^m & = ag + r \cr
f^n & = bh + s \cr
\end{aligned}$$
with $a,b \in A$ and $r,s \in \mathfrak{q}$. But multiplying these two equations, we conclude that $f^{m+n} \in \mathfrak{q}$, which is a contradiction.
\end{proof}

\begin{corollary}\label{cor:closedsets}Let $A$ be a Dirac ring with Dirac--Zariski space $|X|$.
\begin{enumerate}
\item[{\rm (1)}]Every closed subset $Z \subset |X|$ is of the form $Z = V(\mathfrak{a})$ for a unique graded radical ideal $\mathfrak{a} \subset A$.
\item[{\rm (2)}]Every irreducible closed subset $Z \subset |X|$ is of the form $Z = V(\mathfrak{p})$ for a unique graded prime ideal $\mathfrak{p} \subset A$.
\end{enumerate}
\end{corollary}

\begin{proof}To prove~(1), let $\mathfrak{a} \subset A$ be the ideal generated by the family of homogeneous elements $f \in A$ such that $f(x) = 0$ for all $x \in Z$. It is clear that $\mathfrak{a} \subset A$ is a graded ideal and that it is maximal with the property that $V(\mathfrak{a}) = Z$. We claim that $\mathfrak{a} \subset A$ is a radical ideal. Let us write $\mathfrak{p}_x \subset A$ for the graded prime ideal corresponding to $x \in |X|$. By definition, we have $\mathfrak{a} = \bigcap_{x \in Z} \mathfrak{p}_x$, and hence, the claim follows from \cref{lemma:radical_ideals_as_intersection_of_primes}. Finally, if $Z = V(\mathfrak{b})$ for some graded ideal $\mathfrak{b} \subset A$, then $\mathfrak{a} = \sqrt{\mathfrak{b}}$.

To prove~(2), we write $Z = V(\mathfrak{a})$ with $\mathfrak{a} \subset A$ a graded radical ideal. We claim that $\mathfrak{a}$ is a prime ideal. If not, then by \cref{lem:primeideal}, there exist $f,g \in A$ homogeneous with $fg \in \mathfrak{a}$ but with $f \notin \mathfrak{a}$ and $g \notin \mathfrak{a}$. Thus, we obtain the decomposition
$$Z = V(\mathfrak{a}) = V(\mathfrak{a} + (f)) \cup V(\mathfrak{a} + (g))$$
of $Z \subset |X|$ as a union of two proper closed subsets, contradicting the assumption that the closed subset $Z \subset |X|$ is irreducible. This prove the claim.
\end{proof}

We recall from~\cite[08YG]{stacks-project} that a topological space $|X|$ is defined to be spectral if it is quasi-compact and sober, if the intersection of two quasi-compact open subsets is quasi-compact, and if the family of quasi-compact open subsets is a basis for the topology. A topological space $|X|$ is defined to be sober if every irreducible closed subset $Z \subset |X|$ has a unique generic point. A map $p \colon |Y| \to |X|$ between spectral spaces is defined to be spectral if for every quasi-compact open subset $U \subset |X|$, the inverse image $p^{-1}(U) \subset |Y|$ is a quasi-compact open subset.

\begin{proposition}\label{prop:spectralspace}The Dirac--Zariski space $|X|$ of a Dirac ring $A$ is a spectral space, and the map of Dirac--Zariski spaces $p \colon |Y| \to |X|$ induced by a map of Dirac rings $\phi \colon A \to B$ is a spectral map.
\end{proposition}

\begin{proof}We first show that $|X| = \lvert\operatorname{Spec}(A)\rvert$ is a spectral space. By \cref{prop:specquasicompact} and \cref{prop:embedding}, the distinguished open subsets $|X_f| \subset |X|$ are quasi-compact, and since the family of distinguished open subsets is a basis for the Zariski topology on $|X|$, a subset $U \subset |X|$ is quasi-compact and open if and only if it is a finite union of distinguished open subsets. Moreover, if $|X_f|,|X_g| \subset |X|$ are distinguished open subsets, then so is their intersection
$|X_f| \cap |X_g| = |X_{fg}| \subset |X|$. It remains to prove that $|X|$ is sober. So let $Z \subset |X|$ be an irreducible closed subset. We use \cref{cor:closedsets} to write $Z = V(\mathfrak{p})$ for a unique graded prime ideal $\mathfrak{p} \subset A$ and let $x \in |X|$ be the corresponding point. If $x \in V(\mathfrak{a})$, then $\mathfrak{a} \subset \mathfrak{p}$, so $Z = V(\mathfrak{p}) \subset V(\mathfrak{a})$, which shows that $x \in Z$ is a generic point. It is unique, since if also $x' \in Z$ is a generic point corresponding to a graded prime ideal $\mathfrak{p}' \subset A$, then we have $\mathfrak{p} \subset \mathfrak{p}'$ and $\mathfrak{p}' \subset \mathfrak{p}$, which shows that $x = x'$. Finally, we have $p^{-1}(|X_f| = |Y_{\phi(f)}|$, which shows that $p \colon |Y| \to |X|$ is spectral.
\end{proof}

\begin{remark}
\label{remark:stone_duality_and_dirac_zariski_spectrum}
Let $\operatorname{Top}^{\operatorname{spec}} \subset \operatorname{Top}$ be the non-full subcategory spanned by the spectral spaces and the spectral maps between them, and let $\operatorname{Lat}^{\operatorname{dist}} \subset \operatorname{Lat}$ be the full subcategory of the category $\operatorname{Lat}$ of lattices and lattice homomorphisms spanned by the distributive lattices.\footnote{\,We assume, as does~\cite{johnstone1982stone}, that lattices are bounded in the sense that they have both a maximal element $1$ and a minimal element $0$.} By Stone duality, the functor
$$\xymatrix{
{ (\operatorname{Top}^{\operatorname{spec}})^{\operatorname{op}} } \ar[r] &
{ \operatorname{Lat}^{\operatorname{dist}} } \cr
}$$
that to a spectral space assigns its lattice of quasi-compact open subsets is an equivalence. This was proved in~\cite{stone1938topological}, but see also~\cite[Corollary~I.3.4]{johnstone1982stone}. Under this equivalence, the Dirac--Zariski space $|X| = \lvert\operatorname{Spec}(A)\rvert$ of a Dirac ring $A$ is mapped to the distributive lattice $\operatorname{Rad}^{\operatorname{fg}}(A)$ of finitely generated radical graded ideals of $A$. Indeed, by \cref{cor:closedsets}, closed subsets $Z \subset |X|$ correspond to radical graded ideals $\mathfrak{a} \subset A$, and one readily checks that $Z \subset |X|$ has quasi-compact open complement if and only if $\mathfrak{a} \subset A$ is finitely generated.

We may characterize the distributive lattice $\operatorname{Rad}^{\operatorname{fg}}(A)$ of finitely generated radical graded ideals $\mathfrak{a} \subset A$ by a universal property as follows. We define a map
$$\xymatrix{
{ \coprod_{d \in \mathbb{Z}} A_d } \ar[r]^-{d} &
{ D } \cr
}$$
from the set of homogeneous elements in $A$ to a distributive lattice $D$  to be a support function on $A$ if it has the following properties:
\begin{enumerate}
\item $d(1)$ is the maximal element of $D$.
\item $d(0)$ is the minimal element of $D$. 
\item $d(fg) = d(f) \wedge d(g)$.
\item $d(f+g) \leq d(f) \vee d(g)$.
\end{enumerate}
It follows from these axioms that if $f \in \sqrt{(g)}$, then $d(f) \leq d(g)$, so, in particular, if $\sqrt{(f)} = \sqrt{(g)}$, then $d(f) = d(g)$. This, in turn, implies that the map
$$\xymatrix{
{ \coprod_{d \in \mathbb{Z}} A_d } \ar[r]^-{d^{\operatorname{fg}}} &
{ \operatorname{Rad}^{\operatorname{fg}}(A) } \cr
}$$
defined by $d^{\operatorname{fg}}(f) = \sqrt{(f)}$ is the initial support function on $A$.

We may use Stone duality to translate this to a lattice-theoretic characterization of the Dirac--Zariski space $|X| = \lvert\operatorname{Spec}(A)\rvert$ as follows. Let $Y$ be a topological space, and let $Y_{\operatorname{Zar}}$ be its lattice of open subsets. We define a map
$$\xymatrix{
{ \coprod_{d \in \mathbb{Z}}A_d } \ar[r]^-{\operatorname{supp}} &
{ Y_{\operatorname{Zar}} } \cr
}$$
to be a space-valued support function on $A$ if it has the following properties:
\begin{enumerate}
\item $\operatorname{supp}(1) = Y$.
\item $\operatorname{supp}(0) = \emptyset$.
\item $\operatorname{supp}(fg) = \operatorname{supp}(f) \cap \operatorname{supp}(g)$.
\item $\operatorname{supp}(f + g) \subseteq \operatorname{supp}(f) \cup \operatorname{supp}(g)$.
\end{enumerate}
The map $D \colon \coprod_{d \in \mathbb{Z}}A_d \to |X|_{\operatorname{Zar}}$ that to a homogeneous element $f \in A$ assigns the distinguised open subset $D(f) = |X_f| \subset |X|$ is a space-valued support function on $A$ that witnesses the Dirac--Zariski space as the final topological space equipped with a space-valued support function on $A$. 
\end{remark}

\subsection{The category of Dirac schemes}

We define the structure sheaf $\mathcal{O}_X$ on the Dirac--Zariski space $|X|$ of a Dirac ring $A$, and define the category of Dirac schemes.

\begin{theorem}
\label{thm:structuresheaf}
{\rm(1)}~Let $A$ be a Dirac ring, and let $|X| = \lvert\operatorname{Spec}(A)\rvert$. Up to unique isomorphism, there is a unique pair $(\mathcal{O}_X,\epsilon_A \colon A \to \mathcal{O}_X(|X|))$ of a sheaf of Dirac rings on $|X|$ and map of Dirac rings such that for every homogeneous element $f \in A$, the composite map
$$\xymatrix{
{ A } \ar[r]^-{\epsilon_A} &
{ \mathcal{O}_X(|X|) } \ar[r] &
{ \mathcal{O}_X(|X_f|) } \cr
}$$
is a localization with respect to $S = \{1,f,f^2,\dots\} \subset A$.

\noindent~{\rm(2)}~Let $\phi \colon A \to B$ be a map of Dirac rings, let $p \colon |Y| \to |X|$ be the associated map of Dirac--Zariski spaces, and let $\phi \colon \mathcal{O}_X \to p_*(\mathcal{O}_Y)$ be the unique map of sheaves of Dirac rings on $|X|$ that makes the following diagram commute.
$$\xymatrix@C=11mm{
{ A } \ar[r]^-{\epsilon_A} \ar[d]^-{\phi} &
{ \mathcal{O}_X(|X|) } \ar[d]^-{\phi_X} \cr
{ B } \ar[r]^-{\epsilon_B} &
{ \mathcal{O}_Y(|Y|) } \cr
}$$
The pair $(p,\phi) \colon (|Y|,\mathcal{O}_Y) \to (|X|,\mathcal{O}_X)$ is a map of locally Dirac ringed spaces.
\end{theorem}

\begin{proof}To prove~(1), we apply~\cite[Proposition~1.1.4.4]{luriesag}. So let $X_{\operatorname{Zar}}$ be the category of open subsets of $|X|$ and with a single map from $U$ to $V$, if $U \subset V$, and with no maps from $U$ to $V$, otherwise, and let $u \colon D_{\operatorname{Zar}} \to X_{\operatorname{Zar}}$ be the inclusion of the full subcategory spanned by the distinguished open subsets. The family of distinguished open subsets satisfy the conditions~(i)--(iii) of loc.~cit., so it suffices to prove that, up to isomorphism, there is a unique pair $(\mathcal{O}_X \circ u,\epsilon_A \colon A \to (\mathcal{O}_X \circ u)(|X|)$ of a sheaf of Dirac rings on $D_{\operatorname{Zar}}$ and a map of Dirac rings, such that for every homogeneous element $f \in A$, the composite map
$$\xymatrix{
{ A } \ar[r]^-{\epsilon_A} &
{ (\mathcal{O}_X \circ u)(|X|) } \ar[r] &
{ (\mathcal{O}_X \circ u)(|X_f|) } \cr
}$$
is a localization with respect to $S = \{1,f,f^2,\dots\} \subset A$. The uniqueness, up to unique isomorphism, of a solution to this problem is clear, so it remains to prove existence. We first prove~(a) that the desired pair with $\mathcal{O}_X \circ u$ a presheaf of Dirac rings on $D_{\operatorname{Zar}}$ with the desired properties exists, and then prove~(b) that the presheaf $\mathcal{O}_X \circ u$ is a sheaf. If $|X_f| \subset |X_g|$, then the localization $\phi_f \colon A \to A_f$ maps $g$ to a unit. For $g(x) \neq 0$ for all $x \in |X_f|$, so $\phi_f(g) \in A_f$ is not contained in any graded prime ideal, so \cref{lem:zorn} and \cref{cor:primeideal} show that $\phi_f(g) \in A_f$ is a unit. Thus, if $|X_f| \subset |X_g|$, then the localization $\phi_f \colon A \to A_f$ factors uniquely as a composition
$$\xymatrix{
{ A } \ar[r]^-{\phi_g} &
{ A_g } \ar[r]^-{\phi_{f,g}} &
{ A_f, } \cr
}$$
which proves~(a). To prove~(b), we remark that
$$|X_f| \cap |X_g| = |X_{fg}|.$$
Therefore, we must show that if $(|X_{f_i}| \to |X_f|)_{1 \leq i \leq n}$ is a finite jointly surjective family of maps in $D_{\operatorname{Zar}}$, then the diagram
$$\xymatrix{
{ A_f } \ar[r]^-{\alpha} &
{ \prod_{1 \leq i \leq n} A_{f_i} } \ar@<.7ex>[r]^-{\beta} \ar@<-.7ex>[r]_-{\gamma} &
{ \prod_{1 \leq j,k\leq n} A_{f_jf_k} } \cr
}$$
with the maps $\alpha$, $\beta$, and $\gamma$ given by $\operatorname{pr}_i \circ \, \alpha = \phi_{f_i,f}$, $\operatorname{pr}_{j,k} \circ \, \beta = \phi_{f_jf_k,f_j} \circ \operatorname{pr}_j$, and $\operatorname{pr}_{j,k} \circ \, \gamma = \phi_{f_jf_k,f_k} \circ \operatorname{pr}_k$ is an equalizer. But since the map $\alpha$ is faithfully flat, this follows by faithfully flat descent. This proves~(1).

To prove~(2), we first note that $(|X|,\mathcal{O}_X)$ is a locally Dirac ringed space. Indeed, let $i_x \colon \{x\} \to |X|$ is the inclusion of a point. The sheaf $i_x^*(\mathcal{O}_X)$ of Dirac rings on $\{x\}$ is uniquely determined, up to unique isomorphism, by its Dirac ring of global sections $i_x^*(\mathcal{O}_X)(\{x\})$, and that $(|X|,\mathcal{O}_X)$ is a locally Dirac ringed space means that the latter Dirac ring is local for all $x \in |X|$. But
$$i_x^*(\mathcal{O}_X)(\{x\}) \simeq \mathcal{O}_{X,x},$$
which is a local Dirac ring. The map $p \colon |Y| \to |X|$ satisfies $p^{-1}(|X_f|) = |Y_{\phi(f)}|$, so there is a unique map $\phi_{X_f}$ of Dirac rings making the diagram
$$\xymatrix@C=10mm{
{ A } \ar[r]^-{\epsilon_A} \ar[d]^-{\phi} &
{ \mathcal{O}_X(|X|) } \ar[r] \ar[d]^-{\phi_X} &
{ \mathcal{O}_X(|X_f|) } \ar[d]^-{\phi_{X_f}} \cr
{ B } \ar[r]^-{\epsilon_B} &
{ \mathcal{O}_Y(|Y|) } \ar[r] &
{ \mathcal{O}_Y(|Y_{\phi(f)}|) } \cr
}$$
commute for all $f \in A$ homogeneous. This defines the map $\phi \colon \mathcal{O}_X \to p_*(\mathcal{O}_Y)$ of sheaves of Dirac rings on $X$. It remains to show that the induced map of stalks $\phi_y \colon \mathcal{O}_{X,x} \to \mathcal{O}_{Y,y}$ is a local map of local Dirac rings in the sense that $\phi_y(\mathfrak{m}_x) \subset \mathfrak{m}_y$. So let $\mathfrak{p} \subset A$ and $\mathfrak{q} \subset B$ be the graded prime ideals corresponding to $x \in |X|$ and $y \in |Y|$, respectively. By definition of the map $p \colon |Y| \to |X|$, we have $\mathfrak{p} = \phi^{-1}(\mathfrak{q})$, and moreover, the map $\phi_y \colon \mathcal{O}_{X,x} \to \mathcal{O}_{Y,y}$ is canonically identified with the map $\phi_{\mathfrak{q}} \colon A_{\mathfrak{p}} \to B_{\mathfrak{q}}$ induced by $\phi \colon A \to B$, which is indeed local. This proves~(2).
\end{proof}

\begin{definition}\label{def:diracspectrum}Let $A$ be a Dirac ring. Its Dirac spectrum $\operatorname{Spec}(A)$ is the locally Dirac ringed space $X = (|X|,\mathcal{O}_X)$ given by the Dirac--Zariski space of $A$ and the sheaf of Dirac rings $\mathcal{O}_X$ thereon provided by \cref{thm:structuresheaf}.
\end{definition}

\begin{remark}\label{rem:structuresheaf}Let $A$ be a Dirac ring, and let $X = (|X|,\mathcal{O}_X)$ be its Dirac spectrum. One may wonder, where the strange space $|X|$ comes from. The answer was given by Hakim in her thesis~\cite{hakim}. The collection of Dirac ringed topoi is naturally organized into a $2$-category $1$-$\operatorname{Top}_{\operatorname{CAlg}}$; compare~\cite[Definition~1.2.1.1]{luriesag}. If $(\mathcal{X},A)$ is a Dirac ringed topos, then given $U \in \mathcal{X}$ and $f \in A(U)$ homogeneous, we let $U_f \to U$ be the largest subobject such that $A(U) \to A(U_f)$ takes $f$ to a unit. We define $A$ to be local if for every $U \in \mathcal{X}$ and every family $(f_i)_{i \in I}$ of homogeneous elements in $A(U)$ that generates the unit ideal, one has $$\textstyle{ U = \bigcup_{i \in I} U_{f_i}. }$$
Let $(p,\phi) \colon (\mathcal{Y},B) \to (\mathcal{X},A)$ be a map of Dirac ringed topoi with $A$ and $B$ local. We define $\phi \colon A \to p_*(B)$ to be local if for every $U \in \mathcal{X}$ and every homogeneous element $f \in A(U)$, the monomorphism $p^*(U_f) \to p^*(U)_{\phi(f)}$ is an isomorphism. The collection of locally Dirac ringed topoi and local maps between them is organized into a (non-$1$-full) sub-$2$-category $1$-$\operatorname{Top}_{\operatorname{CAlg}}^{\operatorname{loc}}$ of $1$-$\operatorname{Top}_{\operatorname{CAlg}}$. One may prove that the canonical inclusion admits a right adjoint
$$\xymatrix{
{ \text{ $1$-$\operatorname{Top}_{\operatorname{CAlg}}^{\operatorname{loc}}$ } } \ar@<.7ex>[r] &
{ \text{ $1$-$\operatorname{Top}_{\operatorname{CAlg}}$ } } \ar@<.7ex>[l]^-{\operatorname{Spec}} \cr
}$$
and that if $A$ is a Dirac ring in the topos $\operatorname{Set}$, then the adjunct
$$\xymatrix{
{ (\operatorname{Shv}_{\operatorname{Set}}(|X|),\mathcal{O}_X) } \ar[r] &
{ \operatorname{Spec}(\operatorname{Set},A) } \cr
}$$
of the map of Dirac rings $\epsilon_A \colon A \to \mathcal{O}_X(|X|)$ given by \cref{thm:structuresheaf}~(1) is an equivalence of locally Dirac ringed topoi. Thus, the space $|X|$ is the space of points in the underlying topos of the locally Dirac ringed topos $\operatorname{Spec}(\operatorname{Set},A)$.
\end{remark}

\begin{definition}
\label{def:diracscheme}
A Dirac scheme is a locally Dirac ringed space that, locally, is isomorphic to the Dirac-Zariski spectrum of a Dirac ring. The category of Dirac schemes is the full subcategory of the category of locally Dirac ringed spaces spanned by the Dirac schemes.
\end{definition}

\begin{example}
\label{ex:diracscheme}
(1)~If $X = (|X|,\mathcal{O}_X)$ is a Dirac scheme, and if $j \colon |U| \to |X|$ is an open embedding, then $U = (|U|,j^*\mathcal{O}_X)$ is a Dirac scheme. In this situation, we say that $j \colon U \to X$ an open immersion.

\noindent(2)~Let $X = (|X|,\mathcal{O}_X)$ be a Dirac scheme, and let $\phi \colon \mathcal{O}_{X^{\operatorname{ev}}} \to \mathcal{O}_X$ be the inclusion of the subsheaf defined by $\mathcal{O}_{X^{\operatorname{ev}}}(|U|) = \mathcal{O}_X(|U|)^{\operatorname{ev}}$. It is a map of locally Dirac ringed sheaves on $|X|$. Therefore, the pair $X^{\operatorname{ev}} = (|X|,\mathcal{O}_{X^{\operatorname{ev}}})$ is a Dirac scheme and the pair $\eta_X = (\operatorname{id},\phi) \colon X \to X^{\operatorname{ev}}$ is a map of Dirac schemes. 
\end{example}

\begin{proposition}\label{prop:diracscheme}The functor that to a Dirac ring $A$ assigns its Dirac spectrum $X = \operatorname{Spec}(A)$ is right adjoint to the functor that to a locally Dirac ringed space $Y$ assigns the Dirac ring $\mathcal{O}_Y(|Y|)$ of global sections of its structure sheaf.\footnote{\,Here we consider $\operatorname{Spec}$ as a functor from the opposite of the category of Dirac rings to the category of locally Dirac ringed spaces.} Moreover, the counit of the adjunction $\epsilon_A \colon A \to \mathcal{O}_X(|X|)$ is an isomorphism.
\end{proposition}

\begin{proof}We first define the counit and unit of the adjunction. If $A$ is a Dirac ring, then we let $\epsilon_A \colon A \to \mathcal{O}_X(|X|)$ be the map of Dirac rings defined by \cref{thm:structuresheaf}. It is an isomorphism, as stated. If $Y$ is a locally ringed space, if $A = \mathcal{O}_Y(|Y|)$, and if $X = \operatorname{Spec}(A)$, then we let $\eta_Y = (p,\phi) \colon Y \to X$ be the map of locally Dirac ringed spaces defined as follows. The map $p$ assigns to $y \in |Y|$ the element $x \in |X|$ given by the kernel of the composite map
$$\xymatrix{
{ A = \mathcal{O}_Y(|Y|) } \ar[r] &
{ \mathcal{O}_{Y,y} } \ar[r] &
{ k(y). } \cr
}$$
To show that it is continuous, we note that for $f \in A$ homogeneous,
$$y \in p^{-1}(|X_f|) \subset |Y|$$
if and only if its image $f(y) \in k(y)$ by the composite map above is nonzero, or equivalently, if and only if its image $f_y \in \mathcal{O}_{Y,y}$ is a unit. But $\mathcal{O}_{Y,y}$ is the filtered colimit of $\mathcal{O}_Y(|V|)$ with $y \in |V| \subset |Y|$ open, so if $f_y \in \mathcal{O}_{Y,y}$ is invertible, then there exists such a $|V|$ fully contained in
$p^{-1}(|X_f|)$. So $p \colon |Y| \to |X|$ is continuous. The map $\phi \colon \mathcal{O}_X \to p_*\mathcal{O}_Y$ is the unique map of sheaves of Dirac rings on $|X|$ such that $\phi_X = \epsilon_A^{-1}$ and such that for all $f \in A$ homogeneous, the diagram
$$\xymatrix@C=10mm{
{ \mathcal{O}_X(|X|) } \ar[r]^-{\phi_X} \ar[d] &
{ \mathcal{O}_Y(|Y|) } \ar[d] \cr
{ \mathcal{O}_X(|X_f|) } \ar[r]^-{\phi_{X_f}} &
{ \mathcal{O}_Y(p^{-1}(|X_f|)) } \cr
}$$
commutes. Here, the left-hand vertical map is a localization that inverts $f$, and the composition of the top horizontal map and the right-hand vertical map takes $f$ to a unit, so the lower horizontal map making the diagram commute exists and is unique. Moreover, if $y \in |Y|$ and $x = p(y) \in |X|$, then, in the diagram
$$\xymatrix@C=12mm{
{ \mathcal{O}_X(|X|) } \ar[r]^-{\phi_X} \ar[d] &
{ \mathcal{O}_Y(|Y|) } \ar[d] \cr
{ \mathcal{O}_{X,x} } \ar[r]^-{\phi_x} &
{ \mathcal{O}_{Y,y}, } \cr
}$$
the left-hand vertical map is a localization that inverts the homogeneous elements $f \in \mathcal{O}_X(|X|)$ with invertible image under the composition of the top horizontal map and the right-hand vertical map. So the lower horizontal map is a local map of local Dirac rings. This completes the definition of $\eta_Y \colon Y \to X$. Finally, it follows immediately from the definitions that $\epsilon$ and $\eta$ satisfy the triangle identities.
\end{proof}

\begin{theorem}\label{thm:limitscolimitsschemes}The category of Dirac schemes admits finite limits, coproducts indexed by small sets, and quotients by equivalence relations
$$\xymatrix@C=11mm{
{ R } \ar[r]^-{(s,t)} &
{ Y \times Y } \cr
}$$
such that $Y = \coprod_{i \in I}Y_i$ and $R = \coprod_{(i,j) \in I \times I}U_{i,j}$ and such that $s$ and $t$ restrict to open immersions $s|_{U_{i,j}} \colon U_{i,j} \to Y_i$ and $t|_{U_{i,j}} \colon U_{i,j} \to Y_j$.
\end{theorem}

\begin{proof}We first consider colimits. The case of small coproducts is clear. Given an equivalence relation as is the statement, we let $p \colon |Y| \to |X|$ be a coequalizer of $r,s \colon |R\,| \to |Y|$ in the category of topological spaces and continuous maps. We claim that for all $i \in I$, the map $p|_{Y_i} \colon |Y_i| \to p(|Y_i|)$ is a homeomorphism. First, it is a bijection, because the maps $s|_{U_{i,i}} \colon |U_{i,i}| \to |Y_i|$ and $t|_{U_{i,i}} \colon |U_{i,i}| \to |Y_i|$ necessarily are equal. Indeed, they are both open embeddings, and the diagonal map $\Delta \colon Y_i \to Y_i \times Y_i$ factors through $(s,t)|_{U_{i,i}} \colon U_{i,i} \to Y_i \times Y_i$, because $(s,t)$ is an equivalence relation. Second, it is an open map. For if $|V| \subset |Y_i|$ is open, then so is
$$\textstyle{ f^{-1}(f(|V|)) = \coprod_{j \in I} (t \circ s^{-1})(|V| \cap |U_{i,j}|) \subset \coprod_{j \in I} |Y_j| = |Y|. }$$
This proves the claim. Finally, the sheaf of Dirac rings $\mathcal{O}_X$ given by the equalizer
$$\xymatrix@C+=10mm{
{ \mathcal{O}_X } \ar[r]^-{f^{\sharp}} &
{ f_*\mathcal{O}_Y } \ar@<.7ex>[r]^-{f_*s^{\sharp}} \ar@<-.7ex>[r]_-{f_*t^{\sharp}} &
{ h_* \mathcal{O}_R } \cr
}$$
where $h = f \circ s = f \circ t$, makes $(|X|,\mathcal{O}_X)$ a Dirac scheme and makes it the stated quotient in the category of Dirac schemes.

We next consider limits. The opposite of the category of Dirac rings admits all small limits, and \cref{prop:diracscheme} implies that the functor $\operatorname{Spec}$ preserves them. In particular, the statement holds for finite diagrams of affine Dirac schemes, including the empty diagram. So it remains to prove that if $f \colon Y \to X$ and $g \colon X' \to X$ are maps of Dirac schemes, then a fiber product
$$\xymatrix{
{ Y' } \ar[r]^-{g'} \ar@<-.3ex>[d]^-{f'} &
{ Y } \ar[d]^-{f} \cr
{ X' } \ar[r]^{g} &
{ X } \cr
}$$
exists. We already know that this holds, if $X$, $X'$, and $Y$ are all affine. Finally, the general case is proved as in~\cite[Th\'{e}r\`{e}me~I.3.2.6]{egaI} or~\cite[\href{https://stacks.math.columbia.edu/tag/01JM}{Tag 01JM}]{stacks-project} by successively writing $Y$, $X'$, and finally $X$ as quotients of coproducts of affine Dirac schemes by equivalence relations of the form considered above.
\end{proof}

\begin{remark}\label{rem:locallydiracringedspaceslimitscolimits}Alternatively, one may prove that the category of locally Dirac ringed spaces admits all small limits and colimits; compare~\cite[Corollary~5]{gillam} and~\cite[Proposition~I.1.6]{demazuregabriel}. One then checks that the full subcategory of Dirac schemes is closed under finite limits and under the colimits enumerated in \cref{thm:limitscolimitsschemes}.
\end{remark}

\begin{example}\label{ex:dirac_scheme_colimit_of_affine_open_subschemes}Let us show that for every Dirac scheme $X$, the unique map
$$\xymatrix{
{ \varinjlim_{S \subset X} S } \ar[r]^-{f} &
{ X } \cr
}$$
from the colimit of its affine open sub-Dirac schemes is an isomorphism. The colimit in question is given by $Y/R$ with $Y = \coprod_{S \subset X}Y_S$ and $R = \coprod_{S,T \subset X}U_{S,T}$, where $Y_S = S$ and $U_{S,T} = S \cap T$, and where $s|_{U_{S,T}} \colon U_{S,T} \to Y_S$ and $t|_{U_{S,t}} \colon U_{S,T} \to T$ are given by the respective open immersions $S \cap T \to S$ and $S \cap T \to T$. So \cref{thm:limitscolimitsschemes} shows that the colimit exists and we must show that the map $f$ is an isomorphism. But the map of underlying spaces is a homeomorphism, because continuity is a local property, and the relevant diagram of sheaves
$$\xymatrix{
{ \mathcal{O}_X } \ar[r] &
{ g_*\mathcal{O}_Y } \ar@<.7ex>[r] \ar@<-.7ex>[r] &
{ h_*\mathcal{O}_R } \cr
}$$
is a limit diagram, because, taking sections of $U \subset X$ affine open, we get
$$\xymatrix{
{ \mathcal{O}_X(U) } \ar[r] &
{ \prod_{S \subset X}\mathcal{O}_X(U \cap S) } \ar@<.7ex>[r] \ar@<-.7ex>[r] &
{ \prod_{S,T \subset X}\mathcal{O}_X(U \cap S \cap T), } \cr
}$$
which is a limit diagram of sets, since $\mathcal{O}_X$ is a sheaf.
\end{example}

Given \cref{thm:limitscolimitsschemes}, we can adhere to Grothendieck's dictum that for a property to be considered a geometric property, it must be defined for all morphisms and be stable under base-change. Among geometric properties of a map of Dirac schemes, the familiar properties of being an open immersion, a closed immersion, and a flat map are all geometric properties. 

\subsection{Evenness}
\label{subsection:evenness_in_dirac_schemes}

A pleasant property of the Dirac--Zariski space associated with a Dirac ring $A$ is that it only depends on the sub-Dirac ring $A^{\operatorname{ev}} \subset A$ spanned by the homogeneous elements of even degree.

\begin{proposition}\label{prop:evendiraczariskispace}
Let $A$ be a Dirac ring, and let $A^{\operatorname{ev}} \subset A$ be the sub-Dirac ring spanned by the homogeneous elements of even degree, and let $\eta_X \colon X \to X^{\operatorname{ev}}$ be the map of Dirac spectra induced by the canonical inclusion.
\begin{enumerate}
\item[{\rm (1)}]The map of Dirac--Zariski spaces $\eta_X \colon |X| \to |X^{\operatorname{ev}}|$ is a homeomorphism.
\item[{\rm (2)}]For every homogeneous element $f \in A$ of even degree, the canonical map
$$\xymatrix{
{ \mathcal{O}_{X^{\operatorname{ev}}}(|X_f^{\operatorname{ev}}|) \otimes_{\mathcal{O}_{X^{\operatorname{ev}}}(|X|)} \mathcal{O}_X(|X|) } \ar[r] &
{ \mathcal{O}_X(|X_f|) } \cr
}$$
is an isomorphism.
\end{enumerate}
\end{proposition}

\begin{proof}
The inverse of the map in~(1) takes a graded prime ideal $\mathfrak{q} \subset A^{\operatorname{ev}}$ to the radical closure $\mathfrak{p} = \sqrt{\mathfrak{q}A} \subset A$ of the graded ideal $\mathfrak{q}A = A\mathfrak{q} \subset A$, which is proper since $(A \mathfrak{q})^{ev} = \mathfrak{q}$; see \cref{proposition:inclusion_of_even_subring_integral_and_a_universal_homeomorphism} for details. To prove~(2), we must show that the canonical map $A_f^{\operatorname{ev}} \otimes_{A^{\operatorname{ev}}}A \to A_f$ is an isomorphism, but this follows immediately from \cref{prop:calculusoffractions}.
\end{proof}

If a Dirac ring $A$ is even, then the underlying ring $\smash{\widetilde{A}}$ is commutative, and hence, it is meaningful to compare the Dirac spectrum $X = \operatorname{Spec}(A)$ to the prime spectrum $\smash{\widetilde{X}} = \operatorname{Spec}(\smash{\widetilde{A}})$. The grading on $A$ determines and is determined by an action
$$\xymatrix{
{ \mathbb{G}_m \times \smash{\widetilde{X}} } \ar[r]^-{\mu} &
{ \smash{\widetilde{X}} } \cr
}$$
of the multiplicative group $\mathbb{G}_m$ on $\smash{\widetilde{X}}$. We first prove a general lemma concerning the action by a group scheme on a scheme.

\begin{lemma}\label{lem:orbitspace}Let $\mu \colon G \times Y \to Y$ be an action by a group scheme $G$ on a scheme $Y$, and let $p \colon G \times Y \to Y$ be the canonical projection. The image of the map
$$\xymatrix@C=14mm{
{ |G \times Y| } \ar[r]^-{(|p|,|\mu|)} &
{ |Y| \times |Y| } \cr
}$$
is an equivalence relation.
\end{lemma}

\begin{proof}It is clear that the relation is reflexive and symmetric, but we must prove that it is also transitive. Let $[G \backslash Y]$ be the simplicial scheme given by the bar construction of the action. It is a groupoid object in the category of schemes in the sense of~\cite[Definition~6.1.2.7]{luriehtt}. In particular, the diagram
$$\xymatrix{
{ G \times G \times Y } \ar[r]^-{p'} \ar[d]^-(.47){\mu'} &
{ G \times Y } \ar[d]^-{\mu} \cr
{ G \times Y } \ar[r]^-{p} &
{ Y } \cr
}$$
with $p'$ the projection on the second and third factors and with $\mu' = G \times \mu$ is cartesian. The induced diagram of underlying spaces
$$\xymatrix{
{ |G \times G \times Y| } \ar[r]^-{|p'|} \ar[d]^-(.47){|\mu'|} &
{ |G \times Y| } \ar[d]^-{|\mu|} \cr
{ |G \times Y| } \ar[r]^-{|p|} &
{ |Y| } \cr
}$$
is typically not cartesian. However, by~\cite[Proposition~3.4.7]{egaI}, the induced map
$$\xymatrix@C=14mm{
{ |G \times G \times Y| } \ar[r]^-{(|p'|,|\mu'|)} &
{ |G \times Y| \times_{|Y|} |G \times Y| } \cr
}$$
is surjective, and this implies the transitivity of the relation.
\end{proof}

We say that the quotient of the topological space $|Y|$ by the equivalence relation in \cref{lem:orbitspace} is the orbit space of the action $\mu \colon G \times Y \to Y$, and we say that the subspace of $|Y|$ consisting of the points that are equivalent to $y \in |Y|$ is the $G$-orbit through $y \in |Y|$. It is equal to the image of the composite map
$$\xymatrix@C=12mm{
{ |G \times \operatorname{Spec}(k(y))| } \ar[r]^-{|G \times i_y|} &
{ |G \times Y| } \ar[r]^-{|\mu|} &
{ |Y| } \cr
}$$
and is generally not a locally closed subspace.

\begin{proposition}
\label{prop:orbitspace}
Let $A$ be an even Dirac ring with Dirac spectrum $X$, and let $\smash{\widetilde{A}}$ be its underlying ring with spectrum $\smash{\widetilde{X}}$. The map
$$\xymatrix{
{ |\smash{\widetilde{X}}| } \ar[r]^-{p} &
{ |X| } \cr
}$$
that to a prime ideal $\mathfrak{q} \subset \smash{\widetilde{A}}$ assigns the largest graded prime ideal $\mathfrak{p} \subset A$ with $\mathfrak{p} \subset \mathfrak{q}$ exhibits $|X|$ as the orbit space of the action $\mu \colon \mathbb{G}_m \times \smash{\widetilde{X}} \to \smash{\widetilde{X}}$ defined by the grading. There is a canonical isomorphism $\smash{\widetilde{\mathcal{O}}}_X \to p_*(\mathcal{O}_{\smash{\widetilde{X}}})$ of sheaves of rings on $|X|$.
\end{proposition}

\begin{proof}The map $p$ is continuous, since for $f \in A$ homogeneous, we have
$$p^{-1}(|X_f|) = |\smash{\widetilde{X}}_f|,$$
and the inclusion $s \colon |X| \to |\smash{\widetilde{X}}|$ is also continuous, since for any $g \in \smash{\widetilde{A}}$, we have
$$\textstyle{ s^{-1}(|\smash{\widetilde{X}}_g|) = \bigcup_{i \in \mathbb{Z}} |X_{g_i}|, }$$
where we write $g = \sum_{i \in \mathbb{Z}}g_i$ with $g_i \in A$ homogeneous. Indeed, if $\mathfrak{p} \subset A$ is a graded ideal, then $g \notin \mathfrak{p}$ if and only if $g_i \notin \mathfrak{p}$ for some $i \in \mathbb{Z}$. So $s$ is a continuous section of the continuous map $p$, which immediately implies that $p$ is a quotient map. 

We proceed to show that the fiber $p^{-1}(x) \subset |\smash{\widetilde{X}}|$ lying above $x \in |X|$ is equal to the $\mathbb{G}_m$-orbit through $s(x) \in |\smash{\widetilde{X}}|$. Since, by \cref{lem:orbitspace}, the partition of $|\smash{\widetilde{X}}|$ into $\mathbb{G}_m$-orbits is an equivalence relation, this will prove that $p \colon |\smash{\widetilde{X}}| \to |X|$ exhibits $|X|$ as the orbit space of the $\mathbb{G}_m$-action on $\smash{\widetilde{X}}$, as stated.

So let $x \in |X|$. The canonical map of Dirac rings $i_x \colon A \to k(x)$ from $A$ to the Dirac residue field at $x$ induces a $\mathbb{G}_m$-equivariant embedding of schemes
$$\xymatrix{
{ \operatorname{Spec}(\widetilde{k(x)}) } \ar[r]^-{i_x} &
{ \smash{\widetilde{X}}, } \cr
}$$
whose set-theoretic image we now identify with the fiber $p^{-1}(x) \subset |\smash{\widetilde{X}}|$. To this end, we first notice that if $\mathfrak{q} \subset \smash{\widetilde{A}}$ is a prime ideal, then the largest graded ideal $\mathfrak{a} \subset A$ with the property that $\mathfrak{a} \subset \mathfrak{q}$ is automatically a prime ideal. Indeed, if $f,g \in A$ are homogeneous and $fg \in \mathfrak{a}$, then $fg \in \mathfrak{q}$, so either $f \in \mathfrak{q}$ or $g \in \mathfrak{q}$ or both. Therefore, the maximality of $\mathfrak{a} \subset \mathfrak{q}$ implies that $f \in \mathfrak{a}$ or $g \in \mathfrak{a}$ or both, so \cref{lem:primeideal} shows that $\mathfrak{a} \subset A$ is a prime ideal. It follows that if $\mathfrak{p} \subset A$ is the graded prime ideal corresponding to $x \in |X|$, then $y \in p^{-1}(x) \subset |\smash{\widetilde{X}}|$ if and only if the corresponding prime ideal $\mathfrak{q} \subset \smash{\widetilde{A}}$ contains $\mathfrak{p}$ and $\mathfrak{q}$ does not contain any homogeneous elements from $A \smallsetminus \mathfrak{p}$. But these are exactly the points $y \in |\smash{\widetilde{X}}|$ that belong to the set-theoretic image of the embedding $i_x$ above.

Now, if $y = s(x) \in |\smash{\widetilde{X}}|$, then we have a diagram of schemes with $\mathbb{G}_m$-action
$$\xymatrix{
{ \mathbb{G}_m \times \operatorname{Spec}(k(y)) } \ar[rr]^-{\mathbb{G}_m \times i_y} \ar[d] &&
{ \mathbb{G}_m \times \smash{\widetilde{X}} } \ar[d]^-{\mu} \cr
{ \operatorname{Spec}(\widetilde{k(x)}) } \ar[rr]^-{i_x} &&
{ \smash{\widetilde{X}} } \cr
}$$
with the left-hand vertical map induced by the canonical map $\widetilde{k(x)} \to k(y)$ from the underlying ring of the Dirac residue field at $x$ to the residue field at $y = s(x)$. On the one hand, the lower horizontal map is an embedding, whose set-theoretic image we have identified with $p^{-1}(x) \subset |\smash{\widetilde{X}}|$, and on the other hand, the set-theoretic image of the composition of the top horizontal map and the right-hand vertical map is equal to the  $\mathbb{G}_m$-orbit through $y \in |\smash{\widetilde{X}}|$. Thus, it suffices to show that the left-hand vertical map is surjective. The map in question is induced by a map of Dirac rings $h \colon k(x) \to \mathbb{Z}[t^{\pm1}] \otimes k(y)$, where $t$ is a spin $1$ generator. But, by \cref{prop:classification_of_dirac_fields}, the Dirac field $k(x)$ is either a field $k$ or a Laurent polynomial ring $k[s^{\pm1}]$ on a generator of nonzero degree, and in either case, the map $h$ is faithfully flat. So the left-hand map in the diagram is surjective as desired.

Finally, for $f \in A$ homogeneous, we have canonical isomorphisms
$$p_*(\mathcal{O}_{\smash{\widetilde{X}}})(|X_f|) \simeq \mathcal{O}_{\smash{\widetilde{X}}}(|\smash{\widetilde{X}}_f|) \simeq \smash{\widetilde{A}}_f,$$
which proves that $p_*(\mathcal{O}_{\smash{\widetilde{X}}}) \simeq \smash{\widetilde{\mathcal{O}}}_X$.
\end{proof}

\begin{remark}\label{rem:orbitspace}1)~The section $s \colon |X| \to |\smash{\widetilde{X}}|$ of the projection $p \colon |\smash{\widetilde{X}}| \to |X|$ onto the orbit space for the action $\mu \colon \mathbb{G}_m \times \smash{\widetilde{X}} \to \smash{\widetilde{X}}$ assigns to $x \in |X|$ the generic point of the $\mathbb{G}_m$-orbit $p^{-1}(x) \subset |\smash{\widetilde{X}}|$. Thus, it is the spectral nature of these spaces that is the reason that this section exists. Note that since $s$ is a section of $p$ and both are continuous, the map $s$ identifies $|X|$ with a subspace of $|\tilde{X}|$. 

\noindent 2)~The sheaf of Dirac rings $\mathcal{O}_X$ has underlying sheaf of rings given by the direct image $p_*(\mathcal{O}_{\smash{\widetilde{X}}})$. However, we note that, whereas $(|X|,\mathcal{O}_X)$ is a locally Dirac ringed space, its underlying ringed space $(|X|,p_*(\mathcal{O}_{\smash{\widetilde{X}}}))$ is typically not locally ringed.
\end{remark}

\begin{remark}\label{rem:even_subring_has_the_same_spectrum}Let $S$ be an even Dirac ring of non-negative spin. The ideal $S_+ \subset S$ spanned by the homogeneous elements of positive spin is graded, and the canonical maps $S_0 \to S \to S/S_+$ identifies the quotient with the subring of spin zero elements. Thus, the Dirac--Zariski space of $S$ decomposes as the recollement
$$\xymatrix{
{ |\operatorname{Spec}(S_0)| } \ar[r]^-{\phantom{j}i\phantom{j}} &
{ |\operatorname{Spec}(S)| } &
{ |\operatorname{Proj}(S)| } \ar[l]_-{j} \cr
}$$
of the Zariski space of $S_0$, which is a closed subspace, and the projective Zariski space of $S$, which is an open subspace. Note that $\operatorname{Proj}(S)$ is not an open sub-Dirac scheme of $\operatorname{Spec}(S)$. Instead, the structure sheaf of $\operatorname{Proj}(S)$ is defined to be the spin zero part of $\mathcal{O}_{\operatorname{Spec}(S)}|_{\operatorname{Proj}(S)}$. We also note that, while the space $|\operatorname{Spec}(S)|$ is always quasi-compact, the open subspace $|\operatorname{Proj}(S)|$ is generally not so.
\end{remark}

\begin{example}
\label{example:non_quasi_finite_dirac_algebra_with_finite_spectrum}
Let $k$ be a field. The free Dirac $k$-algebra $B = k[\gamma]$ on a generator of non-zero integer spin-$m$ has Dirac--Zariski space
$$|\mathbb{A}_k^1(m)| = |X| = \{\eta,s\}.$$
It is a Sierpi\'{n}ski space with $\eta = (0)$ open and $s = (\gamma)$ closed. So the Dirac--Zariski space $|X|$ is finite even though the $k$-algebra $B$ is not finite, despite being finite type over a field. 
\end{example}

\begin{remark}
\label{remark:sheaf_of_dirac_rings_compared_to_even_case}Let $A$ be a Dirac ring, let $A^{\operatorname{ev}} \subset A$ be the sub-Dirac-ring generated by the homogeneous elements of even degree, and let $r \colon X \to X^{\operatorname{ev}}$ be the induced map of Dirac spectra. If $f \in A^{\operatorname{ev}}$ is a homogeneous element of even degree, then, as we have observed in \cref{prop:evendiraczariskispace}, the square of Dirac schemes
$$\xymatrix@C10mm{
{ X_f } \ar[r]^-{\eta_{X_f}} \ar[d]^-{f} &
{ X_f^{\operatorname{ev}} } \ar[d]^-{f^{\operatorname{ev}}} \cr
{ X } \ar[r]^-{\eta_X} &
{ X^{\operatorname{ev}} } \cr
}$$
is cartesian, and the horizontal maps are homeomorphisms on underlying spaces.
\end{remark}

The next result demonstrates that, in the case of even Dirac rings, the geometry of the Dirac spectrum exerts strong control over the geometry of the much larger prime spectrum of the underlying ring.

\begin{proposition}
\label{prop:open_embeddings_of_even_dirac_rings_detected_on_underlying_rings}
If $\phi \colon A \to B$ is a map between even Dirac rings, then the following are equivalent:
\begin{enumerate}
\item[{\rm (1)}]The map of Dirac spectra $q \colon Y \to X$ associated with $\phi \colon A \to B$ is an open immersion in the sense that the map $q \colon |Y| \to |X|$ is an open embedding and the map $\phi \colon q^*(\mathcal{O}_X) \to \mathcal{O}_Y$ is an isomorphism.
\item[{\rm (2)}] There exists a family $(f_i)_{i \in I}$ of homogeneous elements of $A$ such that $(\phi(f_i))_{i \in I}$ generates the unit ideal of $B$ and such that the induced map $\phi_{f_i} \colon A_{f_i} \to B_{f_i}$ is an isomorphism for all $i \in I$.
\item[{\rm (3)}]The map of prime spectra $\widetilde{q} \colon \widetilde{Y} \to \widetilde{X}$ associated with the map $\phi \colon \widetilde{A} \to \widetilde{B}$ of underlying rings is an open immersion of schemes.
\end{enumerate}
\end{proposition}

\begin{proof}It is clear that~(1) implies~(2) and that~(2) implies~(3), so we assume~(3) and prove~(1). We must show that $p \colon |Y| \to |X|$ is an open embedding. By \cref{rem:orbitspace}, we have a diagram of topological spaces and continuous maps
$$\xymatrix@C+=10mm{
{ |\smash{\widetilde{Y}}| } \ar[r]^-{\widetilde{q}} &
{ |\smash{\widetilde{X}}| } \cr
{ |Y| } \ar[r]^-{q} \ar[u]_-{s} &
{ |X| } \ar[u]_-{s} \cr
}$$
in which the vertical maps are embeddings, and by assumption, the top horizontal map is an open embedding. Therefore, to show that the bottom horizontal map is an open embedding, it suffices to show that the underlying diagram of sets is cartesian. It is clear that the map from $|Y|$ to the fiber product is injective. To prove that it is also surjective, we fix $x \in |X|$ and show that if the fiber $Y_x$ is empty, then so is the fiber $\smash{\widetilde{Y}}_{s(x)}$. But if $\mathfrak{p} \subset A$ is the graded prime ideal corresponding to $x \in |X|$, then the underlying ring of the Dirac residue field $k(x)$ and the residue field $k(s(x))$ are both localizations of the underlying ring of $A/\mathfrak{p}$. Hence, there is a map
$$\xymatrix{
{ \widetilde{k(x)} } \ar[r] &
{ k(s(x)) } \cr
}$$
of $\smash{\widetilde{A}}$-algebras, so if $B \otimes_Ak(x)$ is the zero Dirac ring, then $\smash{\widetilde{B}} \otimes_{\smash{\widetilde{A}}}k(s(x))$ is the zero ring, as we wanted to prove. So $q \colon |Y| \to |X|$ is an open embedding.

It remains to prove that $\phi \colon q^*(\mathcal{O}_X) \to \mathcal{O}_Y$ is an isomorphism. So let $f \in A$ be a homogeneous element such that $|X_f| \subset q(|Y|) \subset |X|$. We claim that
$$|\widetilde{X}_f| \subset \widetilde{q}\,(|\widetilde{Y}|) \subset |\widetilde{X}|.$$
Unwrapping the definitions, we must show that if $\mathfrak{q} \subseteq \smash{\widetilde{A}}$ that does not belong to the image of $|\smash{\widetilde{Y}}|$, then it contains $f$. A prime ideal $\mathfrak{q} \subset \smash{\widetilde{A}}$ does not belong to the image of $|\smash{\widetilde{Y}}|$ if and only if $\phi(\mathfrak{q})\smash{\widetilde{B}} = \smash{\widetilde{B}}$ if and only if we can write
$$\textstyle{ 1 = \sum_{i \in I} f_i g_i \in \smash{\widetilde{B}} }$$
with $f_i \in \mathfrak{q}$ and $g_i \in \widetilde{B}$. Expanding each of these elements as a sum of homogeneous elements, we can assume that they are themselves homogeneous. Thus, we also have $\phi(p(\mathfrak{q})) B = B$, where $p(\mathfrak{q}) \subset A$ is the graded prime ideal generated by the family consisting of the homogeneous elements in $\mathfrak{q}$. It follows that $p(\mathfrak{q})$ is not in the image of $|Y|$ , so $f \in p(\mathfrak{q}) \subset \mathfrak{q}$. This proves the claim.

Finally, since localization commutes with passing to underlying rings, we deduce that if $f \in A$ is as above, then $\phi_f \colon A_{f} \to B_{f}$ is an isomorphism, because the map of underlying rings is so, by assumption.
\end{proof}

We say that a Dirac scheme $X$ is even if for every open subset $U \subseteq |X|$, the Dirac ring $\mathcal{O}_{X}(U)$ is concentrated in even degrees. In \cref{ex:diracscheme}, we defined for every Dirac scheme $X$, a canonical map $\eta_{X} \colon X \to X^{\operatorname{ev}}$ to an even Dirac scheme. This map is characterized by the following universal property:

\begin{lemma}
\label{lem:evenadjunction}
The inclusion of the full subcategory spanned by even Dirac schemes into the category of all Dirac schemes admits the left adjoint functor with unit map $\eta_X \colon X \to X^{\operatorname{ev}}$. Moreover, the projection formula holds: Given a map of even Dirac schemes $p \colon W \to V$ and a map of Dirac schemes $q \colon Y \to V$, the canonical map
$$\xymatrix{
{ (Y \times_VW)^{\operatorname{ev}} } \ar[r] &
{ Y^{\operatorname{ev}} \times_VW } \cr
}$$
is an isomorphism.
\end{lemma}

\begin{proof}
Since the right adjoint is the inclusion of a full subcategory, the counit is necessarily an isomorphism. So the statement that $\eta_X \colon X \to X^{\operatorname{ev}}$ is the unit of an adjunction is equivalent to the statement that for any map $f = (p,\phi) \colon X \to V$ of Dirac schemes with $V$ even, the map $\phi \colon p^*\mathcal{O}_{V} \to \mathcal{O}_{X}$ of sheaves of Dirac rings on $|X|$ factors uniquely through the inclusion $\mathcal{O}_{X^{\operatorname{ev}}} \to \mathcal{O}_X$. But this is clear.

The map in the projection formula is the adjunct of the map
$$\xymatrix@C=14mm{
{ Y \times_VW } \ar[r]^-{\eta_Y \times_VW} &
{ Y^{\operatorname{ev}} \times_VW, } \cr
}$$
and to prove that it is an isomorphism, we may reduce, by the construction of fiber products, to the case, where first $V$ and then $W$ and $Y$ are affine. But in the affine case, the statement is clear.
\end{proof}

The absolute property of being even can be naturally extended to a property of maps. This yields a new geometric property of Dirac schemes.

\begin{definition}
\label{definition:even_morphism_of_dirac_schemes}
A map of Dirac schemes $p \colon Y \to X$ is even if the square
$$\xymatrix@C=10mm{
{ Y } \ar[r]^-{\eta_Y} \ar[d]^-{p} &
{ Y^{\operatorname{ev}} } \ar[d]^-{p^{\operatorname{ev}}} \cr
{ X } \ar[r]^-{\eta_X} &
{ X^{\operatorname{ev}} } \cr
}$$
is cartesian.
\end{definition}

Since the terminal Dirac scheme $\operatorname{Spec}(\mathbb{Z})$ is even, a Dirac scheme $X$ is even if and only if the unique map $X \to \operatorname{Spec}(\mathbb{Z})$ is even.

\begin{lemma}
\label{lemma:even}
If $p \colon Y \to X$ is an even map of Dirac schemes, then so is its base-change $p' \colon Y' \to X'$ along any map of Dirac schemes $q \colon X' \to X$.
\end{lemma}

\begin{proof}We wish to show that, in the diagram
$$\xymatrix@C=12mm{
{ Y' } \ar[r]^-{\eta_{Y'}} \ar[d]^-{p'} &
{ Y'{}^{\operatorname{ev}} } \ar[r]^-{q'{}^{\operatorname{ev}}} \ar[d]^-{p'{}^{\operatorname{ev}}} &
{ Y^{\operatorname{ev}} } \ar[d]^-{p^{\operatorname{ev}}} \cr
{ X' } \ar[r]^-{\eta_{X'}} &
{ X'{}^{\operatorname{ev}} } \ar[r]^-{q^{\operatorname{ev}}} &
{ X^{\operatorname{ev}} } \cr
}$$
the left-hand square is cartesian. This will follow, once we prove that the outer square and the right-hand square both are cartesian. The outer square agrees with the outer square in the diagram
$$\xymatrix@C=12mm{
{ Y' } \ar[r]^-{q'} \ar[d]^-{p'} &
{ Y } \ar[r]^-{\eta_Y} \ar[d]^-{p} &
{ Y^{\operatorname{ev}} } \ar[d]^-{p^{\operatorname{ev}}} \cr
{ X' } \ar[r]^-{q} &
{ X } \ar[r]^-{\eta_X} &
{ X^{\operatorname{ev}} } \cr
}$$
and it is cartesian, since, by assumption, the left-hand square and the right-hand square therein both are cartesian. Finally, the right-hand square in the top diagram is obtained by applying $(-)^{\operatorname{ev}}$ to the left-hand square, or equivalently, the outer square in the bottom diagram. It is cartesian by the projection formula of \cref{lem:evenadjunction}.
\end{proof}

Let us record that the unit map $\eta_X \colon X \to X^{\operatorname{ev}}$ has a number of good properties. We say that a map of Dirac schemes $p \colon Y \to X$ is affine if for every affine open $U \subset X$, the inverse image $V = p^{-1}(U) \subset Y$ is affine open, and it is integral if, in addition, every homogeneous element of $\mathcal{O}_Y(|V|)$ is a root of a homogeneous monic polynomial with coefficients in $\mathcal{O}_X(|U|)$. 

\begin{lemma}
\label{lemma:map_into_even_an_affine_universal_homeomorphism}
Let $X$ be a Dirac scheme. The unit map
$$\xymatrix{
{ X } \ar[r]^-{\eta_X} &
{ X^{\operatorname{ev}} } \cr
}$$
is affine, integral, and a universal homeomorphism.
\end{lemma}

\begin{proof}
If $U \subset X$ is affine open, then $U^{\operatorname{ev}} \subset X^{\operatorname{ev}}$ is affine open and $\eta_X^{-1}(U^{\operatorname{ev}}) = U$. So $\eta_X$ is affine. Since the properties of being integral and a universal homeomorphism are local in the Zariski-topology, these two properties follow from the corresponding statement about Dirac rings, which we prove in \cref{proposition:inclusion_of_even_subring_integral_and_a_universal_homeomorphism}. 
\end{proof}

\begin{corollary}
\label{cor:evenintegralradicial}
Let $p \colon Y \to X$ be a map of Dirac schemes. In the diagram
$$\begin{xy}
(-12,23)*+{ Y }="0";
(0,14)*+{ {}^{\phantom{\operatorname{ev}}}Y\smash{'} }="11";
(0,14)*+{\phantom{i}}="11a";
(20,14)*+{ Y^{\operatorname{ev}} }="12";
(0,0)*+{ {}^{\phantom{\operatorname{ev}}}X }="21";
(20,0)*+{ X^{\operatorname{ev}} }="22";
{ \ar@/^.8pc/^-{\eta_Y} "12";"0";};
{ \ar@/_.7pc/_-{p} "21";"0";};
{ \ar^-(.4){\eta_X{}'} "12";"11";};
{ \ar@<1ex>^-(.46){p^{\operatorname{ev}}{}'} "21";"11";};
{ \ar@<-.5ex>^-(.46){p^{\operatorname{ev}}} "22";"12";};
{ \ar^-{\eta_X} "22";"21";};
{ \ar^-{r} "11a";"0";};
\end{xy}$$
with the square cartesian, the map $r$ is integral and a universal  homeomorphism, and the map $p^{\operatorname{ev}}{}'$ is even.
\end{corollary}

\begin{proof}
By \cref{lemma:map_into_even_an_affine_universal_homeomorphism}, 
The projection formula shows that $r^{\operatorname{ev}} \colon Y^{\operatorname{ev}} \to Y'{}^{\operatorname{ev}}$ is an isomorphism of Dirac schemes, so in particular, the map of underlying spaces is a homeomorphism. To prove that $r$ is integral, we argue as in~\cite[\href{https://stacks.math.columbia.edu/tag/035D}{Tag 035D}]{stacks-project} that it suffices to show that $\eta_Y$ is integral and that $\eta_X'$ is separated. But this follows from \cref{lemma:map_into_even_an_affine_universal_homeomorphism}.
\end{proof}

\section{Commutative algebra}

In this section we generalize several results of classical commutative algebra to the setting of Dirac rings. In many cases, such as Hilbert's basis theorem, the proofs are the same in the classical case and so were certainly known to the experts, even if not written down in precisely this form. An exception here is our proof of Zariski's main theorem, which we believe is a genuinely new result.

\subsection{Properties of modules}
\label{subsection:commutative_algebra_flatness}

Let $A$ be a Dirac ring. The functor that to a graded $A$-module $M$ assigns the graded set $(M_d)_{d \in \mathbb{Z}}$ admits a left adjoint functor that to a graded set $S$ assigns a graded $A$-module $A(S)$. We define a family of homogeneous elements in $M$ to be a map of graded sets $x \colon S \to (M_d)_{d \in \mathbb{Z}}$, and we define it to be linearly independent (resp.\ to generate $M$, resp.\ to be a basis of $M$) if the adjunct map $\tilde{x} \colon A(S) \to M$ is a monomorphism (resp.\ an epimorphism, resp.\ an isomorphism). We say that a graded $A$-module $M$ is free if it admits a basis. We define a graded set $S = (S_d)_{d \in \mathbb{Z}}$ to be finite if the total set $\coprod_{d \in \mathbb{Z}}S_d$ is finite, and we define it to be even if $S_d = \emptyset$ for $d$ odd. Finally, we say that a family $x \colon S \to M$ of homogeneous elements in a graded $A$-module $M$ is finite (resp.\ even) if $S$ is finite (resp.\ even).

\begin{definition}\label{def:finiteevenmodule}Let $A$ be a Dirac ring.
\begin{enumerate}
\item[(1)]A graded $A$-module $M$ is finitely generated if there exists a finite family of homogeneous elements that generates $M$.
\item[(2)]A graded $A$-module $M$ is finitely presented if there exists a finite family of homogeneous elements $x \colon S \to (M_d)_{d \in \mathbb{Z}}$ that generates $M$ and for which the kernel of the adjunct map $\tilde{x} \colon A(S) \to M$ is finitely generated.
\item[(3)]A graded $A$-module $M$ is evenly generated if there exists an even family of homogeneous elements that generates $M$.
\item[(4)]A graded $A$-module $M$ is evenly presented if there exists an even family of homogeneous elements $x \colon S \to (M_d)_{d \in \mathbb{Z}}$ that generates $M$ and for which the kernel of the adjunct map $\tilde{x} \colon A(S) \to M$ is evenly generated.
\end{enumerate}
\end{definition}

\begin{remark}
\label{rem:finiteevenmodule}
Let $A$ be a Dirac ring, and let $M$ be a graded $A$-module. We have defined $A^{\operatorname{ev}} \subset A$ be the sub-Dirac ring generated by the family of even homogeneous elements in $A$, and we similarly define $M^{\operatorname{ev}} \subset M$ to be the graded sub-$A^{\operatorname{ev}}$-module generated by the family of even homogeneous elements in $M$. We remark that $M$ is evenly generated (resp.\ evenly presented) if and only if the canonical map
$$\xymatrix{
{ A \otimes_{A^{\operatorname{ev}}}M^{\operatorname{ev}} } \ar[r] &
{ M } \cr
}$$
is an epimorphism (resp.\ an isomorphism) of graded $A$-modules.
\end{remark}

\begin{remark}
If $A = A^{\operatorname{ev}}$ is an even ring, then an $A$-module $M$ is evenly generated if and only if it is concentrated in even degrees. In this case, it is also automatically evenly presented. 
\end{remark}

\begin{example}
\label{example:evenly_generated_not_evenly_presented_module}
Over a general Dirac ring $A$, there exist evenly generated $A$-modules which are not evenly presented. For example, if $A = k[x]$ is the free Dirac algebra on a generator of odd degree, then the unique $k$-algebra map $\phi \colon A \to k$ is evenly generated, but not evenly presented. 
\end{example}

An important statement about finitely generated modules is the following analogue of the classical Nakayama's lemma \cite[\href{https://stacks.math.columbia.edu/tag/00DV}{Tag 00DV}]{stacks-project}.  

\begin{proposition}[Nakayama's lemma]
\label{lem:nakayama}
Let $A$ be a local Dirac ring with maximal graded ideal $\mathfrak{m} \subset A$ and residue Dirac field $k = A/\mathfrak{m}$. If $M$ is a finitely generated graded $A$-module such that $M \otimes_Ak = 0$, then $M$ is zero.
\end{proposition}

\begin{proof}
We let $(x_1, \ldots, x_n)$ be a finite family of homogeneous elements of $M$ that generates $M$, and assume that the family is minimal with this property. If $n = 0$, then there is nothing to prove, and if  $n > 0$, then there exists, by assumption, a family $(a_1,\dots,a_n)$ of homogeneous elements of $\mathfrak{m}$ such that
$$\textstyle{ x_n = \sum_{1 \leq i \leq n} x_ia_i. }$$
In particular, the element $a_n \in \mathfrak{m}$ is homogeneous of  degree $0$ and $1 - a_n \notin \mathfrak{m}$, so we conclude that $1 -  a_n \in A$ is a unit. It follows that
$$\textstyle{ x_n = \sum_{1 \leq i < n} (1-a_n)^{-1}a_ix_i, }$$
but this contradicts the minimality of the family $(x_1,\dots,x_n)$.
\end{proof}

We now move on to flatness. 

\begin{definition}
\label{def:flat}
Let $A$ be a Dirac ring. A graded $A$-module $M$ is flat if the functor $M \otimes_A -$ from the abelian category of graded $A$-modules to itself is exact, and $M$ is faithful if the functor $M \otimes_A -$ is faithful.
\end{definition}

The properties of being flat and of being faithful are preserved under base change along any map $\phi \colon A \to B$ of Dirac rings. Let us observe that under that faithful flatness can be detected by residue fields. 

\begin{proposition}
\label{proposition:faithful_flatness_detected_by_tensors_with_residue_fields}
Let $A$ be a Dirac ring, and let $M$ be a flat graded $A$-module. The following conditions are equivalent:
\begin{enumerate}
\item[{\rm (1)}] The $A$-module $M$ is faithful.
\item[{\rm (2)}] For every graded prime ideal $\mathfrak{p} \subset A$, $k(\mathfrak{p}) \otimes_{A} M$ is nonzero.
\item[{\rm (3)}] For every maximal graded ideal $\mathfrak{m} \subset A$, $k(\mathfrak{m}) \otimes_{A} M$ is nonzero.
\end{enumerate}
\end{proposition}

\begin{proof}It is clear that~(1) implies~(2) and that~(2) implies~(3). So we assume~(3) and prove~(1). So we let $N$ be an $A$-module such that $M \otimes_{A} N$ is zero and prove that $N$ is zero. If not, then there exists a nonzero homogeneous element $x \in N$. Let $s = \operatorname{spin}(x) \in \mathbb{Z}$. In this situation, the map of graded $A$-modules
$$\xymatrix{
{ A(-s) } \ar[r]^-{h} &
{ N } \cr
}$$
defined by $h(a) = ax$ is nonzero, and therefore, its kernel $I(-s) \subset A(-s)$ is a proper submodule. Since $M$ is flat, the induced map
$$\xymatrix@C+=11mm{
{ M \otimes_A (A/I)(-s) } \ar[r]^-{\bar{h} \otimes M} &
{ M \otimes_{A} N } \cr
}$$
is injective, and since the right-hand side is zero, so is the left-hand side. But this implies that $M \otimes_{A} k(\mathfrak{m}) = 0$ for every maximal graded ideal $\mathfrak{m} \subset A$ that contains the graded ideal $I \subset A$, which contradicts~(3). So we conclude that $N$ is necessarily zero, which proves~(1).
\end{proof}

The following result, due to Lazard, which also appears in the graded context in unpublished work of Davies \cite{davies_lazard_theorem}, gives two alternative characterizations of flatness, both of which we will need in what follows. We refer to~(2) as the equational criterion for flatness.

\begin{theorem}[Lazard]
\label{thm:lazard}
Let $A$ be a Dirac ring, and let $M$ be a graded $A$-module. The following are equivalent.
\begin{enumerate}
\item[{\rm (1)}]The graded $A$-module $M$ is flat.
\item[{\rm (2)}]Given maps of graded $A$-modules
$$\xymatrix{
{ F' } \ar[r]^-{a} &
{ F } \ar[r]^-{x} &
{ M } \cr
}$$
with $F$ and $F'$ finitely generated free and $xa = 0$, there exists a factorization
$$\begin{xy}
(0,12)*+{ {}\phantom{'}F }="11";
(19.5,12)*+{ F'' }="12";
(10,0)*+{ M }="21";
{ \ar@<-.2ex>^-{b} "12";"11";};
{ \ar_-(.35){x} "21";"11";};
{ \ar^-(.35){y} "21";"12";};
\end{xy}$$
with $F''$ finitely generated free and $ba = 0$.
\item[{\rm (3)}]The graded $A$-module $M$ is a filtered colimit of finitely generated free graded $A$-modules.
\end{enumerate}
\end{theorem}

\begin{proof}
We first assume~(1) and prove~(2). We have a diagram of graded $A$-modules
$$\begin{xy}
(1,14)*+{ F\smash{'} }="11";
(16,14)*+{ F }="12";
(32,14)*+{ C }="13";
(44,14)*+{ 0 }="14";
(16,0)*+{ M }="22";
{ \ar^-{\phantom{p}a\phantom{p}} "12";"11";};
{ \ar^-{p} "13";"12";};
{ \ar^-{x} "22";"12";};
{ \ar "14";"13";};
{ \ar@<.3ex>^-(.35){\bar{x}} "22";"13";};
\end{xy}$$
with the top row exact. Moreover, since $F$ and $F'$ are dualizable, and since $M$ is flat, this diagram determines and is determined by the diagram
$$\begin{xy}
(1,14)*+{ M \otimes F'{}^{\vee} }="11";
(28,14)*+{ M \otimes F^{\vee} }="12";
(54,14)*+{ M \otimes C^{\vee} }="13";
(72,14)*+{ 0. }="14";
(28,0)*+{ A }="22";
{ \ar_-(.45){M \otimes a^{\vee}} "11";"12";};
{ \ar_-(.45){M \otimes p^{\vee}} "12";"13";};
{ \ar "13";"14";};
{ \ar_-{\tilde{x}} "12";"22";};
{ \ar@<-.5ex>_-{\tilde{\bar{x}}} "13";"22";};
\end{xy}$$
The map $\tilde{\bar{x}}$ determines and is determined by the element $\tilde{\bar{x}}(1) \in M \otimes C^{\vee}$, which is homogeneous of degree zero. This element, in turn, can be written, non-canonically, as a finite sum pure tensors, and therefore, the map $\tilde{\bar{x}}$ admits a factorization
$$\begin{xy}
(1,14)*+{ M \otimes F'{}^{\vee} }="11";
(28,14)*+{ M \otimes F^{\vee} }="12";
(54,14)*+{ M \otimes C^{\vee} }="13";
(72,14)*+{ 0 }="14";
(28,0)*+{ A }="22";
(54,0)*+{ M \otimes F''{}^{\vee} }="23";
{ \ar_-(.45){M \otimes a^{\vee}} "11";"12";};
{ \ar_-(.45){M \otimes p^{\vee}} "12";"13";};
{ \ar "13";"14";};
{ \ar_-{\tilde{x}} "12";"22";};
{ \ar^-{\tilde{y}} "23";"22";};
{ \ar_-{M \otimes c} "13";"23";};
\end{xy}$$
with $F''$, and hence, its dual finitely generated free. So if $b \colon F \to F''$ the unique map such that $b^{\vee} = p^{\vee}c$, then $x = yb$ with $ba = 0$, which proves~(2).

We next prove that~(2) implies~(3). Let $\operatorname{Mod}_A^{\operatorname{ff}}(\operatorname{Ab}) \subset \operatorname{Mod}_A(\operatorname{Ab})$ be the full subcategory of the category of graded $A$-modules spanned by the finitely generated free graded $A$-modules. Given any graded $A$-module $M$, the diagram
$$\xymatrix{
{ (\operatorname{Mod}_A^{\operatorname{ff}}(\operatorname{Ab})_{/M})^{\triangleright} } \ar[r]^-{\bar{p}} &
{ \operatorname{Mod}_A(\operatorname{Ab}) } \cr
}$$
that to $(F,x \colon F \to M)$ assigns $F$ and to the cone point assigns $M$ is a colimit diagram. But the category $\operatorname{Mod}_A^{\operatorname{ff}}(\operatorname{Ab})_{/M}$ is additive, so~(2) is exactly the statement that it is filtered. Hence,~(2) implies~(3), and it is clear that~(3) implies~(1).
\end{proof}

We also give an equational criterion for a map of Dirac rings $\phi \colon A \to B$ to be faithfully flat.

\begin{addendum}
\label{proposition:criterion_for_faithful_flatness}
Suppose that $\phi \colon A \to B$ is a map of Dirac rings which is both flat and a monomorphism. The following are equivalent:
\begin{enumerate}
\item[{\rm (1)}]The map $\phi \colon A \to B$ is faithfully flat.
\item[{\rm (2)}]The graded $A$-module $B/A$ given by the cokernel of $\phi \colon A \to B$ is flat.
\item[{\rm (3)}]Every map of graded $A$-modules $h \colon M \to B/A$ from a finitely presented graded $A$-module factors through the canonical projection $q \colon B \to B/A$.
\item[{\rm (4)}]Given a solution $(y_k)_{1 \leq k \leq n}$, consisting of homogeneous elements of $B$ to a system of linear equations
$$\textstyle{ \sum_{1 \leq k \leq n} y_kc_{ki} = d_i \hspace{12mm}
      (1 \leq i \leq m) }$$
with $c_{ki}$ and $d_i$ homogeneous elements of $A$, there exists a solution $(x_k)_{1 \leq k \leq n}$ consisting of homogeneous elements of $A$.
\end{enumerate}
\end{addendum}

\begin{proof}The statements~(1) and~(2) are equivalent to the statements that for every graded $A$-module $N$, the boundary map $\partial$ and the left-hand term $\operatorname{Tor}_1^A(N,B/A)$ in the exact sequence of graded $A$-modules
$$\xymatrix@C=7mm{
{ 0 } \ar[r] &
{ \operatorname{Tor}_1^A(N,B/A) } \ar[r]^-{\partial} &
{ N } \ar[r]^-{\phi} &
{ N \otimes_AB } \ar[r]^-{q} &
{ N \otimes_A(B/A) } \ar[r] &
{ 0 } \cr
}$$
are zero, respectively. So the statements~(1) and~(2) are equivalent.

To prove that~(2) and~(3) are equivalent, we employ \cref{thm:lazard}. Suppose first that~(2) holds. Given a map $h \colon M \to B/A$ of graded $A$-modules from a finitely presented graded $A$-module, it admits a factorization
$$\xymatrix{
{ M } \ar[r]^-{f} &
{ F } \ar[r]^-{\bar{g}} &
{ B/A } \cr
}$$
with $F$ a finitely generated free graded $A$-module, and since $F$ is projective, there exists $g \colon F \to B$ such that $\bar{g} = qg$. It follows that $h = qgf$, so~(3) holds.

Conversely, suppose that~(3) holds. We let $h \colon M \to B/A$ be a map from a finitely presented graded $A$-module and wish to show that it admits a factorization $h = \bar{g}f$ as above. By assumption, the map $h$ admits a factorization
$$\xymatrix{
{ M } \ar[r]^-{k} &
{ B } \ar[r]^-{q} &
{ B/A, } \cr
}$$
and since $B$ is flat, the map $k$ admits a further factorization
$$\xymatrix{
{ M } \ar[r]^-{f} &
{ F } \ar[r]^-{l} &
{ B } \cr
}$$
with $F$ finitely generated free. So $h = qlf$ is the desired factorization.

Suppose again that (3) holds. We can visualize the situation in (4) as a diagram of graded $A$-modules with exact rows
$$\xymatrix{
{} &
{ A^m } \ar[r]^-{c} \ar[d]^-{d} &
{ A^n } \ar[r]^-{p} \ar[d]^-{y} &
{ M } \ar[r] \ar[d]^-{\bar{y}} &
{ 0\phantom{.} } \cr
{ 0 } \ar[r] &
{ A } \ar[r]^-{\phi} &
{ B } \ar[r]^-{q} &
{ B/A } \ar[r] &
{ 0. } \cr
}$$
By~(3), there exists $s \colon M \to B$ such that $\bar{y} = qs$, and we wish to show that there exists $x \colon A^n \to A$ such that $d = xc$. But $q(y-sp) = 0$, and hence, there exists a unique $x \colon A^n \to A$ such that $\phi x = y - sp$. Finally, the calculation
$$\phi(d - xc) = yc - yc = 0$$
shows that $d = xc$, so (4) holds.

Conversely, suppose that~(4) holds. Given $\bar{y} \colon M \to B/A$, we produce a diagram of graded $A$-modules as above by using that $M$ is finitely presented and that $A^n$ is projective. By~(4), there exists $x \colon A^n \to A$ such that $xc = d$, and we wish to show that there exists $s \colon M \to B$ such that $qs = \bar{y}$. But $(y-\phi x)c = 0$, so there exists a unique $s \colon M \to B$ such that $y-\phi x = sp$. Therefore, the calculation
$$qsp = q(y-\phi x) = \bar{y}p$$
shows that $\bar{y} = qs$, so (3) holds.
\end{proof}

The module in \cref{example:evenly_generated_not_evenly_presented_module} that is evenly generated but not evenly presented is not flat. We now use \cref{thm:lazard} to show that this was no accident.

\begin{proposition}
\label{proposition:lazard}
Let $A$ be a Dirac ring. If a graded $A$-module $M$ is both flat and evenly generated, then is evenly presented.
\end{proposition}

\begin{proof}
With notation as in the proof of \cref{thm:lazard}, let $\operatorname{Mod}_A^{\operatorname{ff},\operatorname{ev}} \subset \operatorname{Mod}_A^{\operatorname{ff}}$ be the full subcategory spanned by the evenly finitely generated free graded $A$-modules. We claim that if $M$ is flat and evenly generated, then the inclusion
$$\xymatrix{
{ (\operatorname{Mod}_A^{\operatorname{ff},\operatorname{ev}})_{/M} } \ar[r] &
{ (\operatorname{Mod}_A^{\operatorname{ff}})_{/M} } \cr
}$$
is a $\varinjlim$-equivalence. Granting this, the corollary follows, since an evenly generated free graded $A$-module is evenly presented.

It remains to prove the claim. Since the categories in question are filtered, it will suffice to show that every map $h \colon F \to M$ with $F$ in $\operatorname{Mod}_A^{\operatorname{ff}}$ factors through some map $g \colon E \to M$ with $E$ in $\operatorname{Mod}_A^{\operatorname{ff},\operatorname{ev}}$. To this end, we choose families of homogeneous elements $y \colon T \to (F_d)_{d \in \mathbb{Z}}$ and $x \colon S \to (M_d)_{d \in \mathbb{Z}}$ such that $y$ is finite and a basis and such that $x$ is even and generates $M$. Given $j \in T$, we write
$$\textstyle{ h(y_j) = \sum_{i \in S}x_ia_{ij} }$$
with $a_{ij} \in A$. The subset $S' \subset S$ consisting of those $i$ such that $a_{ij} \neq 0$ for some $j \in T$ is finite. So we let $h \colon E \to M$ be the adjunct of $x|_{S'} \colon S' \to (M_d)_{d \in \mathbb{Z}}$, and define $f \colon T \to (E_d)_{d \in \mathbb{Z}}$ by $f(j) = \sum_{i \in S'}e_ia_{ij}$, where $e_i = \eta_{S'}(i) \in E$. It follows immediately from the definitions that $h \colon F \to M$ factors as the composition
$$\xymatrix{
{ F } \ar[r]^-{\tilde{f}} &
{ E } \ar[r]^-{g} &
{ M, } \cr
}$$
which proves the claim.
\end{proof}

\subsection{Properties of algebras}

Let $A$ be a Dirac ring. The functor that to a Dirac $A$-algebra $B$ assigns the graded set $(B_d)_{d \in \mathbb{Z}}$ admits a left adjoint functor that to a graded set $S$ assigns a Dirac $A$-algebra $A[S]$. A family $x \colon S \to (B_d)_{d \in \mathbb{Z}}$ of homogeneous elements in $B$ is algebraically independent (resp.\ generates $B$, resp.\ freely generates $B$) if the adjunct map $\tilde{x} \colon A[S] \to B$ is degreewise injective (resp.\ surjective, resp.\ bijective). A Dirac $A$-algebra $B$ is free if there exists a family of homogeneous elements that freely generates it.

\begin{example}
\label{ex:freealgebra}
If $S$ is a graded set with a single element $X$ of degree $e$, then we write also $A[X]$ for the free Dirac $A$-algebra $A[S]$. We can write every homogeneous element $f \in A[X]_d$ as a homogeneous polynomial
$$f = a_nX^n + \dots + a_1X + a_0$$
with $a_i \in A_{d-ie}$. But if $e$ is odd, then the coefficients $a_i \in A_{d-ie}$ with $i > 1$ are only well-defined modulo $2A_{d-ie}$. So as a graded $A$-module, $A[X]$ may not be free.
\end{example}

\begin{definition}
\label{def:finiteevenalgebra}
Let $A$ be a Dirac ring.
\begin{enumerate}
\item[(1)]A Dirac $A$-algebra $B$ is finitely generated if there exists a finite family of homogeneous elements that generates $B$.
\item[(2)]A Dirac $A$-algebra $B$ is finitely presented if there exists a finite family of homogeneous elements $x \colon S \to (B_d)_{d \in \mathbb{Z}}$ that generates $B$ and for which the kernel of the adjunct map $\tilde{x} \colon A[S] \to B$ is finitely generated.
\item[(3)]A Dirac $A$-algebra $B$ is finite if $B$ considered as a graded $A$-module is finitely generated.
\item[(4)]A Dirac $A$-algebra $B$ is evenly generated if there exists an even family of homogeneous elements that generates $B$.
\item[(5)]A Dirac $A$-algebra $B$ is evenly presented if there exists an even family of homogeneous elements $x \colon S \to (B_d)_{d \in \mathbb{Z}}$ that generates $B$ and for which the kernel of the adjunct map $\tilde{x} \colon A[S] \to B$ is evenly generated.
\item[(6)]A Dirac $A$-algebra $B$ is even if the map $B^{\operatorname{ev}} \otimes_{A^{\operatorname{ev}}}A \to B$ is an isomorphism, or equivalently, if $B$ considered as a graded $A$-module is evenly presented.
\end{enumerate}
\end{definition}
A basic result which we will prove below is a variant on the classical Hilbert's basis theorem: finitely generated Dirac algebras inherit the property of being noetherian in the following sense. 

\begin{definition}
\label{def:noetherian}
A Dirac ring $A$ is noetherian if every graded ideal $\mathfrak{a} \subset A$ is finitely generated.
\end{definition}

\begin{remark}
A Dirac ring $A$ is noetherian in the sense of \cref{def:noetherian} if and only if the Grothendieck abelian category of graded $A$-modules is locally noetherian in the sense of \cite[C.6.8]{luriesagtemp}, in which case, its noetherian objects of are precisely the finitely generated graded $A$-modules. 
\end{remark}

\begin{proposition}[Hilbert's basis theorem]
\label{proposition:hilbert_basis}
If a Dirac ring $A$ is noetherian, then so is every finitely generated $A$-algebra $B$.
\end{proposition}

\begin{proof}
By the obvious inductive argument, we may assume that $B$ is generated as Dirac $A$-algebra by a single homogeneous element $y \in B_e$. So we let
$$\mathfrak{b}_0 \subset \mathfrak{b}_1 \subset \dots \subset \mathfrak{b}_m \subset \dots$$
be an increasing sequence of graded ideals of $B$, and we proceed to show that it stabilizes. Given a homogeneous element $f \in B_d$, we will say that a family of homogeneous elements $(a_0,\dots,a_n)$ with $n \geq -1$ and $a_i \in A_{d-ie}$ is an expression of $f$ as a homogeneous polynomial in $y$ over $A$ if the equality
$$f = a_ny^n + \dots + a_1y + a_0$$
holds. In this case, we will refer to $n$ as the polynomial degree of the expression and to $a_n \in A_{d-ne}$ as its leading coefficient. We note that, in general, neither $n$ nor $a_n$ is uniquely determined by $f$, even if $B$ is a free Dirac $A$-algebra, but as we will see, this is irrelevant for the argument. We define $\mathfrak{a}_{m,n} \subset A$ to be the graded ideal, whose homogeneous elements are the leading coefficients of polynomial expressions of degree $n$ of homogeneous elements $f \in \mathfrak{b}_m$. Clearly, we have $\mathfrak{a}_{m,n} \subset \mathfrak{a}_{m+1,n}$ and $\mathfrak{a}_{m,n} \subset \mathfrak{a}_{m,n+1}$, so by the assumption that $A$ is noetherian, there are only finitely many distinct graded ideals among the $\mathfrak{a}_{m,n}$. Thus, there exists $r \geq 0$ such that if $m \geq r$, then $\mathfrak{a}_{m,n} = \mathfrak{a}_{r,n}$ for all $n \geq 0$. We claim that $\mathfrak{b}_m = \mathfrak{b}_r$ for $m \geq r$. Indeed, suppose that $f \in \mathfrak{b}_m$ is homogeneous and admits an expression as a homogeneous polynomial in $y$ over $A$ of polynomial degree $n$. In this situation, we can find a homogeneous element $g \in \mathfrak{b}_r$ such that $f - g \in \mathfrak{b}_m$ admits an expression as a homogeneous polynomial in $y$ over $A$ of polynomial degree $<n$. Hence, we conclude, by induction on $n$, that $f \in \mathfrak{b}_r$ as desired.
\end{proof}

We will say a map of Dirac rings $\phi \colon A \to B$ is flat (resp.\ faithfully flat) if $B$ is flat (resp.\ faithfully flat) considered as an $A$-module via $\phi$.

\begin{proposition}
\label{proposition:faithfully_flat_if_flat_and_surjective_on_prime_spectra}
A flat map of Dirac rings $\phi \colon A \to B$ is faithfully flat if and only if the induced map of Dirac-Zariski spaces $p \colon |Y| \to |X|$ is surjective.
\end{proposition}

\begin{proof}
This follows immediately from \cref{proposition:faithful_flatness_detected_by_tensors_with_residue_fields}. 
\end{proof}

\begin{proposition}
\label{proposition:local_flat_homomorphisms_are_faithfully_flat}
If a flat map $\phi \colon A \to B$ between local Dirac rings is local, then it is faithfully flat. In particular, it is injective. 
\end{proposition}

\begin{proof}
Since $\phi \colon A \to B$ is flat, it suffices by  \cref{proposition:faithful_flatness_detected_by_tensors_with_residue_fields} to prove that the unique maximal ideal $\mathfrak{m}_{A} \subseteq A$ is in the image of the induced map on Dirac-Zariski spaces. But $\phi^{-1}(\mathfrak{m}_{B}) = \mathfrak{m}_{A}$, because $\phi \colon A \to B$ is local.
\end{proof}

\begin{lemma}
\label{lemma:flatness_and_faithful_flatness_stable_under_base_change}
If $\phi \colon A \to B$ is any map of Dirac rings, then the functor
$$\xymatrix{
{ \operatorname{Mod}_A(\operatorname{Ab}) } \ar[r]^-{\phi^*} &
{ \operatorname{Mod}_B(\operatorname{Ab}) } \cr
}$$
preserves flatness and faithfulness of modules. If $\phi \colon A \to B$ is faithfully flat, then $\phi^*$ also reflects flatness and faithfulness of modules.
\end{lemma}

\begin{proof}
This follows immediately from the definitions and from the fact that the functor $\phi^*$ promotes to a symmetric monoidal functor.
\end{proof}

An important technique in Dirac geometry is passing to the subring of even degree elements. Unfortunately, this can fail many of the desirable properties of a map, as the following two instructive examples show. 

\begin{example}
\label{ex:finiteevenalgebra}
The finiteness conditions~(1)--(3) are generally not preserved by the functor that to $\phi \colon A \to B$ assigns $\phi^{\operatorname{ev}} \colon A^{\operatorname{ev}} \to B^{\operatorname{ev}}$. To wit, let $k$ be a Dirac field of characteristic $p \neq 2$, let $A = k[S]$ be a free graded $k$-algebra generated by an infinite graded set $S$ concentrated in degree $1$, and let $B = A[y]$ be a free graded $A$-algebra on a single generator of degree $1$. In this situtation, the map $\phi \colon A \to B$ is finite, but $\phi^{\operatorname{ev}} \colon A^{\operatorname{ev}} \to B^{\operatorname{ev}}$ is not even finitely generated. Moreover, in the factorization $A \to B' \to B$ of \cref{cor:evenintegralradicial}, neither map is finitely generated.
\end{example}

\begin{example}
\label{ex:maptoresiduefield}
Since $A^{\operatorname{ev}} \to A$ induces a homeomorphism of Dirac--Zariski spaces, one might expect that the canonical map $A \to k(\mathfrak{p})$, where $\mathfrak{p} \subset A$ is a graded prime ideal, is even. However, while the map is evenly generated, it may not be even, as the example $A = k[x]$ with $x$ homogeneous of degree $1$ shows.
\end{example}

We will now show that if we assume that $\phi \colon A \to B$ is even, then the pathological behavior illustrated by the above two examples does not occur. We begin with the following simple observation. 

\begin{lemma}
\label{lem:finiteevenalgebra}
If a map of Dirac rings $\phi \colon A \to B$ is finitely and evenly generated, then $\phi^{\operatorname{ev}} \colon A^{\operatorname{ev}} \to B^{\operatorname{ev}}$ is finitely generated.
\end{lemma}

\begin{proof}
By assumption, there exists families $x \colon S \to (B_d)_{d \in \mathbb{Z}}$ and $y \colon T \to (B_d)_{d \in \mathbb{Z}}$ with $x$ finite and $y$ even, both of which generate $B$ as an $A$-algebra. In particular, for every $i \in S$, we may write $x_i$ as a polynomial
$$\textstyle{ x_i = \sum_{n \in \mathbb{Z}_{\geq 0}}\sum_{j \colon n \to T}\sum_{k \colon n \to \mathbb{Z}_{\geq 0}} a_{i,j,k}y_{j_1}^{k_1} \dots y_{j_n}^{k_n}, }$$
where all but finitely many of the coefficients $a_{i,j,k}$ are zero. So let $T' \subset T$ be the graded subset consisting of those $j'$ for which there exists $i \in S$, $j \colon n \to T$, and $k \colon n \to \mathbb{Z}_{\geq 0}$ such that $a_{i,j,k} \neq 0$ and $j' = j_s$ for some $s \in n$. In this situation, the restriction $y' \colon T' \to (B_d)_{d \in \mathbb{Z}}$ of $y$ to $T' \subset T$ generates $B$ as an $A$-algebra and is finite and even. It follows that $y'$ also generates $B^{\operatorname{ev}}$ as an $A^{\operatorname{ev}}$-algebra.
\end{proof}

\begin{proposition}
\label{proposition:properties_of_even_extensions_reflected_by_base_change}
Let $\phi \colon A \to B$ be an even map of Dirac rings. The following properties hold for $\phi \colon A \to B$ if and only if they hold for $\phi^{\operatorname{ev}} \colon A^{\operatorname{ev}} \to B^{\operatorname{ev}}$.
\begin{enumerate}
\item[{\rm (1)}]Being finitely generated.
\item[{\rm (2)}]Being finitely presented.
\item[{\rm (3)}]Being flat.
\item[{\rm (4)}]Being a localization that inverts a single homogeneous element.
\item[{\rm (5)}]Being degreewise surjective.
\end{enumerate}
\end{proposition}

\begin{proof}
The properties in question are all stable under base change. In particular, if they hold for $\phi^{\operatorname{ev}}$, then they hold for $\phi$. We must prove the converse.

The case of~(1) follows from \cref{lem:finiteevenalgebra}, and in the case of~(2), we must show that $B^{\operatorname{ev}}$ is a compact object of the category of $A^{\operatorname{ev}}$-algebras. We claim that for any $A^{\operatorname{ev}}$-algebra $C$, there is a natural bijection
$$\operatorname{Map}_{A^{\operatorname{ev}}}(B^{\operatorname{ev}},C) \simeq \operatorname{Map}_A(B,C^{\operatorname{ev}} \otimes_{A^{\operatorname{ev}}}A).$$
Granting this, the statement in the case of~(2) follows. Indeed, if $C \simeq \varinjlim_i C_i$ is a filtered colimit of $A^{\operatorname{ev}}$-algebras, then $C^{\operatorname{ev}} \simeq \varinjlim_iC_i^{\operatorname{ev}}$, since filtered colimits are calculated degreewise, so the claim identifies the canonical map
$$\xymatrix{
{ \operatorname{Map}_{A^{\operatorname{ev}}}(B^{\operatorname{ev}},C) } \ar[r] &
{ \varinjlim_i \operatorname{Map}_{A^{\operatorname{ev}}}(B^{\operatorname{ev}},C_i) } \cr
}$$
that we would like to show is a bijection with the canonical map
$$\xymatrix{
{ \operatorname{Map}_A(B,C^{\operatorname{ev}}) } \ar[r] &
{ \varinjlim_i \operatorname{Map}_{A^{\operatorname{ev}}}(B^{\operatorname{ev}},C_i^{\operatorname{ev}} \otimes_{A^{\operatorname{ev}}}A), } \cr
}$$
which is a bijection by the assumption that $\phi \colon A \to B$ is finitely presented. It remains to prove the claim. First, the canonical map
$$\xymatrix{
{ C^{\operatorname{ev}} \simeq C^{\operatorname{ev}} \otimes_{A^{\operatorname{ev}}}A^{\operatorname{ev}} } \ar[r] &
{ (C^{\operatorname{ev}} \otimes_{A^{\operatorname{ev}}}A)^{\operatorname{ev}} } \cr
}$$
is an isomorphism, by the projection formula in \cref{lem:evenadjunction}. Therefore, by using the adjunction in loc.\ cit.\, we obtain a natural bijection
$$\xymatrix{
{ \operatorname{Map}_{A^{\operatorname{ev}}}(B^{\operatorname{ev}},C) } \ar[r] &
{ \operatorname{Map}_{A^{\operatorname{ev}}}(B^{\operatorname{ev}},C^{\operatorname{ev}} \otimes_{A^{\operatorname{ev}}}A). } \cr
}$$
Finally, since $f \colon A \to B$ is even, the adjunction between extension of scalars and restriction of scalars gives a natural bijection
$$\xymatrix{
{ \operatorname{Map}_A(B,C^{\operatorname{ev}} \otimes_{A^{\operatorname{ev}}}A) } \ar[r] &
{ \operatorname{Map}_{A^{\operatorname{ev}}}(B^{\operatorname{ev}},C^{\operatorname{ev}} \otimes_{A^{\operatorname{ev}}}A), } \cr
}$$
which proves the claim.

In the case of~(3), we must show that for every graded $A^{\operatorname{ev}}$-module $M$, the graded $A^{\operatorname{ev}}$-module $\pi_j(M \otimes_{A^{\operatorname{ev}}}^LB^{\operatorname{ev}})$ vanishes for $j > 0$. We may further assume that $M$ is even, since every graded $A^{\operatorname{ev}}$-module decomposes as the sum of its even part and its odd part. Now, since $\phi \colon A \to B$ is flat, the spectral sequence
$$E_{i,j}^2 = \pi_i(\pi_j(M \otimes_{A^{\operatorname{ev}}}^LA) \otimes_A^LB) \Rightarrow \pi_{i+j}(M \otimes_{A^{\operatorname{ev}}}^LB)$$
collapses, so the edge homomorphism
$$\xymatrix{
{ \pi_j(M \otimes_{A^{\operatorname{ev}}}^LA) \otimes_AB } \ar[r] &
{ \pi_j(M \otimes_{A^{\operatorname{ev}}}^LB) } \cr
}$$
is an isomorphism, and since $\phi \colon A \to B$ is even, also the base change map
$$\xymatrix{
{ \pi_j(M \otimes_{A^{\operatorname{ev}}}^LA) \otimes_{A^{\operatorname{ev}}}B^{\operatorname{ev}} } \ar[r] &
{ \pi_j(M \otimes_{A^{\operatorname{ev}}}^LA) \otimes_AB } \cr
}$$
is an isomorphism. The composite isomorphism restricts to an isomorphism
$$\xymatrix{
{ \pi_j(M \otimes_{A^{\operatorname{ev}}}^LA^{\operatorname{ev}}) \otimes_{A^{\operatorname{ev}}} B^{\operatorname{ev}} } \ar[r] &
{ \pi_j(M \otimes_{A^{\operatorname{ev}}}^LB^{\operatorname{ev}}) } \cr
}$$
of the even parts, which shows that $\phi^{\operatorname{ev}} \colon A^{\operatorname{ev}} \to B^{\operatorname{ev}}$ is flat.

Finally,~(4) follows from the fact that we can restrict attention to localizations at homogeneous elements of even degree, and~(5) is an immediate consequence of the fact that both extension of scalars and passing to even parts both preserve degreewise surjections.
\end{proof}

\subsection{Faithfully flat descent}

In this section, we prove faithfully flat descent for module categories of Dirac rings. The proof is analogous to the one for commutative rings, but to emphasize the similarity, we give an abstract argument, which applies to presentably symmetric monoidal abelian categories in general.

\begin{proposition}[Grothendieck]
\label{prop:flat_descent_for_abelian_categories_of_modules}
If $\phi \colon A \to B$ is a faithfully flat map of Dirac rings, then the diagram $\operatorname{Mod}_{B^{\otimes_A[-]}}(\operatorname{Ab}) \colon \Delta_+ \to \operatorname{LPr}$ is a limit diagram.
\end{proposition}

\begin{proof}
We have cocartesian diagrams of Dirac rings
$$\xymatrix{
{ A } \ar[r]^-{\phi_n} \ar[d]^-{\phi} &
{ B^{\otimes_A[n]} } \ar[d]^-{\phi'} \cr
{ B } \ar[r]^-{\phi_n'} &
{ (B \otimes_A B)^{\otimes_B[n]} } \cr
}$$
and we will use the fact that the diagrams of categories of graded modules
$$\xymatrix@C+8mm{
{ \operatorname{Mod}_A(\operatorname{Ab}) } \ar[r]^-{\phi_n^*} \ar[d]^-{\phi^*} &
{ \operatorname{Mod}_{B \otimes_A[n]}(\operatorname{Ab}) } \ar[d]^-{\phi'{}^*} \cr
{ \operatorname{Mod}_B(\operatorname{Ab}) } \ar[r]^-{\phi_n'{}^{\!\!*}} &
{ \operatorname{Mod}_{(B \otimes_AB)^{\otimes_B[n]}}(\operatorname{Ab}) } \cr
}$$
are right adjointable in the sense that the composite natural transformations
$$\xymatrix{
{ \phi^*\phi_{n*} } \ar[r] &
{ \phi_n'{}_*\phi_n'{}^{\!\!*}\phi^*\phi_{n*} \simeq \phi_n'{}_*\phi'{}^*\phi_n^*\phi_{n*} } \ar[r] &
{ \phi_n'{}_*\phi'{}^* } \cr
}$$
are natural isomorphisms. The map $\phi' \colon B \to B \otimes_AB$ is again a faithfully flat map of Dirac rings, and $\operatorname{Mod}_{(B \otimes_AB)^{\otimes_B[-]}}(\operatorname{Ab}) \colon \Delta_+ \to \operatorname{LPr}$ is a split cosimplicial diagram. Therefore, it is a limit diagram by~\cite[Lemma~6.1.3.16]{luriehtt}. We will show that this implies that also $\operatorname{Mod}_{B^{\otimes_A[-]}}(\operatorname{Ab}) \colon \Delta_+ \to \operatorname{LPr}$ is a limit diagram.

We recall from \cite[Proposition~5.5.3.13]{luriehtt} that the subcategory $\operatorname{LPr} \subset \widehat{\operatorname{Cat}}_{\infty}$ is closed under limits. In particular, the augmentation
$$\xymatrix{
{ \operatorname{Mod}_A(\operatorname{Ab}) } \ar[r]^-{F} &
{ \varprojlim_{[n] \in \Delta} \operatorname{Mod}_{B^{\otimes_A[n]}}(\operatorname{Ab}), } \cr
}$$
which we wish to prove is an equivalence, is a left adjoint. Its right adjoint
$$\xymatrix{
{ \varprojlim_{[n] \in \Delta} \operatorname{Mod}_{B^{\otimes_A[-]}}(\operatorname{Ab}) } \ar[r]^-{G} &
{ \operatorname{Mod}_A(\operatorname{Ab}) } \cr
}$$
is given, informally, by the limit
$$\textstyle{ G((M_n)_{[n] \in \Delta}) \simeq \varprojlim_{[n] \in \Delta} \phi_{n*}(M_n). }$$
This is the limit of a cosimplicial diagram in the $1$-category $\operatorname{Mod}_A(\operatorname{Ab})$. We recall the general fact, which we prove in \cref{prop:totalization}, that for any cosimplicial diagram $X \colon \Delta \to \mathcal{C}$ in a $1$-category, the restriction
$$\xymatrix{
{ \varprojlim_{\Delta} X } \ar[r] &
{ \varprojlim_{\Delta_{\leq 1}} X|_{\Delta_{\leq 1}} } \cr
}$$
along the inclusion $\Delta_{\leq 1} \to \Delta$ is an isomorphism. In particular, if a functor between $1$-categories preserves finite limits, then it preserves limits of cosimplicial diagrams. We must prove that the unit map
$$\xymatrix{
{ N } \ar[r]^-{\eta} &
{ \varprojlim_{[n] \in \Delta} \phi_{n*}\phi_n^*(N) } \cr
}$$
and the counit map
$$\xymatrix{
{ (\phi_m^*(\varprojlim_{[n] \in \Delta} \phi_{n*}(M_n))_{[m] \in \Delta}) } \ar[r]^-{\epsilon} &
{ (M_m)_{[m] \in \Delta} } \cr
}$$
are isomorphisms. Moreover, it follows from \cref{lem:conservative} that the projection
$$\xymatrix{
{ \varprojlim_{[n] \in \Delta} \operatorname{Mod}_{B^{\otimes_A[n]}}(\operatorname{Ab}) } \ar[r] &
{ \operatorname{Mod}_B(\operatorname{Ab}) } \cr
}$$
is conservative, so the counit is an isomorphism if and only if its image
$$\xymatrix{
{ \phi_0^*(\varprojlim_{[n] \in \Delta} \phi_{n*}(M_n)) } \ar[r]^-{\epsilon_0} &
{ M_0 } \cr
}$$
by the projection is so. Now, since $\phi \colon A \to B$ and $\phi' \colon B \to B \otimes_AB$ are flat and faithful, the extension of scalars functors $\phi^*$ and $\phi'{}^*$ are exact and faithful, and hence, conservative. Therefore, it suffices to show that $\phi^*(\eta)$ and $\phi'{}^*(\epsilon_0)$ are isomorphisms. But the composition of $\phi^*(\eta)$ and the isomorphisms
$$\begin{xy}
(0,8)*+{ \phi^*(\varprojlim_{[n] \in \Delta}\phi_{n*}\phi_n^*(N)) }="1";
(44,8)*+{ \varprojlim_{[n] \in \Delta}\phi^*\phi_{n*}\phi_n^*(N) }="2";
(-17,0)*+{}="3";
(11,0)*+{ \varprojlim_{[n] \in \Delta}\phi_n'{}_*\phi'{}^*\phi_n^*(N) }="4";
(55,0)*+{ \varprojlim_{[n] \in \Delta}\phi_n'{}_*\phi_n'{}^*\phi^*(N) }="5";
{ \ar "2";"1";};
{ \ar "4";"3";};
{ \ar "5";"4"; }
\end{xy}$$
is equal to the unit map
$$\xymatrix{
{ \phi^*(N) } \ar[r]^-{\eta'} &
{ \varprojlim_{[n] \in \Delta}\phi_n'{}_*\phi_n'{}^*(\phi_*(N)), } \cr
}$$
which we already know is an isomorphism. Similarly, the map $\phi'{}^*(\epsilon_0)$ is equal to the composition of the isomorphisms
$$\begin{xy}
(0,8)*+{ \phi'{}^*\phi_0^*(\varprojlim_{[n] \in \Delta} \phi_{n*}(M_n)) }="1";
(48,8)*+{ \phi_0'{}^*\phi^*(\varprojlim_{[n] \in \Delta} \phi_{n*}(M_n)) }="2";
(-19,0)*+{}="3";
(11,0)*+{ \phi_0'{}^*(\varprojlim_{[n] \in \Delta} \phi^*\phi_{n*}(M_n)) }="4";
(60,0)*+{ \phi_0'{}^*(\varprojlim_{[n] \in \Delta} \phi_n'{}_*\phi'{}^*(M_n)) }="5";
{ \ar "2";"1";};
{ \ar "4";"3";};
{ \ar "5";"4"; }
\end{xy}$$
and the counit map
$$\xymatrix{
{ \phi_0'{}^*(\varprojlim_{[n] \in \Delta} \phi_n'{}_*\phi'{}^*(M_n)) } \ar[r]^-{\epsilon_0'} &
{ \phi'{}^*(M_0), } \cr
}$$
which we know is an isomorphism. This completes the proof.
\end{proof}

\begin{addendum}
\label{addendum:faithfully_flat_descent_for_modules}If $\phi \colon A \to B$ is a faithfully flat map of Dirac rings, then
$$\xymatrix@C+=12mm{
{ \Delta_+ } \ar[r]^-{B^{\otimes_A[-]}} &
{ \operatorname{CAlg}(\operatorname{Ab}) } \cr
}$$
is a limit diagram of Dirac rings.
\end{addendum}

\begin{proof}The restriction of scalars along the unique map $\pi \colon \mathbb{Z} \to A$ promotes to a lax symmetric monoidal functor, and being a right adjoint, it preserves limits. Moreover, the diagram in question factors as the composition
$$\xymatrix@C+4mm{
{ \Delta_+ } \ar[r]^-{B^{\otimes_A[-]}} &
{ \operatorname{CAlg}(\operatorname{Mod}_A(\operatorname{Ab})) } \ar[r]^-{\pi_*} &
{ \operatorname{CAlg}(\operatorname{Ab}), }
}$$
so it suffices to prove that the left-hand functor is a limit diagram. But the forgetful functor from $\operatorname{CAlg}(\operatorname{Mod}_A(\operatorname{Ab}))$ to $\operatorname{Mod}_A(\operatorname{Ab})$ creates limits, so it suffices to show that the underlying diagram of graded $A$-modules is a limit diagram, or equivalently, that the comparison map $\eta \colon A \to \varprojlim_{[n] \in \Delta} B^{\otimes_A[n]}$ is an isomorphism. This map is the unit map of the adjunction appearing in the proof of \cref{prop:flat_descent_for_abelian_categories_of_modules}, where we proved that it is an isomorphism. 
\end{proof}

It is often also useful to know that we have descent of Dirac-Zariski spaces. In the following, let us denote by $X_{\operatorname{Zar}}$ the poset of the open subsets of the Dirac-Zariski space $|X|$ of a Dirac ring $A$.

\begin{proposition}
\label{proposition:faithfully_flat_maps_of_dirac_rings_gives_equalizer_of_open_sets}
Let $\phi \colon A \to B$ be a faithfully flat map of Dirac rings, and let $p \colon Y \to X$ be the induced map of Dirac spectra. The diagram
$$\xymatrix@C+=18mm{
{ \Delta_+ } \ar[r]^-{(Y^{\times_X[-]})_{\operatorname{Zar}}} &
{ \operatorname{Set} } \cr
}$$
of partially ordered sets of open subsets in the Dirac--Zariski space is a limit diagram.
\end{proposition}

\begin{proof}
The argument is analogous to \cite[1.6.2.2]{luriesag}, but for completeness, we include it here. Because sets form an ordinary category, it follows from \cref{prop:totalization} that the restriction map
$$\xymatrix{
{ \varprojlim_{[n] \in \Delta} (Y^{\times_X[n]})_{\operatorname{Zar}} } \ar[r] &
{ \varprojlim_{[n] \in \Delta_{\leq 1}} (Y^{\times_X[n]})_{\operatorname{Zar}} } \cr
}$$
is an isomorphism. Hence, it suffices to show that the diagram
$$\xymatrix{
{ X_{\operatorname{Zar}} } \ar[r]^{p^{-1}} &
{ Y_{\operatorname{Zar}} } \ar@<.7ex>[r]^-{p_1^{-1}} \ar@<-.7ex>[r]_-{p_2^{-1}} &
{ (Y \times_XY)_{\operatorname{Zar}} }
}$$
is an equalizer. The left-hand map is injective as a consequence of  \cref{proposition:faithfully_flat_if_flat_and_surjective_on_prime_spectra}. So it suffices to show that if $V \subseteq |Y|$ is an open subset such that $p_1^{-1}(V) = p_2^{-1}(V)$, then $V = p^{-1}(U)$ for some open subset $U \subset |X|$. We let $J \subset B$ be the radical graded ideal with $V = |Y| \smallsetminus V(J)$ and define $I = \phi^{-1}(J) \subset A$. We claim that $J \subset B$ is the radical of the ideal $\phi(I)B \subset B$, so that $U = |X| \smallsetminus V(I)$ will work. We write $A' = A/I$ and $B' = B/\phi(I)B = B \otimes_AA'$. The claim is equivalent to the statement that the image of any $y \in J$ by the canonical map $B \to B'$ is nilpotent. We also write $J' \subset B'$ for the image of $J \subset B$. The definition of $I \subset A$ implies that the map $\bar\phi \colon A'  \to B'/J'$ is injective, which will be important later.

Now, the graded ideals $J' \otimes_{A'}B' \subset B' \otimes_{A'}B'$ and $B' \otimes_{A'}J' \subset B' \otimes_{A'}B'$ have the same radical, because of our assumption that $p_1^{-1}(V) = p_2^{-1}(V)$. Hence, for every $y' \in J'$, some power of $y' \otimes 1$ maps to zero in $B' \otimes_{A'}(B'/J')$. But $y' \otimes 1$ is the image of $y' \in B'$ under the injective map
$$\xymatrix{
{ B' \simeq B' \otimes_{A'}A' } \ar[r] &
{ B' \otimes_{A'} (B'/J'), }
}$$
so we conclude that the image $y' \in B'$ of $j \in J$ is nilpotent, as desired.
\end{proof}

\begin{remark}
In \cref{proposition:faithfully_flat_maps_of_dirac_rings_gives_equalizer_of_open_sets}, it is possible to replace the category $\operatorname{Set}$ by the $\infty$-category $\operatorname{Cat}_{\infty}$, where we now consider each poset as a category in the usual way. That is, this diagram of sets promotes to a limit diagram of $\infty$-categories.
\end{remark}

\subsection{Integrality}\label{section:integrality_and_zariski's_main_theorem}

We revisit integrality in more detail in order to give a proof of Zariski's main theorem for Dirac rings. Suppose that $A[X]$ is a free Dirac $A$-algebra on a generator $X$ of degree $e$. We say that a homogeneous element $f(X) \in A[X]$ of the form
$$f(X) = X^n + a_{n-1}X^{n-1} + \dots + a_1X + a_0$$
is a monic homogeneous polynomial of degree $n$ over $A$. We remark that if $e$ is odd, then the terms $a_iX^i$ with $i > 1$ are of order $2$. Let $\phi \colon A \to B$ is a Dirac algebra. Given a homogeneous element $y \in B$ of degree $e$, there is a unique map of Dirac $A$-algebras $\operatorname{ev}_y \colon A[X] \to B$ that to $X$ assigns $y$. We write $f(y)$ for the image of $f(X)$ by this map and say that $y$ is a root in $f(X)$ if $f(y) = 0$.

\begin{definition}\label{def:integral}Let $\phi \colon A \to B$ be a Dirac algebra. A homogeneous element $y \in B$ is integral if it is a root in a monic homogeneous polynomial $f(X) \in A[X]$.
\end{definition}

We next give the usual equivalent characterization of integrality. We include the proof here, because most standard texts make use of Cayley--Hamilton, which is not available in the Dirac setting.

\begin{proposition}\label{prop:integral}Let $\phi \colon A \to B$ be a Dirac algebra, and let $y \in B$ be a homogeneous element. The following are equivalent.
\begin{enumerate}
\item[{\rm (1)}]The homogeneous element $y$ is integral.
\item[{\rm (2)}]The sub-$A$-algebra $A[y] \subset B$ generated by $y \in B$ is finite.
\item[{\rm (3)}]There exists a finite sub-$A$-algebra $B' \subset B$ such that $y \in B'$.
\end{enumerate}
\end{proposition}

\begin{proof}We first prove that~(1) and~(2) are equivalent. If $y \in B$ is a root in the monic homogeneous polynomial
$$f(X) = X^n + a_{n-1}X^{n-1} + \dots + a_1X + a_0,$$
then $(y^i)_{0 \leq i < n}$ generates $A[y]$ as a graded $A$-module. This shows that $A \to A[y]$ is finite. Conversely, suppose that $A \to A[y]$ is finite. We claim that for some $n \geq 1$, $(y^i)_{0 \leq i < n}$ generates $A[y]$ as a graded $A$-module. For $(y^i)_{i \geq 0}$ generates $A[y]$, so if some finite family $(z_1,\dots,z_m)$ generates $A[y]$, then each $z_j$ can be written as a linear combination of $(y^i)_{i \geq 0}$, and hence, of $(y^i)_{0 \leq i < n}$ for some $n \geq 1$. We now write $y^n$ as a linear combination of $(y^i)_{0 \leq i < n}$ to obtain a monic homogeneous polynomial $f(X) \in A[X]$ with $y$ as root. This shows that $y$ is integral.

It is clear that~(2) implies~(3), so we assume~(3) and prove~(1). We may assume that $B' = B$. So we let $y_1 = y$ and choose $y_2,\dots,y_m \in B$ homogeneous such that the family $(y_1,\dots,y_m)$ generates $B$ as a graded $A$-module. We can find homogeneous elements $a_{i,j,k} \in A$ such that for all $1 \leq i,j \leq m$,
$$\textstyle{ y_iy_j = \sum_{1 \leq k \leq m} a_{i,j,k}y_k. }$$
Let $A' \subset A$ be the Dirac subring generated by the family $(a_{i,j,k})_{1 \leq i,j,k \leq m}$. It is noetherian by  \cref{proposition:hilbert_basis} and by the fact that $\mathbb{Z}$ is noetherian. Now, the sub-$A'$-module $B' \subset B$ generated by $(y_1,\dots,y_m)$ is a finite Dirac $A'$-algebra that contains $y$, and since $A'$ is noetherian, so is the sub-$A'$-algebra $A'[y]$. Therefore, the first part of the proof shows that $y \in B'$ is integral over $A'$. This shows in particular that $y \in B$ is integral over $A$. This completes the proof.
\end{proof}

\begin{corollary}\label{cor:integral}
If $\phi \colon A \to B$ is any Dirac algebra, then the graded sub-$A$-module $B' \subset B$ spanned by the family of integral homogeneous elements $y \in B$ is a Dirac subalgebra.
\end{corollary}

\begin{proof}Suppose that $y_1,y_2 \in B$ are integral homogeneous elements. It follows from \cref{prop:integral} that the Dirac subalgebra $A[y_{1}] \subset B$ is finite over $A$, and that the Dirac subalgebra $A[y_{1}, y_{2}] \subset B$ is finite over $A[y_{1}]$, so $A[y_{1}, y_{2}] \subset B$ is finite over $A$. We conclude from \cref{prop:integral} that $y_1 \cdot y_2 \in B$ is integral. If $y_1$ and $y_2$ have the same degree, then we also conclude that $y_1 + y_2 \in B$ is integral.
\end{proof}

\begin{definition}
\label{def:integralalgebra}
Let $\phi \colon A \to B$ be a Dirac algebra. The integral closure of $A$ in $B$ is the Dirac subalgebra $B' \subset B$ spanned by the family of integral homogeneous elements $y \in B$. The Dirac algebra $\phi \colon A \to B$ is integral if $B' = B$.
\end{definition}

\begin{corollary}
\label{cor:integralalgebra}
Let $\phi \colon A \to B$ be a Dirac algebra, and let $B' \subset B$ be a Dirac subalgebra that is integral over $A$. If a homogeneous element $y \in B$ is integral over $B'$, then it is also integral over $A$.
\end{corollary}

\begin{proof}
If a homogeneous element $y \in B$ is integral over $B'$, then it is also integral over a finitely generated Dirac subalgebra $B'' \subset B'$. But since $B'$ is integral over $A$, an inductive argument based on \cref{prop:integral} shows that $B''$ is finite over $A$. It follows that $B''[y]$ is finite over $A$, since $B''[y]$ is finite over $B''$ and $B''$ is finite over $A$. So \cref{prop:integral} shows that $y$ is integral over $A$.
\end{proof}

\begin{remark}
\label{rem:integral}
Let $B$ be a Dirac $A$-algebra. It follows from \cref{prop:integral} that a homogeneous element $y \in B$ is integral if and only if $y^2 \in B$ is integral. Moreover, the homogeneous element $y^2 \in B$, which is  even, is integral over $A$ if and only if it is integral over $A^{\operatorname{ev}}$. So, in this sense, being integral is an even property.
\end{remark}

A non-trivial fact of commutative algebra is that in the even context Dirac integrality is closely related to integrality for the underlying commutative ring. 

\begin{proposition}
\label{proposition:integrality_compatible_with_underlying_rings}
Let $\phi \colon A \to B$ be a Dirac algebra, where both $A$ and $B$ are even, and let $\smash{\widetilde{\phi}} \colon \smash{\widetilde{A}} \to \smash{\widetilde{B}}$ be the  underlying algebra. An element $b \in \smash{\widetilde{B}}$ is integral over $\smash{\widetilde{A}}$ if and only if its homogeneous components $b_i \in B$ all are integral over $A$.
\end{proposition}

\begin{proof}
Indeed, this is a rephrasing of the compatibility of normalization with smooth base-change; see \cite[\href{https://stacks.math.columbia.edu/tag/03GG}{Lemma 03GG}]{stacks-project}.
\end{proof}

One useful property of integral maps is that the going-up theorem holds.

\begin{lemma}
\label{lemma:integral_homomorphisms_satisfy_going_up}
Let $\phi \colon A \to B$ be an integral Dirac algebra. Given  graded prime ideals $\mathfrak{p} \subset \mathfrak{p}' \subset A$ and $\mathfrak{q} \subset B$ such that $\phi^{-1}(\mathfrak{q}) = \mathfrak{p}$, there exists a graded prime ideal $\mathfrak{q} \subset \mathfrak{q}' \subset B$ such that $\phi^{-1}(\mathfrak{q}') = \mathfrak{p}'$.
\end{lemma}

\begin{proof}
The proof is analogous to~\cite[\href{https://stacks.math.columbia.edu/tag/00GU}{Tag 00GU}]{stacks-project}, but the argument is short, so we include it here. By replacing $A$ by $A/\mathfrak{p}$ and $B$ by $B/\mathfrak{q}$, it suffices to show that if $\phi \colon A \to B$ is integral and injective, then the induced map on Dirac--Zariski spaces is surjective. Thus, we show that for every graded prime ideal $\mathfrak{p} \subset A$, the inclusion $\mathfrak{p} B_{\mathfrak{p}} \neq B_{\mathfrak{p}}$. By further replacing $A$ and $B$ by the localizations at $\mathfrak{p} \subset A$, we can assume that $A$ is a local Dirac ring, and it suffices to show $\mathfrak{p} B \neq B$. So we assume that $\mathfrak{p}B = B$ and write $1 = \sum_{1 \leq i \leq m} x_i b_i$ with $x_i \in \mathfrak{p}$ and $b_i \in B$ homogeneous. The Dirac sub-$A$-algebra $B' \subset B$ generated by $(b_1,\dots,b_m)$ is finitely generated and integral, and hence, finite. Moreover, it satisfies $\mathfrak{p}B' = B'$, so Nakayama's lemma shows that $B'$ is zero, contradicting the assumption that $\phi \colon A \to B$ is injective.
\end{proof}

\begin{corollary}
\label{corollary:integral_extensions_are_closed}
If $\phi \colon A \to B$ is an integral Dirac algebra, then the induced map
$$\xymatrix{
{ | \operatorname{Spec}(B) | } \ar[r] &
{ | \operatorname{Spec}(A) | } \cr
}$$
of Dirac--Zariski spaces is closed.
\end{corollary}

\begin{proof}
This is a consequence of the going-up property established in \cref{lemma:integral_homomorphisms_satisfy_going_up}, see \cite[\href{https://stacks.math.columbia.edu/tag/00HZ}{Tag 00HZ}]{stacks-project}.
\end{proof}

An important example of an integral extension is the inclusion of the sub-Dirac ring of homogeneous elements of even degree, as we now prove.

\begin{proposition}
\label{proposition:inclusion_of_even_subring_integral_and_a_universal_homeomorphism}
Let $S$ be the Dirac--Zariski spectrum of a Dirac ring $A$. The unit map $\eta_S \colon S \to S^{\operatorname{ev}}$ is integral and a universal homeomorphism.
\end{proposition}

\begin{proof}
The inclusion $A^{\operatorname{ev}} \to A$ is integral, because every element of $A$ squares to an element of even degree, so the map $\eta_S$ is integral. It remains to show that it is a universal homeomorphism. We must show that the base-change $\eta_S' \colon T' \to T$ of $\eta_S$ along any map $f \colon T \to S^{\operatorname{ev}}$ is a homeomorphism of underlying spaces, and since this property is local on the base, we may assume that $T \simeq \operatorname{Spec}(B)$ is affine. Thus, we wish to show that the map of Dirac--Zariski spaces
$$\xymatrix@C=10mm{
{ | \operatorname{Spec}(A \otimes_{A^{\operatorname{ev}}}B) | } \ar[r]^-{|\eta_S'|} &
{ | \operatorname{Spec}(B) | } \cr
}$$
is a homeomorphism. Since $\eta_S$ is integral, so is $\eta_S'$, and therefore, we conclude from \cref{corollary:integral_extensions_are_closed} that the map $|\eta_S'|$ is closed. To see that it is surjective, we observe that the fiber over $\mathfrak{p} \in |\operatorname{Spec}(B)|$ is identified with $| \operatorname{Spec}(A \otimes_{A^{ev}} k(\mathfrak{p}))|$, which is non-empty, since the $k(\mathfrak{p})$-algebra 
$$A \otimes_{A^{\operatorname{ev}}} k(\mathfrak{p}) \simeq (A^{\operatorname{ev}} \oplus A^{\operatorname{odd}}) \otimes_{A^{\operatorname{ev}}} k(\mathfrak{p}) \simeq k(\mathfrak{p}) \oplus (A^{\operatorname{odd}} \otimes_{A^{\operatorname{ev}}} k(\mathfrak{p}))$$
is non-zero. Thus, it remains to show that the map $|\eta_S'|$ is injective, and by arguing as in \cite[\href{https://stacks.math.columbia.edu/tag/01S4}{Tag 01S4}]{stacks-project}, it suffices to check that every map $\psi \colon A^{\operatorname{ev}} \to K$ to a Dirac field $K$ admits at most one extension to a map $\psi' \colon A \to K$. If $K$ is even, then this is automatic, so by \cref{prop:classification_of_dirac_fields}, only the case $\operatorname{char}(K) = 2$ needs proof. We claim that if $f \in A$ is any homogeneous element, then the square
$$\psi'(f)^{2} = \psi'(f^{2}) = \psi(f^{2})$$
determines $\psi'(f)$. Indeed, since $\operatorname{char}(K) = 2$, so the equation $x = \psi(f^{2})$ has at most a one solution in $K$, because the difference between any two solutions is nilpotent. This proves the claim, and hence, the proposition.
\end{proof}

\subsection{Zariski's main theorem}

In this section, we prove the Dirac analogue of one of the fundamental theorems of commutative algebra, namely, Zariski's main theorem. Classically, is the fact that quasi-finite algebras of finite type can be well approximated by integral, or even finite, subalgebras. As \cref{example:non_quasi_finite_dirac_algebra_with_finite_spectrum} show, there exist finitely generated algebras over a Dirac field, which finite Dirac--Zariski space, but which are not finite as algebras. So in the Dirac context, quasi-finiteness is a more subtle phenomenon, which cannot be detected by the Dirac--Zariski space alone. Instead, we make the following definition. 

\begin{definition}
\label{def:quasifinite}
A Dirac algebra $\phi \colon A \to B$ is quasi-finite if it is finitely generated and if for every $x \in \lvert\operatorname{Spec}(A)\rvert$, the fiber
$$\xymatrix{
{ k(x) } \ar[r] &
{ B \otimes_Ak(x) } \cr
}$$
is a finite Dirac algebra.
\end{definition}

We will now study the interaction of quasi-finiteness and passing to even subrings. It follows from \cref{prop:evendiraczariskispace} that every graded prime ideal of $\mathfrak{q} \subset A^{\operatorname{ev}}$ is of the form $\mathfrak{q} = \mathfrak{p}^{\operatorname{ev}}$ for some graded prime ideal $\mathfrak{p} \subseteq A$. 

\begin{lemma}
\label{lem:tensoring_with_even_residue_fields}
Let $A$ be a Dirac ring, and let $\mathfrak{p} \subset A$ be a graded prime ideal. If $\phi \colon A \to C$ is an evenly generated algebra, then the canonical map
$$\xymatrix{
{ C^{\operatorname{ev}} \otimes_{A^{\operatorname{ev}}} k(\mathfrak{p}^{\operatorname{ev}}) } \ar[r] &
{ (C \otimes_A k(\mathfrak{p}))^{\operatorname{ev}} } \cr
}$$
is an isomorphism.
\end{lemma}

\begin{proof}
In the composition
$$\xymatrix{
{ C^{\operatorname{ev}} \otimes_{A^{\operatorname{ev}}} \mathfrak{p}^{\operatorname{ev}} } \ar[r] &
{ (C^{\operatorname{ev}} \otimes_{A^{\operatorname{ev}}} \mathfrak{p})^{\operatorname{ev}} } \ar[r] &
{ (C \otimes_A \mathfrak{p})^{\operatorname{ev}}, } \cr
}$$
the left-hand map is an isomorphism, by the projection formula, and the right-hand map is surjective, by the assumption that $\phi$ is evenly generated. It follows that the canonical map $(\mathfrak{p}^{\operatorname{ev}}C)^{\operatorname{ev}} \to (\mathfrak{p}C)^{\operatorname{ev}}$ is an isomorphism, and hence, so is
$$\xymatrix{
{ (C/\mathfrak{p}^{\operatorname{ev}}C)^{\operatorname{ev}} } \ar[r] &
{ (C/\mathfrak{p}C)^{\operatorname{ev}}. } \cr
}$$
The lemma follows by inverting all even homogeneous elements of $C$ that are not contained in $\mathfrak{p}$ and using \cref{prop:calculusoffractions}.
\end{proof}

\begin{remark}
\label{rem:tensoring_with_even_residue_fields}
In \cref{lem:tensoring_with_even_residue_fields}, the assumption that the algebra $\phi \colon A \to C$ be evenly generated is necessary, as the example $A = k[x]$ and $C = A[y]$ with $x$ and $y$ homogeneous of degree $1$ shows.
\end{remark}

\begin{corollary}
\label{cor:tensoring_with_even_residue_fields}Let $A$ be a Dirac ring, and let $\mathfrak{p} \subseteq A$ be a graded prime ideal. The canonical map $k(\mathfrak{p}^{\operatorname{ev}}) \to k(\mathfrak{p})^{\operatorname{ev}}$ is an isomorphism.
\end{corollary}

\begin{proof}
This follows from \cref{lem:tensoring_with_even_residue_fields} with $C = A$.
\end{proof}

\begin{proposition}
\label{proposition:quasi_finiteness_passes_to_even_extensions}
Suppose that $\phi \colon A \to B$ is an evenly generated Dirac algebra. If $\phi \colon A \to B$ is quasi-finite, then so is $\phi^{\operatorname{ev}} \colon A^{\operatorname{ev}} \to B^{\operatorname{ev}}$.
\end{proposition}

\begin{proof}
Since $\phi$ is both evenly generated and finitely generated, it is generated by finitely many even degree elements by \cref{lem:finiteevenalgebra}. These elements then also generate $B^{ev}$ over $A^{ev}$, so $\phi^{even}$ is finite type, and we also have to verify that the base-changes to residue fields are finite. 

A combination of \cref{lem:tensoring_with_even_residue_fields} and \cref{cor:tensoring_with_even_residue_fields} shows that, in the diagram
$$\xymatrix{
{ k(\mathfrak{p}^{\operatorname{ev}}) } \ar[r] \ar[d] &
{ k(\mathfrak{p})^{\operatorname{ev}} } \ar[r] \ar[d] &
{ k(\mathfrak{p})^{\operatorname{ev}} } \ar[d] \cr
{ B^{\operatorname{ev}} \otimes_{A^{\operatorname{ev}}} k(\mathfrak{p}^{\operatorname{ev}}) } \ar[r] &
{ B^{\operatorname{ev}} \otimes_{A^{\operatorname{ev}}} k(\mathfrak{p})^{\operatorname{ev}} } \ar[r] &
{ (B \otimes_A k(\mathfrak{p}))^{\operatorname{ev}}, } \cr
}$$
the horizontal inclusions all are isomorphisms. Moreover, the right-hand vertical map is finite, because $k(\mathfrak{p}) \to B \otimes_Ak(\mathfrak{p})$ is evenly generated and finite. Hence, the left-hand vertical map is finite, as we wanted to show.
\end{proof}

\begin{proposition}
\label{proposition:quasi_finiteness_passes_to_underlying_rings}
If $\phi \colon A \to B$ is a quasi-finite map between even Dirac rings, then the map of underlying rings $\smash{\widetilde{\phi}} \colon \smash{\widetilde{A}} \to \smash{\widetilde{B}}$ is quasi-finite.
\end{proposition}

\begin{proof}We let $p \colon |\smash{\widetilde{X}}| \to |X|$ be the map in \cref{prop:orbitspace} from the Zariski space of $\smash{\widetilde{A}}$ to the Dirac--Zariski space of $A$, and let $y \in |\smash{\widetilde{X}}|$ with image $x = p(y) \in |X|$. In this situation, the set-theoretic image of the $\mathbb{G}_m$-equivariant embedding
$$\xymatrix{
{ \operatorname{Spec}(\widetilde{k(x)}) } \ar[r]^-{i_x} &
{ \smash{\widetilde{X}} } \cr
}$$
is equal to the fiber $p^{-1}(x) \subset |\smash{\widetilde{X}}|$. In particular, we have a map of rings
$$\xymatrix{
{ \widetilde{k(x)} } \ar[r]^{i_y} &
{ k(y), } \cr
}$$
because $y \in p^{-1}(x)$. Hence, we may identify the map $k(y) \to \smash{\widetilde{B}} \otimes_{\smash{\widetilde{A}}}k(y)$, which we wish to prove is finite, with the base-change along $i_y$ of the map of underlying of rings associated with the map of Dirac rings $k(x) \to B \otimes_Ak(x)$, which by assumption is finite. This concludes the proof.
\end{proof}

\begin{theorem}[Zariski's Main Theorem]
\label{theorem:zariski_main}
Let $\phi \colon A \to B$ be a quasi-finite map of Dirac rings, and let $B' \subset B$ be the integral closure of $A$ in $B$. In this situation, the induced map of Dirac spectra
$$\xymatrix{
{ Y \simeq \operatorname{Spec}(B) } \ar[r]^-{j} &
{ Y' \simeq \operatorname{Spec}(B') } \cr
}$$
is an open immersion.
\end{theorem}

\begin{proof}We prove this in three steps:
\begin{enumerate}
    \item The case, where Dirac rings $A$ and $B$ both are even.
    \item The case, where the map of Dirac rings $\phi \colon A \to B$ is even.
    \item The general case.
\end{enumerate}
In case~(1), \cref{proposition:quasi_finiteness_passes_to_underlying_rings} shows that the map $\phi \colon \smash{\widetilde{A}} \to \smash{\widetilde{B}}$ of underlying rings is quasi-finite, and \cref{proposition:integrality_compatible_with_underlying_rings} shows that the subring
$$\smash{\widetilde{B}}' \subset \smash{\widetilde{B}}$$
is the integral closure of $\smash{\widetilde{A}}$ in $\smash{\widetilde{B}}$. By Zariski's main theorem~\cite[\href{https://stacks.math.columbia.edu/tag/00QB}{Tag 00QB}]{stacks-project}, the induced map of prime spectra is an open immersion, and therefore, we conclude from \cref{prop:open_embeddings_of_even_dirac_rings_detected_on_underlying_rings} that so is the map of Dirac spectra $j \colon Y \to Y'$.

In case~(2), the assumption that $\phi \colon A \to B$ is generated by a finite family of even elements implies that the same is true of the map $\phi^{\operatorname{ev}} \colon A^{\operatorname{ev}} \to B^{\operatorname{ev}}$ of even Dirac rings. Moreover, it follows from  \cref{proposition:quasi_finiteness_passes_to_even_extensions} that $\phi^{\operatorname{ev}}$ is quasi-finite, so that what we have already proved applies. Since the Dirac--Zariski space is insensitive to passing to the even part, we deduce from \cref{prop:evendiraczariskispace} that the continuous map $j \colon |Y| \to |Y'|$ is an open embedding.

Now, suppose that $f \in B'$ is such that $\smash{ |Y_f'| \subset |Y| \subset |Y'| }$. By replacing $f$ by its square, we can assume that $f$ is of even degree. We claim that the map $B_f' \to B_f$ is an isomorphism. Indeed, it is injective, because it is the localization of an injective map. And it is surjective, because $(B')_f^{\operatorname{ev}} \to B_f^{\operatorname{ev}}$ is an isomorphism, by what we have already proved, and because $\phi \colon A \to B$ is evenly generated. This ends the argument in case (2). 

In the general case, we factor $\phi \colon A \to B$ as a composition
$$\xymatrix{
{ A } \ar[r]^-{\psi} &
{ C } \ar[r]^-{\pi} &
{ B, } \cr
}$$
where the left-hand map is even and finitely generated, and where the right-hand map is finite and obtained by attaching square roots. We claim that the map $\psi$ is quasi-finite. Given a graded prime ideal $\mathfrak{p} \subset A$, the composition of the maps
$$\xymatrix{
{ k(\mathfrak{p}) } \ar[r]^-{\bar{\psi}} &
{ C \otimes_A k(\mathfrak{p}) } \ar[r]^-{\bar{\pi}} &
{ B \otimes_A k(\mathfrak{p}) } \cr
}$$
is finite, and we wish to show that $\bar{\psi}$ is finite. The map $\bar{\pi}$ induces a homeomorphism of Dirac--Zariski spaces, so every homogeneous element of $I = \ker(\bar{\pi})$ is nilpotent. But $I$ is a finitely generated ideal, because $C \otimes_Ak(\mathfrak{p})$ is finitely generated over a Dirac field, and hence, noetherian. It follows that $I$ is a nilpotent ideal, so an inductive argument using the $I$-adic filtration proves the claim. This shows that the map $\psi$ satisfies the hypotheses of~(2), so if we let $C' \subset C$ be the integral closure of $A$ in $C$, then the map induced map of Dirac spectra
$$\xymatrix{
{ V \simeq \operatorname{Spec}(C) } \ar[r]^-{h} &
{ V' \simeq \operatorname{Spec}(C') } \cr
}$$
is an open immersion.

Now, the map $\pi \colon C \to B$ is an $F$-isomorphism in the sense of Quillen. This means, we recall, that its kernel is nilpotent (in our case zero) and that its image contains some power of every homogeneous element in $B$. The map $\pi' \colon C' \to B'$ is an $F$-isomorphism as well. Indeed, it is injective, and if $f \in B'$, then some power of it belongs to $C' = C \cap B'$. Thus, in the diagram
$$\xymatrix@C+=11mm{
{ |Y| } \ar[r]^-{j} \ar[d] &
{ |Y'|\phantom{,} } \ar@<-.5ex>[d] \cr
{ |V| } \ar[r]^-{h} &
{ |V'|, } \cr
}$$
the vertical maps are homeomorphisms, and since $h$ is an open embedding, so is $j$.

Now, suppose that $f \in B'$ is such that $|Y_f'| \subset |Y| \subset |Y'|$. Replacing $f$ by some power thereof, if necessary, we can assume that $f \in C'$. We know from~(2) that the induced map $C_f' \to C_f$ is an isomorphism, and we claim that so is $B_f' \to B_f$. The map in question is injective, and to prove that it is also surjective, we must show that for every $g \in B$ homogeneous, there exists $k \geq 0$ such that $f^kg$ is integral. To this end, we choose $n \geq 1$ such that $g^n \in C$. Since $C_f' \to C_f$ is an isomorphism, there exists $k \geq 0$ such that $(f^kg)^n = f^{kn}g^n$ is integral. But then $f^kg$ is integral, as we wanted to prove. This completes the proof of Zariski's main theorem.
\end{proof}

\begin{corollary}
\label{corollary:finite_open_immersion_factorization_of_qf_morphism}
If $\phi \colon A \to B$ is a quasi-finite map of Dirac rings, then the induced map of Dirac prime spectra  $f \colon Y \to X$ admits a factorization as the composition
$$\xymatrix{
{ Y } \ar[r]^-{j_0} &
{ Y_0 } \ar[r]^-{p_0} &
{ X } \cr
}$$
of an open immersion $j_0$ and a finite map $p_0$.
\end{corollary}

\begin{proof}
This follows from \cref{theorem:zariski_main} as in \cite[\href{https://stacks.math.columbia.edu/tag/00QB}{Tag 00QB}]{stacks-project}, but let us include the argument here for completeness. Let $B' \subset B$ be the integral closure of $A$ in $B$. By \cref{theorem:zariski_main}, the induced map of Dirac spectra $j \colon Y \to Y'$ is an open immersion. Therefore, there exists a finite family $(|Y_{f_i}'|)_{i \in I}$ of distinguished open subsets of $|Y'|$ that covers $j(|Y|) \subset |Y'|$. In particular, the $A$-algebra $B_{f_i}'$ is finitely generated, because $B_{f_i}' \simeq B_{f_i}$. It follows that for each $i \in I$, we may choose a finite family $(g_{i,j})_{j \in J_i}$ of homogeneous elements that generates $B_{f_i}'$ as an $A_{f_i}$-algebra. We let $B_0 \subset B'$ to be the sub-$A$-algebra generated by all of the $f_i$ and $g_{i,j}$, and let
$$\xymatrix{
{ Y } \ar[r]^-{j_0} &
{ Y_0 } \ar[r]^-{p_0} &
{ X } \cr
}$$
be the maps of Dirac spectra induced by $\phi \colon A \to B_0 \subset B$. The map $p_0$ is finite by definition, and to prove that $j_0$ is an open immersion, it suffices to show that for every $i \in I$, the map $(B_0)_{f_i} \to B_{f_i}'$ induced by the canonical inclusion is an isomorphism. But the map is injective, since it is a filtered colimit of injective maps, and it is surjective, since its image contains a family that generates the target.
\end{proof}

\section{\'{E}tale maps}

This section is devoted to \'{e}tale maps of Dirac rings.

\subsection{Differentials} 
\label{subsection:module_of_kahler_differentials}
Let $k$ be a Dirac ring, let $\phi \colon k \to A$ be a Dirac algebra, and let $M$ be a graded $A$-module. We define the trivial square-zero extension of $A$ by $M$ to be the Dirac $k$-algebra with underlying graded $k$-module $A \oplus M$ and with multiplication given by
$$(a,x) \cdot (b,y) = (ab,ay+xb) = (ab,ay + (-1)^{\deg(b)\deg(x)}bx).$$
The projection $\pi \colon A \oplus M \to A$ is a map of Dirac $k$-algebras, and the pair $(A \oplus M,\pi)$ is an object of the slice category $\operatorname{CAlg}(\operatorname{Mod}_k)_{/A}$. The Dirac $k$-module structure of $M$ gives rise to an abelian group structure on $(A \oplus M,\pi)$, and conversely, one can show that every abelian group object in $\operatorname{CAlg}(\operatorname{Mod}_k)_{/A}$ is of this form.

\begin{definition}\label{def:firstdefinitionofdifferentials}
Let $k$ be a Dirac ring, let $\phi \colon k \to A$ be a Dirac algebra, and let $M$ be a graded $A$-module. A $k$-linear derivation of $A$ with values in $M$ is a section of Dirac $k$-algebra map $\pi \colon A \oplus M \to A$.
\end{definition}

Unwinding the definition, we see that $\sigma \colon A \to A \oplus M$ is a $k$-linear derivation of $A$ with values in $M$ if and only if $\sigma(a) = (a,d(a))$, where $d \colon A \to M$ is a map of graded $k$-modules that satisfies the Leibniz rule
$$d(ab) = d(a)b + ad(b) = (-1)^{\deg(a)\deg(b)}bd(a) + ad(b)$$
for all homogeneous $a,b \in A$. We stress that we treat $d$ as a spin $0$ operator. Below, we will abuse language and refer to both $\sigma \colon A \to A \oplus M$ and to $d \colon A \to M$ as a $k$-linear derivation of $A$ with values in $M$. We write
$$\operatorname{Der}_k(A,M)$$
for the set of $k$-linear derivations of $A$ with values in $M$. Up to unique isomorphism, there exists a unique $k$-linear derivation $d \colon A \to \Omega_{A/k}$ such that the map
$$\xymatrix{
{ \operatorname{Map}(\Omega_{A/k},M) } \ar[r] &
{ \operatorname{Der}_k(A,M) } \cr
}$$
given by composition with $d$ is a bijection for every graded $A$-module $M$.

\begin{definition}
\label{definition:module_of_kahler_differentials}If $\phi \colon k \to A$ is a Dirac algebra, then the $k$-linear derivation
$$\xymatrix{
{ A } \ar[r]^-{d} &
{ \Omega_{A/k} } \cr
}$$
is called the universal $k$-linear derivation of $A$, and the graded $A$-module $\Omega_{A/k}$ is called the module of differentials of $A$ over $k$.
\end{definition}

Concretely, if $I$ is the kernel of the multiplication $s^0 \colon A^{\otimes_k[1]} \to A^{\otimes_k[0]}$, then the map $d \colon A \to I/I^2$ defined by $df = d^0f - d^1f + I^2$ is a universal $k$-linear derivation.

\begin{example}
\label{example:kahler_differentials_of_a_free_dirac_algebra}
Let $k$ be a Dirac ring, $N$ a $k$-module, and $A = \operatorname{Sym}_{k}(N)$ the free Dirac $k$-algebra generated by $N$. We consider the following diagram.
$$\xymatrix{
{ \operatorname{Der}_k(A,M) } \ar[r] \ar[d] &
{ \operatorname{Map}(A,A \oplus M) } \ar[r] \ar[d]^-{\pi} &
{ \operatorname{Map}(N,A \oplus M) } \ar[r] \ar[d]^-{\pi} &
{ \operatorname{Map}(N,M) } \ar[d] \cr
{ \{\operatorname{id} \} } \ar[r] &
{ \operatorname{Map}(A,A) } \ar[r] &
{ \operatorname{Map}(N,A) } \ar[r] &
{ \{0\} } \cr
}$$
By definition, the left-hand square, whose right-hand column is formed by the respective sets of maps of Dirac $k$-algebras, is cartesian. The right-hand square, which is formed by the respective sets of maps of graded $k$-modules, and where the left-hand vertical map and the top horizontal map are induced by the respective projections, is also cartesian. Finally, the middle square, whose horizontal maps are the adjunction isomorphisms, commutes by naturality, and hence, is cartesian. So we conclude that the composition of the top horizontal maps is a bijection. This, in turn, implies that the $A$-linear map
$$\xymatrix{
{ A \otimes_kN } \ar[r] &
{ \Omega_{A/k} } \cr
}$$
given by the adjunct of the composite $k$-linear map
$$\xymatrix{
{ N } \ar[r]^-{\eta} &
{ A } \ar[r]^-{d} &
{ \Omega_{A/k} } \cr
}$$
is an isomorphism.
In particular, if $N$ is the free graded $k$-module generated by the graded set $S = \{x_i \mid i \in I\}$ so that $A$ is the free Dirac algebra $k[x_i \mid i \in I]$, then the graded $A$-module $\Omega_{A/k}$ is free with basis $(dx_i)_{i \in I}$, where $\deg(dx_i) = \deg(x_i)$. So given $f \in A$, we can write $df$ uniquely as a sum
$$\textstyle{ df = \sum_{i \in I} (\partial f/\partial x_i)dx_i, }$$
and we call $\partial f/\partial x_i$ the partial derivative of $f$ with respect to $x_i$. Beware that this involves signs! For instance, if $f = xy$ with both $x$ and $y$ of half-integer spin, then
$$df = dx \cdot y + x \cdot dy = -y \cdot dx + x \cdot dy,$$
so that $\partial f/\partial x = -y$, whereas $\partial f/\partial y = x$.
\end{example}

\begin{remark}
\label{remark:kahler_differentials_stable_under_base_change}
Let $\phi \colon k \to A$ and $\psi \colon k \to k'$ be two maps of Dirac rings, and let $\phi' \colon k' \to A'$ be the cobase-change of $\phi$ along $\psi$. It follows immediately from the definitions that the canonical map
$$\xymatrix{
{ k' \otimes_k\Omega_{A/k} \simeq A' \otimes_A\Omega_{A/k} } \ar[r] &
{ \Omega_{A'/k'} } \cr
}$$
is an isomorphism. In particular, if $\phi \colon k \to A$ is even, then
$$\xymatrix{
{ k \otimes_{k^{\operatorname{ev}}}\Omega_{A^{\operatorname{ev}}/k^{\operatorname{ev}}} } \ar[r] &
{ \Omega_{A/k} } \cr
}$$
is an isomorphism, and $d \colon A^{\operatorname{ev}} \to \Omega_{A^{\operatorname{ev}}/k^{\operatorname{ev}}}$ is given by the universal $k^{\operatorname{ev}}$-linear derivation of the underlying rings equipped with the induced (even) grading, where again $d$ is considered as a spin $0$ operator and not as a spin $\nicefrac{1}{2}$ operator.
\end{remark}

\begin{proposition}
\label{prop:differentialsexactsequences}
If $\phi \colon A \to B$ and $\psi \colon B \to C$ are composable maps of Dirac rings, then the induced sequence of graded $C$-modules
$$\xymatrix{
{ C \otimes_B\Omega_{B/A} } \ar[r] &
{ \Omega_{C/A} } \ar[r] &
{ \Omega_{C/B} } \ar[r] &
{ 0 } \cr
}$$
is exact. Moreover, if $\psi \colon B \to C$ is surjetive with kernel $I \subset B$, then this sequence extends to an exact sequence of graded $C$-modules
$$\xymatrix{
{ I/I^2 } \ar[r] &
{ C \otimes_B\Omega_{B/A} } \ar[r] &
{ \Omega_{C/A} } \ar[r] &
{ 0, } \cr
}$$
where the right-hand map takes $f + I^2$ to $1 \otimes df$. Finally, if $\psi \colon B \to C$ admits an $A$-algebra section, then the latter sequence is short exact.
\end{proposition}

\begin{proof}
Omitted; see e.g. \cite[\href{https://stacks.math.columbia.edu/tag/00RM}{Tag 00RM}]{stacks-project}.
\end{proof}

\begin{example}
\label{example:calculation_of_kahler_differentials_of_a_local_algebra}
If $k$ is a Dirac field and $(A,\mathfrak{m})$ a local Dirac $k$-algebra with residue field $k$, then the map $\mathfrak{m}/\mathfrak{m}^2 \to k \otimes_A\Omega_{A/k}$ that to $f + \mathfrak{m}^2$ assigns $1 \otimes df$ is an isomorphism.
\end{example}

\subsection{Smooth, unramified, and \'{e}tale maps}

\begin{definition}
\label{def:smoothunramifiedetale}
A map of Dirac rings $\phi \colon A \to B$ is smooth (resp.\ unramified, resp.\ \'{e}tale) if it is finitely presented and if for every diagram of Dirac rings
$$\xymatrix@C+=9mm{
{ A } \ar[r]^-{\psi_0} \ar[d]^-{\phi} &
{ C } \ar[d]^-{\pi} \cr
{ B } \ar[r]^-{\psi} &
{ C/I, } \cr
}$$
where $\pi$ is the projection onto the quotient by a square-zero graded ideal $I \subset C$, there exists a (resp.\ there exists at most one, resp.\ there exists a unique) map of Dirac rings $\varphi \colon B \to C$ such that $\pi \circ \varphi = \psi$ and $\varphi \circ \phi = \psi_0$.
\end{definition}

\begin{warning}[Smooth not necessarily flat]
\label{warning:smoothisnotflat}
A smooth map of Dirac rings need not be flat! Indeed, the unique map $\phi \colon \mathbb{Z} \to \mathbb{Z}[e]$, where $e$ is a generator of half-integer spin, is smooth, but it is not flat. We will show, however, that every \'{e}tale map of Dirac rings is flat.
\end{warning}

\begin{warning}[Smooth and finite not necessarily \'{e}tale]
A smooth and finite map of Dirac rings need not be \'{e}tale. For example, if $k$ is a Dirac ring in which $2$ is a unit, then the unique $k$-algebra map $\phi \colon k \to k[e]$, where $e$ is a generator of half-integer spin, is smooth and finite, but it is not \'{e}tale, since \cref{example:kahler_differentials_of_a_free_dirac_algebra} shows that $\Omega_{k[e]/k} \simeq k[e] \cdot de$ is nonzero.
\end{warning}

\begin{lemma}
\label{lemma:etale_maps_stable_under_base_change}The composition of two smooth (resp.\ unramified, resp.\ \'{e}tale) maps of Dirac rings is a smooth (resp.\ unramified, resp.\ \'{e}tale) map, and the cobase-change of a smooth (resp.\ unramified, resp.\ \'{e}tale) map of Dirac rings along any map of Dirac rings is a smooth (resp.\ unramified, resp.\ \'{e}tale) map.
\end{lemma}

\begin{proof}This follows immediately from the definition.
\end{proof}

\begin{remark}We will define the cotangent complex $L_{B/A}$ in a sequel to this paper and show that $\phi \colon A \to B$ is smooth (resp.\ unramified, resp.\ \'{e}tale) if and only if it is finitely presented and $L_{B/A}$ is perfect with $\operatorname{Tor}$-amplitude in $[0,0]$ (resp.\ is connected, resp. is contractible).
\end{remark}

\begin{addendum}\label{add:differentialssmooth}If $\phi \colon A \to B$ and $\psi \colon B \to C$ are composable maps of Dirac rings with $\psi$ smooth, then the induced sequence of graded $C$-modules
$$\xymatrix{
{ 0 } \ar[r] &
{ C \otimes_B\Omega_{B/A} } \ar[r] &
{ \Omega_{C/A} } \ar[r] &
{ \Omega_{C/B} } \ar[r] &
{ 0 } \cr
}$$
is short exact and splittable.
\end{addendum}

\begin{proof}Omitted; see~\cite[Theorem~25.1]{matsumura}.
\end{proof}

\begin{lemma}\label{lem:smoothunramifiedetaleeven}
A map of Dirac rings $\phi \colon A \to B$ between even Dirac rings is smooth (resp.\ unramified, resp.\ \'{e}tale) if and only if the map of underlying rings is smooth (resp.\ unramified, resp. \'{e}tale).
\end{lemma}

\begin{proof}If $P$ is an evenly generated free $A$-algebra, then $\widetilde{P}$ is a free $\widetilde{A}$-algebra, so the statement for smooth maps follows from \cite[\href{https://stacks.math.columbia.edu/tag/00TL}{Tag 00TL}]{stacks-project}. The statement for unramified maps follows from \cite[\href{https://stacks.math.columbia.edu/tag/00UO}{Tag 00UO}]{stacks-project} and \ref{remark:kahler_differentials_stable_under_base_change}.
\end{proof}

\begin{definition}\label{definition:standard_smooth}
A map of Dirac rings $\phi \colon A \to B$ is standard smooth, if there exists $c \geq 0$ homogeneous polynomials $f_1,\dots,f_c \in A[x_1,\dots,x_n]$ in $n \geq c$ homogeneous variables $x_1,\dots,x_n$ such that $B \simeq A[x_1,\dots,x_n]/(f_1,\dots,f_c)$ and such that the Jacobian matrix $(\partial f_i / \partial x_j)$ has maximal rank $c$ over $B$.
\end{definition}

We recall that the partial derivatives $\partial f_i/\partial x_j$ were defined in \cref{example:kahler_differentials_of_a_free_dirac_algebra} and that the definition involves signs. That the Jacobian matrix has maximal rank $c$ means that the associated $B$-linear map $B^n \to B^c$ is surjective.

\begin{lemma}\label{lemma:standard_smooth_is_smooth}
A standard smooth map of Dirac rings $\phi \colon A \to B$ is smooth.
\end{lemma}

\begin{proof}
Arguing as in \cite[\href{https://stacks.math.columbia.edu/tag/031I}{Tag 031I}]{stacks-project}, it suffices to show that the sequence
$$\xymatrix{
{ 0 } \ar[r] &
{ I/I^2 } \ar[r] &
{ B \otimes_A \Omega_{A[x_1,\dots,x_n]/A} } \ar[r] &
{ \Omega_{B/A} } \ar[r] &
{ 0 } \cr
}$$
with $I = (f_{1}, \ldots, f_{c})$ is exact, and by  \cref{prop:differentialsexactsequences}, only it is only the injectivity of the left-hand map that is at issue. The domain and target of this map are both free graded $B$-modules, and the matrix that represents it with respect to the bases $(f_1+I^2,\dots,f_c+I^2)$ of the domain and $(dx_1,\dots,dx_n)$ of the target is the transpose of the Jacobian matrix. Thus, the map is injective if and only if the Jacobian matrix has maximal rank.
\end{proof}

\begin{remark}\label{remark:standard_smooth}Conversely, one can prove that, Zariski locally on the source, every smooth map of Dirac rings $\phi \colon A \to B$ is standard smooth. We will not need this result, so we will not give the (rather lengthy) proof here.
\end{remark}

Let $K$ be a Dirac field, and let $K_0 \subset K$ be the subfield consisting of the homogeneous elements of degree $0$. We recall from \cref{prop:classification_of_dirac_fields} that either 
\begin{enumerate}
\item $K = K_{0}$, or
\item $K = K_{0}[t^{\pm 1}]$ for every homogeneous $t \in K^{\times}$ with $|\operatorname{spin}(t)| \neq 0$ minimal.
\end{enumerate}
Moreover, if $\operatorname{char}(K) \neq 2$, then, in the second case, the generator $t$ necessarily is of integer spin. In all cases, the underlying ring is commutative.

\begin{lemma}
\label{lemma:classification_of_finite_ext_of_dirac_fields}
If $K \to L$ is a finite extension of Dirac fields, then $K_0 \to L_0$ is a finite extension of fields, and the following hold:
\begin{enumerate}
\item[{\rm (1)}] If $K = K_{0}$, then $L = L_{0}$.
\item[{\rm (2)}] If $K = K_0[t^{\pm1}]$, then $L = L_0[u^{\pm1}]$ with $\lambda_{0} u^{e} = t$ for some $e \geq 1$ and $\lambda_{0} \in L_{0}$.
\end{enumerate}
Moreover, in the latter case, $e = [L:L']$, where $L' \simeq L_0 \otimes_{K_0}K \to L$.
\end{lemma}

\begin{proof}
If $K = K_0$, but $L = L_0[u^{\pm}]$, then the extension $K \to L$ is infinite, so~(1) is clear. If $K = K_0[t^{\pm1}]$ with $\operatorname{spin}(t) > 0$, then we write $L = L_0[u^{\pm1}]$ with $\operatorname{spin}(u) > 0$. In particular, we can write $t \in K \subset L$ uniquely as $t = \lambda_0u^e$ with $\lambda_0 \in L_0$ and $e \geq 1$. If $\operatorname{spin}(t) < 0$, then we argue analogously.
\end{proof}

\begin{warning}\label{warning:diracfieldextensions}In \cref{lemma:classification_of_finite_ext_of_dirac_fields}~(2), it is not always possible to choose $t \in K$ and $u \in L$ such that $t = u^e$. For instance, let $K = \mathbb{Q}[t^{\pm1}] \subset L = \mathbb{Q}(i)[u^{\pm1}]$ with $\operatorname{spin}(u) = 1$ and $u^2 = it$. In this case, there is no $v \in L$ of spin $1$ with $v^2 = qt \in K$. For if there were, then $iq^{-1} = (uv^{-1})^2 \in L_0 = \mathbb{Q}(i)$ would be a square root.
\end{warning}

\begin{proposition}\label{proposition:unramified_extension_of_Dirac_fields}
An unramified extension of Dirac fields $K \to L$ is even, finite, and \'{e}tale.
\end{proposition}

\begin{proof}If $L = L_0$, then $K = K_0$, and there is nothing to prove. So we assume that $L = L_0[u^{\pm1}]$. In this case, we have $K = K_0[t^{\pm1}]$. For \cref{example:kahler_differentials_of_a_free_dirac_algebra} shows that $\Omega_{L/L_0} \simeq L \cdot du$ is nonzero, so if $K = K_0$, then \cref{prop:differentialsexactsequences} implies that $\Omega_{L/K} \simeq \Omega_{L/K_0}$ is nonzero, contradicting the assumption that $K \to L$ is unramified; see~\cite[\href{https://stacks.math.columbia.edu/tag/00UO}{Tag 00UO}]{stacks-project}. We now factor the extension in question as the composition
$$\xymatrix{
{ K } \ar[r] &
{ L' \simeq L_0 \otimes_{K_0}K } \ar[r] &
{ L } \cr
}$$
and note that $L = L'[u]/(u^e+\lambda t)$ with $e = [L:L']$ and $\lambda \in L_0^{\times}$. If $K \to L$ is not even, then neither is $L' \to L$. But this can happen only if the degree of $u$ is odd and the degree of $t$ is even, in which case, the characteristic of $L$ is $2$ and $e = [L:L']$ is even. This, in turn, implies that
$$\Omega_{L/L'} \simeq L/(eu^{e-1}) \cdot du \simeq L \cdot du$$
is nonzero, which, by \cref{prop:differentialsexactsequences}, contradicts that $\Omega_{L/K}$ is zero.

Finally, since $K \to L$ is even, the fact that $K \to L$ is unramified implies that also $K^{\operatorname{ev}} \to L^{\operatorname{ev}}$ is unramified. So we conclude from \cref{lem:smoothunramifiedetaleeven} and \cite[\href{https://stacks.math.columbia.edu/tag/090W}{Tag 090W}]{stacks-project} that $K^{\operatorname{ev}} \to L^{\operatorname{ev}}$ is finite and \'{e}tale, and therefore, so is $K \to L$.
\end{proof}

We proceed to prove that any \'{e}tale algebra over a Dirac field decomposes as a finite product of \'{e}tale Dirac field extensions. Our proof is an adaptation of the usual argument, with care taken whenever one has to manipulate elements. 

\begin{lemma}
\label{lemma:maximal_ideal_of_local_artin_algebra_is_nilpotent}Let $A$ be an artinian local Dirac ring. Its maximal ideal $\mathfrak{m} \subset A$ is nilpotent.
\end{lemma}

\begin{proof}Since $A$ is artinian, there exists $n \geq 1$ such that $\mathfrak{m}^i = \mathfrak{m}^{i+1}$ for all $i \geq n$. If $\mathfrak{m}^n$ is nonzero, then we let $I \subset A$ be a (possibly non-proper) minimal ideal such that $I\mathfrak{m}^n$ is nonzero. But $(I\mathfrak{m}^n)\mathfrak{m}^n = I\mathfrak{m}^{2n} = I\mathfrak{m}^n$ is again nonzero, so by the minimality of $I \subset A$, we conclude that $I\mathfrak{m}^n = I$, and hence, that $I\mathfrak{m} = I$, so by Nakayama's lemma, \cref{lem:nakayama}, we conclude that $I$ is zero, which is a contradiction.
\end{proof}

\begin{proposition}\label{proposition:etale_algebra_over_dirac_field_is_finite_and_even}
If $k$ is a Dirac field, and let $\phi \colon k \to A$ be an \'{e}tale map. The $k$-algebra $\phi \colon k \to A$ decomposes as a product of a finite family $(k \to k_i)_{i \in I}$ of \'{e}tale Dirac field extensions. In particular, it is finite and even.
\end{proposition}

\begin{proof}We first show that $\phi \colon k \to A$ is finite. To this end, we may assume that $k$ is algebraically closed, since an algebraic closure $k \to \bar{k}$ is faithfully flat, and hence, reflects the property of being finite. So if $\mathfrak{m} \subset A$ is a maximal graded ideal, then $k \to k(\mathfrak{m})$ is an isomorphism, and hence, the canonical map
$$\xymatrix{
{ \mathfrak{m}/\mathfrak{m}^2 } \ar[r] &
{ k(\mathfrak{m}) \otimes_{A_{\mathfrak{m}}} \Omega_{A_{\mathfrak{m}}/k} } \cr
}$$
is an isomorphism by \cref{example:calculation_of_kahler_differentials_of_a_local_algebra}. The right-hand side is zero, by the assumption that $\phi \colon k \to A$ is \'{e}tale, and hence, unramified, so we conclude from Nakayama's lemma, \cref{lem:nakayama}, that $A_{\mathfrak{m}} \to k(\mathfrak{m})$ is an isomorphism. Therefore, every maximal graded ideal in $A$ is a minimal graded prime ideal. Since $A$ is noetherian, this implies that $A$ is artinian by the classical argument \cite[\href{https://stacks.math.columbia.edu/tag/00KH}{Tag 00KH}]{stacks-project}. This proves that $\phi \colon k \to A$ is finite.

We now let $k$ be a general Dirac field and choose an algebraic closure $k \to \bar{k}$. The classical argument \cite[\href{https://stacks.math.columbia.edu/tag/00JA}{Tag 00JA}]{stacks-project} shows that $\phi \colon k \to A$ decomposes as a finite product of artinian local $k$-algebras. So we may assume that $\phi \colon k \to A$ is local. By \cref{lemma:maximal_ideal_of_local_artin_algebra_is_nilpotent}, the maximal ideal $\mathfrak{m} \subset A$ is nilpotent. Now, as a map of graded $k$-vector spaces, the extension $k \to \bar{k}$ admits a retraction, so the induced map
$$\xymatrix{
{ A \simeq A \otimes_kk } \ar[r] &
{ A\otimes_k\bar{k} } \cr
}$$
is (split) injective. The target is an \'{e}tale $\bar{k}$-algebra, and hence, it is a finite product of copies of $\bar{k}$ by the argument above. So we conclude that $A$ has no nonzero nilpotent homogeneous elements, and therefore, that $\mathfrak{m} \subset A$ is the zero ideal. It follows that $A \to k(\mathfrak{m})$ is an isomorphism, so we conclude that $\phi \colon k \to A \simeq k(\mathfrak{m})$ is an \'{e}tale extension of fields, as desired. Finally, it follows from \cref{proposition:unramified_extension_of_Dirac_fields} that $\phi \colon k \to A$ is even.
\end{proof}

\begin{corollary}
\label{corollary:etale_extensions_are_quasifinite}
An \'{e}tale map of Dirac rings $\phi \colon A \to B$ is quasi-finite. 
\end{corollary}

\begin{proof}
This is an consequence of \cref{proposition:etale_algebra_over_dirac_field_is_finite_and_even} and the fact that \'{e}tale maps of Dirac rings are stable under cobase-change.
\end{proof}

\begin{corollary}
\label{corollary:etale_morphisms_unramified_in_terms_of_ideals}
If $\phi \colon A \to B$ is an \'{e}tale map of Dirac rings, and if $\mathfrak{q} \subset B$ is a graded prime ideal lying over a graded prime ideal $\mathfrak{p} \subset A$, then $\mathfrak{p}B_{\mathfrak{q}} = \mathfrak{q}B_{\mathfrak{q}}$.
\end{corollary}

\begin{proof}
By \cref{proposition:etale_algebra_over_dirac_field_is_finite_and_even}, the Dirac ring $k(\mathfrak{p}) \otimes_{A} B$ is a finite product of Dirac fields, and hence, the localization $k(\mathfrak{p}) \otimes_{A} B_{\mathfrak{q}}$ is a Dirac field. It follows that the kernel $\mathfrak{p}B_{\mathfrak{q}}$ of $B_{\mathfrak{q}} \to k(\mathfrak{p}) \otimes_AB_{\mathfrak{q}}$ is the maximal graded ideal, as stated.
\end{proof}

The following standard property of unramified maps will be useful later. 

\begin{proposition}
\label{proposition:unramified_finite_type_morphisms_have_flat_multiplication}
Let $A \to B$ be an unramified morphism of Dirac rings of finite type. Then the multiplication map $B \otimes_{A} B \to B$ is flat. 
\end{proposition}

\begin{proof}
The standard argument works here, but let us recall it for the convenience of the reader. Let us write $I$ for the kernel of multiplication. Since $B$ is unramified over $A$, we have 
$$I/I^{2} \simeq \Omega^{1}_{B/A} = 0$$
so that $I = I^{2}$. 

We claim that $I$ is finitely generated. To see this, let $x_{1}, \ldots, x_{n}$ be generators of $B$ as an $A$-algebra, and let $I'$ denote the ideal of $B \otimes_{A} B$ generated by $x_{i} \otimes 1 - 1 \otimes x_{i}$ for $1 \leq i \leq n$. Now suppose we have two maps $f, g: B \to C$ into another Dirac $A$-algebra $C$ such that the unique induced homomorphism $B \otimes_{A} B \to C$ annihilates $I'$. It follows that $f(x_{i}) \cdot 1 - 1 \cdot g(x_{i}) = 0$, so that $f(x_{i}) = g(x_{i})$ and hence $f$ and $g$ coincide. Thus, any such homomorphism factors through $B \simeq (B \otimes_{A} B)/I$ and $I = I'$, showing that it's finitely generated. 

To check that $B$ is flat over $B \otimes_{A} B$ it is enough to verify that for any maximal ideal $\mathfrak{m}$ of the latter, the induced map 
$$(B \otimes_{A} B)_{\mathfrak{m}} \to B_{\mathfrak{m}}$$
is flat. The kernel of this map is given by the finitely generated ideal $I_{\mathfrak{m}}$ which satisfies $I_{\mathfrak{m}} = I^{2}_{\mathfrak{m}} = (I_{\mathfrak{m}})^{2}$. If $I_{\mathfrak{m}} = (B \otimes_{A} B)_{\mathfrak{m}}$, then $B_{\mathfrak{m}} = 0$ which is flat. Otherwise, we have $I_{\mathfrak{m}} \subseteq \mathfrak{m}$, so that since 
$$I_{\mathfrak{m}} \otimes (B \otimes_{A} B)_{\mathfrak{m}}/I_{\mathfrak{m}} \simeq I_{\mathfrak{m}} / I_{\mathfrak{m}}^{2} = 0,$$
also $I_{\mathfrak{m}} \otimes (B \otimes_{A} B)_{\mathfrak{m}}/\mathfrak{m} = 0$. Nakayama lemma of \cref{lem:nakayama} then implies that $I_{\mathfrak{m}} = 0$, hence $B_{\mathfrak{m}} \simeq (B \otimes_{A} B)_{\mathfrak{m}}$ which is also flat. 
\end{proof}

\subsection{Evenness of \'{e}tale extensions}

We have seen in \cref{proposition:etale_algebra_over_dirac_field_is_finite_and_even} that every \'{e}tale extension of a Dirac field is even. We proceed to generalize this to arbitraty maps of Dirac rings. We define
$$\operatorname{CAlg}(\operatorname{Mod}_A(\operatorname{Ab}))^{\operatorname{\acute{e}t}} \subset \operatorname{CAlg}(\operatorname{Mod}_A(\operatorname{Ab}))$$
to be the full subcategory spanned by the \'{e}tale algebras.

\begin{theorem}
\label{theorem:equivalence_of_categories_between_etale_extensions_of_a_and_even_subring}
Let $A$ be a Dirac ring. The functor
$$\xymatrix{
{ \operatorname{CAlg}(\operatorname{Mod}_{A^{\operatorname{ev}}}(\operatorname{Ab}))^{\text{\rm\'{e}t}} } \ar[r] &
{ \operatorname{CAlg}(\operatorname{Mod}_A(\operatorname{Ab}))^{\text{\rm\'{e}t}} } \cr
}$$
given by extension of scalars along $A^{\operatorname{ev}} \to A$ is an equivalence of categories.
\end{theorem}

\begin{corollary}
An \'{e}tale map of Dirac rings $\phi \colon A \to B$ is flat. 
\end{corollary}

\begin{proof}
By \cref{theorem:equivalence_of_categories_between_etale_extensions_of_a_and_even_subring}, it suffices to show that $\phi^{\operatorname{ev}} \colon A^{\operatorname{ev}} \to B^{\operatorname{ev}}$ is flat. Now, by \cref{lem:smoothunramifiedetaleeven}, the map $\phi^{\operatorname{ev}} \colon A^{\operatorname{ev}} \to B^{\operatorname{ev}}$ is \'{e}tale if and only if the map of underlying rings is \'{e}tale. But an \'{e}tale map of rings is flat by \cite[\href{https://stacks.math.columbia.edu/tag/00TA}{Tag 00TA}]{stacks-project}. 
\end{proof}

We will prove \cref{theorem:equivalence_of_categories_between_etale_extensions_of_a_and_even_subring} below, after a few preliminary results. 

\begin{lemma}
\label{lemma:an_extension_of_even_map_is_etale_iff_the_map_was}
An even map of Dirac rings $\phi \colon A \to B$ is \'{e}tale if and only if the induced map of even Dirac rings $\phi^{\operatorname{ev}} \colon A^{\operatorname{ev}} \to B^{\operatorname{ev}}$ is \'{e}tale.
\end{lemma}

\begin{proof}If $\phi^{\operatorname{ev}} \colon A^{\operatorname{ev}} \to B^{\operatorname{ev}}$ is \'{e}tale, then so is $\phi \colon A \to B \simeq A \otimes_{A^{\operatorname{ev}}}B^{\operatorname{ev}}$, since \'{e}tale maps are stable under cobase-change. Conversely, if $\phi \colon A \to B$ is \'{e}tale, then it follows from \cref{proposition:properties_of_even_extensions_reflected_by_base_change} that $\phi^{\operatorname{ev}} \colon A^{\operatorname{ev}} \to B^{\operatorname{ev}}$ is finitely presented, so we only have to check that it has the unique lifting property with respect to square-zero extensions. Now, if we are given the outer square
$$\xymatrix@C+=9mm{
{ A^{\operatorname{ev}} } \ar[r] \ar[d] &
{ C^{\operatorname{ev}} } \ar[r] \ar[d] &
{ C } \ar[d] \cr
{ B^{\operatorname{ev}} } \ar[r] &
{ (C/I)^{\operatorname{ev}} } \ar[r] &
{ C/I } \cr
}$$
with $I \subset C$ a square-zero graded ideal, then it factors uniquely as indicated, and the middle vertical map again is a square-zero extension. In the square
$$\xymatrix@C+=9mm{
{ A } \ar[r] \ar[d] &
{ A \otimes_{A^{\operatorname{ev}}}C^{\operatorname{ev}} } \ar[d] \cr
{ B } \ar[r] &
{ A \otimes_{A^{\operatorname{ev}}}(C/I)^{\operatorname{ev}} } \cr
}$$
obtained by extension of scalars along $A^{\operatorname{ev}} \to A$, the right-hand vertical map is a square-zero extension of Dirac rings, and the left-hand vertical map is \'{e}tale, so there exists a unique map of Dirac rings $B \to A \otimes_{A^{\operatorname{ev}}}C^{\operatorname{ev}}$ making the diagram commute. This map restricts to a lifting $B^{\operatorname{ev}} \to C^{\operatorname{ev}}$ in the top left-hand square, which, in turn, gives a lifting in the given outer square. Finally, this lifting is unique, since $A^{\operatorname{ev}} \to A$ is faithful, as $A \simeq A^{\operatorname{ev}} \oplus A^{\operatorname{odd}}$ as a graded $A^{\operatorname{ev}}$-module.
\end{proof}

The following result is key to establishing \cref{theorem:equivalence_of_categories_between_etale_extensions_of_a_and_even_subring}.

\begin{proposition}
\label{proposition:etale_algebra_is_evenly_generated}
An \'{e}tale map of Dirac rings $\phi \colon A \to B$ is evenly generated. 
\end{proposition}

\begin{proof}
We wish to show that the canonical map
$$\xymatrix{
{ A \otimes_{A^{\operatorname{ev}}} B^{\operatorname{ev}} } \ar[r] &
{ B } \cr
}$$
is surjective. This is a map of $B^{\operatorname{ev}}$-modules, so it suffices to show that for every graded prime ideal $\mathfrak{q}^{\operatorname{ev}} \subset B^{\operatorname{ev}}$, the induced map of localizations
$$\xymatrix{
{ (A \otimes_{A^{\operatorname{ev}}} B^{\operatorname{ev}})_{\mathfrak{q}^{\operatorname{ev}}} } \ar[r] &
{ B_{\mathfrak{q}^{\operatorname{ev}}} } \cr
}$$
is surjective. We recall from \cref{prop:evendiraczariskispace} that, as our notation indicates, every graded prime ideal $\mathfrak{q}^{\operatorname{ev}} \subset B^{\operatorname{ev}}$ is the even part of a unique prime ideal $\mathfrak{q} \subset B$, and, moreover, we may identify the map in question with the canonical map
$$\xymatrix{
{ A \otimes_{A^{\operatorname{ev}}}(B_{\mathfrak{q}})^{\operatorname{ev}} } \ar[r] &
{ B_{\mathfrak{q}}. } \cr
}$$
Therefore, it will suffice to show that for every graded prime ideal $\mathfrak{q} \subset B$, the local ring $B_{\mathfrak{q}}$ is evenly generated as an $A$-algebra. 

Now, it follows from \cref{corollary:etale_extensions_are_quasifinite} that the map $\phi \colon A \to B$ is quasi-finite, so by Zariski's main theorem, \cref{corollary:finite_open_immersion_factorization_of_qf_morphism}, it admits a factorization
$$\xymatrix{
{ A } \ar[r]^-{\phi_0} &
{ B_0 } \ar[r]^-{\eta} &
{ B } \cr
}$$
with $\phi_0$ is finite and with the map of Dirac spectra induced by $\eta$ an open immersion. Let $\mathfrak{p} = \phi^{-1}(\mathfrak{q}) \subset A$ and $\mathfrak{q}_0 = \eta^{-1}(\mathfrak{q}) \subset B_0$. Since $\eta$ induces an open immersion of Dirac spectra, it induces an isomorphism of local rings
$$\xymatrix{
{ (B_0)_{\mathfrak{q}_0} } \ar[r] &
{ B_{\mathfrak{q}}, } \cr
}$$
and hence, it suffices to show that $A \to (B_0)_{\mathfrak{q}_0}$ is evenly generated. The assumption that  $\phi \colon A \to B$ is \'{e}tale implies that also $\phi(\mathfrak{p}) \colon k(\mathfrak{p}) \to k(\mathfrak{p}) \otimes_AB$ is \'{e}tale, and therefore, we conclude from \cref{proposition:etale_algebra_over_dirac_field_is_finite_and_even} that $\phi(\mathfrak{p})$ is even. Moreover, it follows from  \cref{corollary:etale_morphisms_unramified_in_terms_of_ideals} that $k(\mathfrak{p}) \otimes_AB_{\mathfrak{q}} \simeq k(\mathfrak{q})$, which, combined with the isomorphism above, shows that
$$k(\mathfrak{p}) \otimes_A(B_0)_{\mathfrak{q}_0}  \simeq 
k(\mathfrak{q}_0).$$
This proves the following facts, which we will use below:
\begin{enumerate}
    \item[(a)]The graded ideal $\mathfrak{p}(B_0)_{\mathfrak{q}_0} \subset (B_0)_{\mathfrak{q}_0}$ is the maximal graded ideal.
    \item[(b)]The $A$-algebra $k(\mathfrak{q}_0)$ is evenly generated.
\end{enumerate}
We further factor $\phi_0$ as
$$\xymatrix{
{ A } \ar[r]^-{\phi_1} &
{ B_1 \simeq A \otimes_{A^{\operatorname{ev}}}B_0^{\operatorname{ev}} } \ar[r]^-{\psi} &
{ B_0 } \cr
}$$
and let $\mathfrak{q}_1 = \psi^{-1}(\mathfrak{q}_0) \subset B_1$. Since $\psi^{\operatorname{ev}} \colon B_1^{\operatorname{ev}} \to B_0^{\operatorname{ev}}$ is an isomorphism, so is
$$\xymatrix{
{ (B_0)_{\mathfrak{q}_1} } \ar[r] &
{ (B_0)_{\mathfrak{q}_0}. } \cr
}$$
We now consider the map
$$\xymatrix@C+=11mm{
{ (B_1)_{\mathfrak{q}_1} } \ar[r]^-{\psi_{\mathfrak{q}_1}} &
{ (B_0)_{\mathfrak{q}_1} \simeq (B_0)_{\mathfrak{q}_0}. } \cr
}$$
Since $\phi_0 \colon A \to B_0$ is finite, so is $\psi \colon B_1 \to B_0$, and hence, the map $\psi_{\mathfrak{q}_1}$ is a map of finite $(B_1)_{\mathfrak{q}_1}$-modules. Moreover, it follows from~(a) above $\mathfrak{q}_1(B_0)_{\mathfrak{q}_0} \subset (B_0)_{\mathfrak{q}_0}$ is the maximal graded ideal. Therefore, we conclude that the diagram
$$\xymatrix@C+=12mm{
{ (B_1)_{\mathfrak{q}_1} } \ar[r]^-{\psi_{\mathfrak{q}_1}} \ar[d] &
{ (B_0)_{\mathfrak{q}_1} } \ar[d] \cr
{ k(\mathfrak{q}_1) } \ar[r]^-{\psi(\mathfrak{q}_1)} &
{ k(\mathfrak{q}_0) } \cr
}$$
is cocartesian. The map $\psi(\mathfrak{q}_1)$ induces an isomorphism of the even parts, so it follows from~(b) above that it is surjective. (It is also injective, but we this is not relevant for us.) By Nakayama's lemma, \cref{lem:nakayama}, so is the map $\psi_{\mathfrak{q}_1}$. This shows that $(B_0)_{\mathfrak{q}_0}$ admits a surjective map from an evenly generated $A$-algebra, and therefore, we conclude that it is itself an evenly generated $A$-algebra, as desired.
\end{proof}

\begin{lemma} 
\label{lemma:evenly_generated_etale_algebra_is_even}
An evenly generated \'{e}tale map of Dirac rings $\phi \colon A \to B$ is even.
\end{lemma}

\begin{proof}
The canonical map $A \otimes_{A^{\operatorname{ev}}} B^{\operatorname{ev}} \to B$ is surjective, by assumption, and we must show that it is an isomorphism. Its kernel $I \subset A \otimes_{A^{\operatorname{ev}}} B^{\operatorname{ev}}$ is concentrated in odd degree, and hence, is a square-zero ideal. Since $\phi \colon A \to B$ is \'{e}tale, we conclude that the canonical map of Dirac rings $A \otimes_{A^{\operatorname{ev}}}B^{\operatorname{ev}} \to B$ admits a section. This section, in particular, is a map of graded $A$-modules, and hence, induces a map
$$\xymatrix{
{ A \otimes_{A^{\operatorname{ev}}} B^{\operatorname{ev}} } \ar[r]^-{r} &
{ I } \cr
}$$
 of graded $A$-modules, which is a retraction of the canonical inclusion. The map $r$ determines and is determined by the map $r^{\operatorname{ev}} \colon B^{\operatorname{ev}} \to I^{\operatorname{ev}}$ of graded $A^{\operatorname{ev}}$-modules, which is necessarily zero, since $I^{\operatorname{ev}} \simeq 0$. So the map $r$ is zero, but it is also surjective, and therefore, we conclude that $I$ is zero, as desired.
\end{proof}

We are now ready to prove the main result of this section. 

\begin{proof}[Proof of \cref{theorem:equivalence_of_categories_between_etale_extensions_of_a_and_even_subring}]
We let $A$ be a Dirac ring, and define
$$\operatorname{CAlg}(\operatorname{Mod}_A(\operatorname{Ab}))^{\operatorname{ev}} \subset \operatorname{CAlg}(\operatorname{Mod}_A(\operatorname{Ab}))$$
to be the full subcategory spanned by the even $A$-algebras. By the definition of even algebras, we have an adjoint equivalence
$$\xymatrix{
{ \operatorname{CAlg}(\operatorname{Mod}_{A^{\operatorname{ev}}}(\operatorname{Ab}))^{\operatorname{ev}} } \ar@<.7ex>[r] &
{ \operatorname{CAlg}(\operatorname{Mod}_A(\operatorname{Ab}))^{\operatorname{ev}} } \ar@<.7ex>[l] \cr
}$$
with the left adjoint given by extension of scalars along $A^{\operatorname{ev}} \to A$ and with the right adjoint given by the functor $(-)^{\operatorname{ev}}$ that takes even parts. By \cref{lemma:an_extension_of_even_map_is_etale_iff_the_map_was}, an even map $\phi \colon A \to B$ is \'{e}tale if and only if the map $\phi^{\operatorname{ev}} \colon A^{\operatorname{ev}} \to B^{\operatorname{ev}}$ is \'{e}tale. Therefore, the adjoint equivalence above restricts to an adjoint equivalence
$$\xymatrix{
{ \operatorname{CAlg}(\operatorname{Mod}_{A^{\operatorname{ev}}}(\operatorname{Ab}))^{\operatorname{\acute{e}t},\operatorname{ev}} } \ar@<.7ex>[r] &
{ \operatorname{CAlg}(\operatorname{Mod}_A(\operatorname{Ab}))^{\operatorname{\acute{e}t},\operatorname{ev}} } \ar@<.7ex>[l] \cr
}$$
between the full subcategories of the categories in the statement of the theorem spanned by the even \'{e}tale algebras. Thus, it is suffices to show that every \'{e}tale algebra is even, and this follows from \cref{proposition:etale_algebra_is_evenly_generated} and \cref{lemma:evenly_generated_etale_algebra_is_even} above.
\end{proof}

\subsection{\'{E}tale maps of ring spectra}
\label{subsection:etale_maps_of_ring_spectra}

Lurie, in~\cite[Definition~7.5.0.4]{lurieha}, defines a map of algebras in spectra $R \to A$ to be \'{e}tale if 
\begin{enumerate}
    \item[(1)] the induced map $\pi_0(R) \to \pi_0(A)$ is an \'{e}tale map of commutative rings, and
    \item[(2)] the induced map $\pi_0(A) \otimes_{\pi_0(R)}\pi_n(R) \to \pi_n(A)$ is an isomorphism for all $n \in \mathbb{Z}$.
\end{enumerate}
We will refer to Lurie's notion as ``$\pi_0$-\'{e}tale'' and will reserve the term ``\'{e}tale'' for the more general notion that we introduce in \cref{definition:etale_map_of_ek_rings} below.

At the level of homotopy groups, a $\pi_0$-\'{e}tale map $R \to A$ is determined completely by its induced map on $\pi_0$, but it turns out that a much stronger statement is true. This statement is Lurie's rigidity theorem for $\pi_0$-\'{e}tale maps, which we now recall. If $R$ is an $\mathbf{E}_{k+1}$-algebra in spectra with $1 \leq k \leq \infty$, then $\operatorname{Mod}_R(\operatorname{Sp})$ promotes to an $\mathbf{E}_k$-monoidal $\infty$-category, and therefore, it is meaningful to consider the $\infty$-category $\operatorname{Alg}_{\mathbf{E}_k}(\operatorname{Mod}_R(\operatorname{Sp}))$ of $\mathbf{E}_k$-algebras. Now, Lurie's rigidity theorem for $\pi_0$-\'{e}tale maps, which is \cite[Theorem~7.5.0.6]{lurieha}, states that the functor
$$\xymatrix{
{ \operatorname{Alg}_{\mathbf{E}_k}(\operatorname{Mod}_R(\operatorname{Sp}))^{\text{$\pi_0$-\'{e}t}} } \ar[r]^-{\pi_0} &
{ \operatorname{CAlg}(\operatorname{Mod}_{\pi_0(R)}(\operatorname{Ab}_0))^{\operatorname{\acute{e}t}} } \cr
}$$
is an equivalence from the full subcategory of $\operatorname{Alg}_{\mathbf{E}_k}(\operatorname{Mod}_R(\operatorname{Sp}))$ spanned by the $\pi_0$-\'{e}tale algebras to the full subcategory of $\operatorname{CAlg}(\operatorname{Mod}_{\pi_0(R)}(\operatorname{Ab}_0))$ spanned by the \'{e}tale algebras. (Here $\operatorname{Ab}_0$ is the symmetric monoidal category of abelian groups.) 

A goal of Dirac geometry is to provide a language to describe phenomena relating to $\pi_*(E)$ for algebras in spectra $E$, which may not be captured by $\pi_0(E)$. In this spirit, it is natural to make the following definition. 

\begin{definition}
\label{definition:etale_map_of_ek_rings}
If $R$ is an $\mathbf{E}_{k+1}$-algebra in spectra with $1 \leq k \leq \infty$, then an $\mathbf{E}_k$-algebra $A$ in $\operatorname{Mod}_R(\operatorname{Sp})$ is \'{e}tale if the map $\pi_*(R) \to \pi_*(A)$ induced by the unit is an \'{e}tale map of Dirac algebras.
\end{definition}

\begin{remark}
\label{remark:etale_map_of_e1_rings}
If $A$ is an $\mathbf{E}_1$-algebra in $\operatorname{Mod}_R(\operatorname{Sp})$, then the graded ring $\pi_*(A)$ may not be a Dirac ring. So for $k = 1$, the requirement that $\pi_*(A)$ be a Dirac ring is part of the definition of $A$ being \'{e}tale.
\end{remark}

\begin{example}
\label{example:homotopy_fixed_points_in_p_complete_k_theory_etale}
Let $p$ be a prime number, and let $KU_p$ be the commutative algebra in spectra representing $p$-complete complex $K$-theory. It carries an action through Adams operations of the group $\operatorname{Aut}(\mu_{p^{\infty}}) = \mathbb{Z}_p^{\times}$. Hence, the forgetful map
$$\xymatrix{
{ KU_{p}^{\,h\operatorname{Aut}(\mu_p)} } \ar[r]^-{\theta} &
{ KU_{p} } \cr
}$$
from the homotopy fixed points of the subgroup $\operatorname{Aut}(\mu_p) \subset \operatorname{Aut}(\mu_{p^{\infty}})$, which is cyclic of order $p-1$, to the underlying spectrum is a map of commutative algebras in spectra. The group $\pi_2(KU_p)$ is a free $\mathbb{Z}_p$-module of rank $1$, on which the induced action by $\operatorname{Aut} (\mu_{p^{\infty}}) = \mathbb{Z}_p^{\times}$ is given by scalar multiplication. Hence, if we choose a generator $\beta \in \pi_2(KU_p)$, then we may identify $\pi_*(\theta)$ with the canonical inclusion
$$\xymatrix{
{ \mathbb{Z}_p[\beta^{\pm(p-1)}] } \ar[r] &
{ \mathbb{Z}_p[\beta^{\pm1}], } \cr
}$$
which is an \'{e}tale map of Dirac rings. Accordingly, the map $\theta$ is \'{e}tale in our sense, but it is not $\pi_0$-\'{e}tale.
\end{example}

The following generalizes Lurie's rigidity theorem~\cite[Theorem~7.5.0.6]{lurieha}, as well as work of Baker--Richter and Rognes; see \cref{remark:work_of_baker_richter_rognes} below.

\begin{theorem}[\'{E}tale rigidity]
\label{theorem:dirac_etale_extensions_of_ek_rings}
If $R$ is an $\mathbf{E}_{k+1}$-algebra in spectra with $1 \leq k \leq \infty$, then taking homotopy groups
$$\xymatrix@C=10mm{
{ \operatorname{Alg}_{\mathbf{E}_k}(\operatorname{Mod}_R(\operatorname{Sp}))^{\text{{\rm \'{e}t}}} } \ar[r]^-{\pi_*} &
{ \operatorname{CAlg}(\operatorname{Mod}_{\pi_*(R)}(\operatorname{Ab}))^{\text{{\rm\'{e}t}}} } \cr
}$$
is an equivalence from the full subcategory of $\operatorname{Alg}_{\mathbf{E}_k}(\operatorname{Mod}_R(\operatorname{Sp}))$ spanned by the \'{e}tale algebras to the full subcategory of $\operatorname{CAlg}(\operatorname{Mod}_{\pi_*(R)}(\operatorname{Ab}))$ spanned by the \'{e}tale algebras.
\end{theorem}

\begin{remark}
\label{remark:work_of_baker_richter_rognes}If $k = \infty$, then \cref{theorem:dirac_etale_extensions_of_ek_rings} is a statement about commutative algebras in spectra, which we may divide into the following two statements:
\begin{enumerate}
\item If $B$ is an \'{e}tale $\pi_*(R)$-algebra, then there exists an \'{e}tale $R$-algebra $A$ and an isomorphism of $\pi_*(R)$-algebras $\pi_*(A) \simeq B$.
\item If $A$ and $A'$ is a map of \'{e}tale $R$-algebras, then the map
$$\xymatrix{
{ \operatorname{Map}(A,A') } \ar[r]^-{\pi_*} &
{ \operatorname{Map}(\pi_*(A),\pi_*(A')) } \cr
}$$
from the mapping anima of commutative algebras in $\operatorname{Mod}_R(\operatorname{Sp})$ to the mapping anima of commutative algebras in $\operatorname{Mod}_{\pi_*(R)}(\operatorname{Ab})$ is an equivalence.
\end{enumerate}
Under the additional assumption that the $\pi_*(R)$-module $B$ be projective, the first statement can be deduced from the work of Baker--Richter~\cite{baker2007realizability}, whereas the second statement can be deduced from the work of Rognes~\cite[Corollary~10.1.5]{rognes2008galois}. Both of these sources work with the more restrictive notion of a $G$-Galois extension in the sense of Rognes, but one can show that the obstruction groups that they consider vanish more generally for \'{e}tale maps in our sense.
\end{remark}

The rest of this section is devoted to the proof of \cref{theorem:dirac_etale_extensions_of_ek_rings}. In the proof of his rigidity theorem, Lurie first passes to the connective cover and then inducts his way up the Postnikov tower, using the theory of square-zero extensions, from the case of $\pi_0$. This strategy will not work here, since an \'{e}tale extension in our sense is not determined by its behavior on $\pi_0$.

Instead, we employ the obstruction theory of Goerss--Hopkins~\cite{goerss2005moduli} in the refined form of the second author and VanKoughnett~\cite{pstrkagowski2022abstract}. This, in turn, is based on the second author's synthetic deformation of the stable $\infty$-category $\operatorname{Mod}_A(\operatorname{Sp})$, which has the derived $\infty$-category $\mathcal{D}(\pi_*(A))$ as its special fiber. The Goerss--Hopkins obstruction theory then amounts to giving successive obstructions to extending an object from the special fiber $\mathcal{D}(\pi_*(A))$ to the entire deformation, and hence, to the generic fiber, which is identified with $\operatorname{Mod}_A(\operatorname{Sp})$. In the case that we consider, the relevant obstructions are controlled by the $\mathbf{E}_k$-cotangent complex. The standard reference for the latter is~\cite[Section~7.3.5]{lurieha}, but see also~\cite{francis2013tangent}.

\begin{proposition}
\label{proposition:etale_morphsisms_of_dirac_rings_have_vanishing_cotangent_complex}
Let $\phi \colon R \to A$ be an \'{e}tale map of Dirac rings. Considering $A$ as an $\mathbf{E}_{k}$-algebra with $2 \leq k \leq \infty$ in the symmetric monoidal derived $\infty$-category $\mathcal{D}(R)$ of graded $R$-modules, its $\mathbf{E}_{k}$-cotangent complex $L_{A/R}^{\mathbf{E}_k}$ is zero.
\end{proposition}

\begin{proof}
The map $\phi \colon R \to A$ is even, by \cref{theorem:equivalence_of_categories_between_etale_extensions_of_a_and_even_subring}, so it is the base-change of the map $\phi^{\operatorname{ev}} \colon R^{\operatorname{ev}} \to A^{\operatorname{ev}}$ along the inclusion $j \colon R^{\operatorname{ev}} \to R$. Since $\phi^{\operatorname{ev}}$ is \'{e}tale, and hence, flat, this is also the derived base-change. The $\mathbf{E}_k$-cotangent complex satisfies base-change along symmetric monoidal left adjoint functors between symmetric monoidal stable $\infty$-categories such as $j^* \simeq R \otimes_{R^{\operatorname{ev}}}^L-$; see~\cite[Lemma~7.6]{pstrkagowski2022abstract}.\footnote{\,The cited reference proves the base-change statement for the $\mathbf{E}_{\infty}$-cotangent complex, but the argument for the $\mathbf{E}_k$-cotangent complex is identical.} In the case at hand, we conclude that the base-change map
$$\xymatrix{
{ j^*(L_{A^{\operatorname{ev}}/R^{\operatorname{ev}}}^{\mathbf{E}_k}) } \ar[r] &
{ L_{A/R}^{\mathbf{E}_k} } \cr
}$$
is an equivalence, so we may assume without loss of generality that $\phi \colon R \to A$ is an \'{e}tale map of even Dirac rings.

If $R$ is even, let us denote by $\mathcal{D}_{\operatorname{ev}}(R)$ the full subcategory of the derived $\infty$-category of graded $R$-modules spanned by those objects, whose homology graded $R$-modules are even. The functor that to a graded $R$-module $M$ assigns its underlying ungraded $\smash{\widetilde{R}}$-module $\smash{\widetilde{M}}$ induces a symmetric monoidal (note that it is symmetric monoidal by the evenness assumption, otherwise it would merely be lax symmetric monoidal) left adjoint functor
\begin{equation}
\label{equation:underlying_module_on_derived_infty_cats}
\xymatrix@C=10mm{
{ \mathcal{D}_{\operatorname{ev}}(R) } \ar[r]^-{\phantom{j}\widetilde{(-)}\phantom{j}} &
{ \mathcal{D}(\widetilde{R}) } \cr
}
\end{equation}
to the derived $\infty$-category of $\widetilde{R}$. Therefore, we may again apply base-change for the $\mathbf{E}_k$-cotangent complex to conclude that
$$\widetilde{L_{A/R}^{\mathbf{E}_{k}}} \simeq L_{\widetilde{A}/\widetilde{R}}^{\mathbf{E}_k}.$$
Moreover, since $M \mapsto \widetilde{M}$ is exact, we have
$$H_n(\widetilde{X}) \simeq \widetilde{H_n(X)}$$
for all $X \in \mathcal{D}_{\operatorname{ev}}(R)$ and $n \in \mathbb{Z}$, so, in particular, the functor (\ref{equation:underlying_module_on_derived_infty_cats}) is conservative. Thus, it suffices to show that the $\mathbf{E}_k$-cotangent complex of the map of underlying rings is zero. 

So we let $\phi \colon R \to A$ be an \'{e}tale map of commutative rings, and proceed to show that its $\mathbf{E}_k$-cotangent complex is zero for $2 \leq k \leq \infty$. By~\cite[Theorem~7.3.5.1]{lurieha}, there is a canonical cofiber sequence
$$\xymatrix{
{ \int_{S^{k-1}} A } \ar[r] &
{ A } \ar[r] &
{ \Sigma^k L^{\mathbf{E}_{k}}_{A/R} } \cr
}$$
of $\mathbf{E}_{k}$-$A$-modules in $\mathcal{D}(R)$. The left term is the factorization homology of $A$ in this $\mathbf{E}_k$-monoidal $\infty$-category, and since $A$ is commutative, it follows from~\cite[7.3.5.6]{lurieha} that we have an equivalence
$$\textstyle{ \int_{S^{k-1}} A \simeq A^{\otimes_R S^{k-1}} }$$
with the colimit in $\operatorname{CAlg}(\mathcal{D}(R))$ of the constant diagram with value $A$ indexed by the anima $S^{k-1}$. We wish to show that the map
$$\xymatrix{
{ A^{\otimes_RS^{k-1}} } \ar[r] &
{ A^{\otimes_R1} \simeq A } \cr
}$$
induced by the unique map $S^{k-1} \to 1$ is an equivalence. For $2 \leq k < \infty$, we have
$$S^{k-1} \simeq \Sigma S^{k-2} \simeq \varinjlim\,(1 \leftarrow S^{k-2} \to 1),$$
and the functor that to $X$ assigns $A^{\otimes_RX}$ preserves colimits, since it is left adjoint to the functor $\operatorname{CAlg}(\mathcal{D}(R)) \to \mathcal{S}$ corepresented by $A$. For $k = 2$, we find that
$$A^{\otimes_RS^1} \simeq \varinjlim\,(A \leftarrow A \otimes_RA \to A) \simeq A \otimes_{A \otimes_RA}A \simeq \operatorname{HH}(A/R)$$
is the Hochschild homology of $A/R$, and since $\phi \colon R \to A$ is \'{e}tale, we have
$$\operatorname{HH}(A/R) \simeq A$$
by~\cite[Theorem~0.1]{gellerweibel}. For $k < 2 < \infty$, we find inductively that
$$A^{\otimes_RS^{k-1}} \simeq \varinjlim\,(A \leftarrow A^{\otimes_RS^{k-2}} \to A) \simeq \varinjlim\,(A \leftarrow A \to A) \simeq A,$$
where the middle equivalence is a consequence of the inductive hypothesis. Finally, in the case $k = \infty$, we use~\cite[Remark~7.3.5.6]{lurieha} to conclude that the canonical map
$$\xymatrix{
{ \varinjlim_k\, L_{A/R}^{\mathbf{E}_k} } \ar[r] &
{ L_{A/R}^{\mathbf{E}_{\infty}} } \cr
}$$
is an equivalence.
\end{proof}

We are now ready to prove the classification of Dirac \'{e}tale morphisms.

\begin{proof}[Proof of \cref{theorem:dirac_etale_extensions_of_ek_rings}]
Let $R$ be an $\mathbf{E}_{k+1}$-ring for $1 \leq k \leq \infty$. We must show that the functor $A \mapsto \pi_*(A)$ induces an equivalence between the $\infty$-categories of \'{e}tale $\mathbf{E}_k$-algebras in $\operatorname{Mod}_R(\operatorname{Sp})$ and \'{e}tale commutative algebras in $\operatorname{Mod}_{\pi_*(R)}(\operatorname{Ab})$. We will employ Goerss--Hopkins theory following \cite{pstrkagowski2022abstract}. This takes as an input an appropriate complete prestable $\infty$-category, which appears in \cite[Section~6.5]{patchkoria2021adams}, but we recall the construction here.

Let $\operatorname{Mod}_R^{\operatorname{ff}}(\operatorname{Sp}) \subset \operatorname{Mod}_R(\operatorname{Sp})$ be the full stable subcategory generated by the monoidal unit. It promotes canonically to an $\mathbf{E}_{k}$-monoidal subcategory. We say that a presheaf of anima $X \colon \operatorname{Mod}_R^{\operatorname{ff}}(\operatorname{Sp})^{\operatorname{op}} \to \mathcal{S}$ is a (connective) synthetic $R$-module if it preserves finite products, and we write
$$\operatorname{Syn}_R = \mathcal{P}_{\Sigma}(\operatorname{Mod}_R^{\operatorname{ff}}(\operatorname{Sp}))$$
for the $\infty$-category of (connective) synthetic $R$-modules. The $\infty$-category $\operatorname{Syn}_R$ is a complete Grothendick prestable $\infty$-category in the sense of~\cite[Section~C.1.4]{luriesag}, and it inherits an $\mathbf{E}_k$-monoidal structure from that of $\operatorname{Mod}_R^{\operatorname{ff}}(\operatorname{Sp})$ through Day convolution. An argument using Lawvere theories gives a canonical equivalence
$$\operatorname{Syn}_R^{\heartsuit} \simeq \operatorname{Mod}_{\pi_*(R)}(\operatorname{Ab})$$
between the heart of $\operatorname{Syn}_R$ and the abelian category of graded $\pi_*(R)$-modules.

The additive monoid of pairs $(i,j)$ of an integer $i$ and a half-integer $j$ such that the Chow degree $i - 2j$ is non-negative acts on $\operatorname{Syn}_R$ through left adjoint functors\footnote{\,We use motivic indexing, but the motivic weight $j$ is allowed to be a half-integer, as opposed to an integer. In terms of the indexing used in~\cite{pstrkagowski2022abstract}, we have $\Sigma^{i,j}X = \Sigma^{i-2j}X[2j]$.}
$$\xymatrix@C=10mm{
{ \operatorname{Syn}_R } \ar[r]^-{\Sigma^{i,j}} &
{ \operatorname{Syn}_R. } \cr
}$$
The generators $\Sigma^{1,0}$ and $\Sigma^{1,\nicefrac{1}{2}}$ of the action are given respectively by the suspension in the prestable $\infty$-category $\operatorname{Syn}_R$ and by the unique extension
$$\xymatrix@C=10mm{
{ \operatorname{Mod}_R^{\operatorname{ff}}(\operatorname{Sp}) } \ar[r]^-{\Sigma} \ar@<-.2ex>[d]^-{\nu} &
{ \operatorname{Mod}_R^{\operatorname{ff}}(\operatorname{Sp}) } \ar@<-.2ex>[d]^-{\nu} \cr
{ \operatorname{Syn}_R } \ar[r]^-{\Sigma^{1,\nicefrac{1}{2}}} &
{ \operatorname{Syn}_R } \cr
}$$
of the suspension in the stable $\infty$-category $\operatorname{Mod}_R^{\operatorname{ff}}(\operatorname{Sp})$ along the Yoneda embedding. Thus, the functor $\Sigma^{1,\nicefrac{1}{2}}$ is an auto-equivalence, whereas the functor $\Sigma^{1,0}$ is not. The monoidal unit is $R = \nu(R)$ and the colimit-interchange map $\Sigma \nu(R) \to \nu(\Sigma R)$ is a map $\Sigma^{1,0}R \to \Sigma^{1,\nicefrac{1}{2}}R$, or equivalently, a map\footnote{\,The map $\tau$ is a tensor square-root of the map denoted $\tau$ in motivic homotopy theory.}
$$\xymatrix{
{ \Sigma^{0,-\nicefrac{1}{2}}R } \ar[r]^-{\tau} &
{ R. } \cr
}$$
It gives $R$ the structure of a shift algebra in the sense of \cite[Definition 2.13]{pstrkagowski2022abstract}, and the Yoneda embedding induces a canonical $\mathbf{E}_k$-monoidal equivalence
$$\xymatrix{
{ \operatorname{Mod}_R(\operatorname{Sp}) } \ar[r] &
{ \operatorname{Syn}_R^{\operatorname{per}} \subset \operatorname{Syn}_R } \cr
}$$
from the $\mathbf{E}_k$-monoidal $\infty$-category of $R$-modules in spectra to the full subcategory of the $\mathbf{E}_k$-monoidal $\infty$-category of synthetic $R$-modules spanned by the synthetic $R$-modules $X$ that are $\tau$-periodic in the sense that the map
$$\xymatrix{
{ \Sigma^{0,-\nicefrac{1}{2}}X } \ar[r]^-{\tau} &
{ X } \cr
}$$
is a $1$-connective cover in the prestable $\infty$-category $\operatorname{Syn}_R$.

By \cite[Section~5]{pstrkagowski2022abstract} (which assumes the background prestable $\infty$-category $\operatorname{Syn}_R$ to be symmetric monoidal, but $\mathbf{E}_{k}$-monoidal is sufficient for our purposes), we have a tower of $\infty$-categories of potential $n$-stages
$$\xymatrix@C=6mm{
{ \mathcal{M}_{\infty} } \ar[r] &
{ \dots } \ar[r] &
{ \mathcal{M}_{n} } \ar[r] &
{ \dots } \ar[r] &
{ \mathcal{M}_{1} } \ar[r] &
{ \mathcal{M}_{0} } \cr
}$$
with the following properties:
\begin{enumerate}
\item [(1)]There is a canonical equivalence $\mathcal{M}_{0} \simeq \operatorname{Alg}_{\mathbf{E}_k}(\operatorname{Mod}_{\pi_*(R)}(\operatorname{Ab}))$.
\item [(2)]There is a canonical equivalence $\mathcal{M}_{\infty} \simeq \operatorname{Alg}_{\mathbf{E}_{k}}(\operatorname{Mod}_R(\operatorname{Sp}))$.
\end{enumerate}
In particular, an object $X \in \mathcal{M}_n$ with $0 \leq n \leq \infty$ determines a point in $\mathcal{M}_0$, and hence, by~(1) an $\mathbf{E}_k$-algebra in graded $\pi_*(R)$-modules that we denote by $\pi_*(X)$. If $k = \infty$, then $\pi_*(X)$ agrees, under the equivalence~(2), with the homotopy groups of the $\mathbf{E}_k$-algebra in $R$-modules in spectra corresponding to $X \in \mathcal{M}_{\infty}$; whence, the choice of notation. The following additional properties describe how these different $\infty$-categories relate to each other: 
\begin{enumerate}
\item [(3)]The canonical map of $\infty$-categories $\mathcal{M}_{\infty} \to \varprojlim_n \mathcal{M}_n$ is an equivalence.
\item [(4)]Given an $X \in \mathcal{M}_{n-1}$ with $1 \leq n < \infty$, there is an obstruction
$$\theta_X \in \operatorname{Ext}^2_{\pi_*(X)}(L^{\mathbf{E}_k}_{\pi_*(X)/ \pi_*(R)}, \pi_*(X)(\nicefrac{n}{2}))$$
where the $\operatorname{Ext}$-groups are calculated in the stable $\infty$-category $$\operatorname{Mod}_{\pi_*(X)}^{\mathbf{E}_k}(\mathcal{D}(\pi_*(R)))$$
of $\mathbf{E}_k$-$\pi_*(X)$-modules in the symmetric monoidal derived $\infty$-category of graded $\pi_*(R)$-modules, and where the twist ``$(\nicefrac{n}{2})$'' is the Serre twist of Definition~\ref{definition:Serre_twist}. The obstruction $\theta_X$ vanishes if and only the given object $X \in \mathcal{M}_{n-1}$ can be lifted to an object of $\mathcal{M}_n$.
\item [(5)]Given $X, Y \in \mathcal{M}_n$ and a map $f' \colon X' \to Y'$ between their images in $\mathcal{M}_{n-1}$, there is an obstruction 
$$\theta_{f'} \in \operatorname{Ext}^1_{\pi_*(X)}(L^{\mathbf{E}_k}_{\pi_*(X)/ \pi_*(R)}, \pi_*(Y)(\nicefrac{n}{2})).$$
The obstruction $\theta_{f'}$ vanishes if and only if the map $f' \colon X' \to Y'$ in $\mathcal{M}_{n-1}$ lifts to a map $f \colon X \to Y$ in $\mathcal{M}_n$. Suppose that a lift $f \colon X \to Y$ exists and let
$$F_{f'} = \operatorname{Map}_{\mathcal{M}_{n}}(X,Y) \times_{\operatorname{Map}_{\mathcal{M}_{n-1}}(X',Y')} \{ f' \}$$
be the anima of lifts. In this situation, there is a canonical bijection
$$\pi_i(F_{f'},f) \simeq \operatorname{Ext}^{-i}_{\pi_*(X)}(L^{\mathbf{E}_k}_{\pi_*(X)/\pi_*(R)}, \pi_*(Y)(\nicefrac{n}{2})).$$
\end{enumerate}
Let us write $\mathcal{M}_{0}^{\operatorname{\acute{e}t}} \subset \mathcal{M}_0$ for the full subcategory of potential $0$-stages that correspond to Dirac rings that are \'{e}tale over $\pi_*(R)$, so that we have an induced subtower 
\begin{equation}
\label{equation:etale_goerss_hopkins_tower}
\xymatrix@C=5mm{
{ \mathcal{M}_{\infty}^{\operatorname{\acute{e}t}} } \ar[r] &
{ \dots } \ar[r] &
{ \mathcal{M}_n^{\operatorname{\acute{e}t}} } \ar[r] &
{ \dots } \ar[r] &
{ \mathcal{M}_1^{\operatorname{\acute{e}t}} } \ar[r] &
{ \mathcal{M}_0^{\operatorname{\acute{e}t}} } \cr
}
\end{equation}
with $\mathcal{M}_n^{\operatorname{\acute{e}t}} = \mathcal{M}_n \times_{\mathcal{M}_0} \mathcal{M}_0^{\operatorname{\acute{e}t}}$. We claim that all maps in~(\ref{equation:etale_goerss_hopkins_tower}) are equivalences of $\infty$-categories. This will follow immediately from the properties (3)--(5) above, if we prove that the abelian groups in which the various obstructions live all are zero. If $2 \leq k \leq \infty$, then the map $\pi_*(R) \to \pi_*(X)$ induced by the unit is an \'{e}tale map of Dirac rings, and hence, it follows from \cref{proposition:etale_morphsisms_of_dirac_rings_have_vanishing_cotangent_complex} that the $\mathbf{E}_k$-cotangent complex $L^{\mathbf{E}_k}_{\pi_*(X)/ \pi_*(R)}$ is zero.

If $k = 1$, then we need to show the obstruction groups vanish, since the cotangent complex may not be zero. By \cite[Theorem~7.3.5.1]{lurieha}, we have an identification 
\begin{equation}
\label{equation:e1_cotangent_complex_of_an_etale_map}
L^{\mathbf{E}_{1}}_{\pi_*(X)/\pi_*(R)} \simeq I,
\end{equation}
where $I$ is the kernel of the multiplication $\pi_*(X) \otimes_{\pi_*(R)} \pi_*(X) \to \pi_*(X)$, considered as an object of the heart of the derived $\infty$-category of $\pi_*(X)$-bimodules. Here we use the fact that $\pi_*(R) \to \pi_*(X)$ is flat, which implies that the tensor products in the derived $\infty$-category and in the heart agree. Now, in the situation of~(5) (and of~(4) in the special case of $\pi_*(Y) = \pi_*(X)$), the map $\pi_*(X) \to \pi_*(Y)$ is a map of Dirac rings by assumption. Therefore, the left and right multiplication coincide and the bimodule structure on $\pi_*(Y)$ is obtained by restriction of scalars along the multiplication map $\pi_*(X) \otimes_{\pi_*(R)} \pi_*(X) \to \pi_*(X)$. But this map is flat by \cref{proposition:unramified_finite_type_morphisms_have_flat_multiplication} and using (\ref{equation:e1_cotangent_complex_of_an_etale_map}), it follows that we can rewrite the relevant obstruction groups as
$$\operatorname{Ext}_{\pi_*(X) \otimes_{\pi_*(R)} \pi_*(X))}^i(I,\pi_*(Y)(\nicefrac{n}{2})) \simeq \operatorname{Ext}_{\pi_*(X)}^i(I/I^2, \pi_*(Y)(\nicefrac{n}{2})).$$
But $I/I^2 \simeq \Omega_{\pi_*(X)/\pi_*(R)}^1 \simeq 0$,
since $\pi_*(R) \to \pi_*(X)$ is an \'{e}tale map of Dirac rings, and hence, obstruction groups all are zero, also for $k = 1$. This proves that all maps in~(\ref{equation:etale_goerss_hopkins_tower}) are equivalences of $\infty$-categories, as claimed.

Finally, properties~(1)--(2) identify the functor $\mathcal{M}_{\infty}^{\operatorname{\acute{e}t}} \to \mathcal{M}_{0}^{\operatorname{\acute{e}t}}$ with the functor
$$\xymatrix@C=10mm{
{ \operatorname{Alg}_{\mathbf{E}_k}(\operatorname{Mod}_R(\operatorname{Sp}))^{\text{{\rm \'{e}t}}} } \ar[r]^-{\pi_*} &
{ \operatorname{CAlg}(\operatorname{Mod}_{\pi_*(R)}(\operatorname{Ab}))^{\text{{\rm\'{e}t}}} } \cr
}$$
in the statement of theorem, so we conclude that this functor is an equivalence of $\infty$-categories, as stated.
\end{proof}

\appendix

\section{Geometric realization in an $n$-category}

Let $n \geq -1$ be an integer. We say that an $\infty$-category $\mathcal{C}$ is an $n$-category if the mapping anima $\operatorname{Map}(y,x)$ are $(n-1)$-truncated for all $x,y \in \mathcal{C}$. We let $\Delta_{\leq n} \subset \Delta$ be the full subcategory spanned by the objects $[m]$ with $0 \leq m \leq n$ and record the following fundamental property of the category $\Delta$.

\begin{proposition}\label{prop:totalization}Let $n \geq -1$ be an integer, and let $\mathcal{C}$ be an $n$-category.
\begin{enumerate}
\item[{\rm (1)}]If $\mathcal{C}$ admits finite colimits, then it admits geometric realizations, and for every simplicial object $X \colon \Delta^{\operatorname{op}} \to \mathcal{C}$, the canonical map
$$\xymatrix{
{ \varinjlim_{\Delta_{\leq n}^{\operatorname{op}}} X|_{\Delta_{\leq n}^{\operatorname{op}}} } \ar[r] &
{ \varinjlim_{\Delta^{\operatorname{op}}} X } \cr
}$$
is an equivalence.
\item[{\rm (2)}]If $\mathcal{C}$ admits finite limits, then it admits totalizations, and for every cosimplicial object $X \colon \Delta \to \mathcal{C}$, the canonical map
$$\xymatrix{
{ \varprojlim_{\Delta} X } \ar[r] &
{ \varprojlim_{\Delta_{\leq n}} X|_{\Delta_{\leq n}} } \cr
}$$
is an equivalence.
\end{enumerate}
\end{proposition}

\begin{proof}
An $\infty$-category $\mathcal{C}$ is an $n$-category if and only if its opposite $\infty$-category $\mathcal{C}^{\operatorname{op}}$ is an $n$-category, so the statements~(1) and~(2) are equivalent. We will prove~(1), following~\cite[Lemma~14.8]{luriedagI}. By passing to a larger universe, if necessary, we may assume that $\mathcal{C}$ is small. By~\cite[Theorem~5.5.1.1]{luriehtt}, the $\infty$-category $\mathcal{D} \simeq \operatorname{Ind}(\mathcal{C})$ is presentable, and the Yoneda embedding
$$\xymatrix{
{ \mathcal{C} } \ar[r]^-{j} &
{ \mathcal{D} \simeq \operatorname{Ind}(\mathcal{C}) } \cr
}$$
is fully faithful; see~\cite[Remark~5.3.5.2]{luriehtt}. The $\infty$-category $\mathcal{D}$ is again an $n$-category, because the full subcategory $\tau_{\leq n-1}\mathcal{S} \subset \mathcal{S}$ spanned by the $(n-1)$-truncated anima is closed under filtered colimits, and the Yoneda embedding $j \colon \mathcal{C} \to \mathcal{D}$ preserves finite colimits by~\cite[Proposition~5.3.5.14]{luriehtt}. So it preserves $\smash{ \Delta_{\leq n}^{\operatorname{op}} }$-indexed colimits, since $\smash{ \Delta_{\leq n}^{\operatorname{op}} }$ is $\varinjlim$-equivalent to a finite $\infty$-category by~\cite[Lemma~1.2.4.17]{lurieha}.

We now claim that it suffices to prove the result for $\mathcal{D}$. Indeed, suppose that for a diagram $X \colon \Delta^{\operatorname{op}} \to \mathcal{C}$, the induced map
$$\xymatrix{
{ \varinjlim_{\Delta_{\leq n}^{\operatorname{op}}} (j \circ X)|_{\Delta_{\leq n}^{\operatorname{op}}} } \ar[r] &
{ \varinjlim_{\Delta^{\operatorname{op}}} (j \circ X) } \cr
}$$
is an equivalence. The left-hand side is in the essential image of $j$, by the assumption that $\mathcal{C}$ finite small colimits, and therefore, so is the right-hand side. Since $j$ is fully faithful, we conclude that $X \colon \Delta^{\operatorname{op}} \to \mathcal{C}$ admits a colimit and that, in the diagram
$$\xymatrix{
{ j(\varinjlim_{\Delta_{\leq n}^{\operatorname{op}}} X|_{\Delta_{\leq n}^{\operatorname{op}}}) } \ar[r] \ar[d] &
{ j(\varinjlim_{\Delta^{\operatorname{op}}} X)\phantom{,} } \ar[d] \cr
{ \varinjlim_{\Delta_{\leq n}^{\operatorname{op}}} (j \circ X)|_{\Delta_{\leq n}^{\operatorname{op}}} } \ar[r] &
{ \varinjlim_{\Delta^{\operatorname{op}}} (j \circ X), } \cr
}$$
the vertical maps are equivalences. Hence, the top horizontal map is an equivalence, and since $j$ is fully faithful, so is the map
$$\xymatrix{
{ \varinjlim_{\Delta_{\leq n}^{\operatorname{op}}} X|_{\Delta_{\leq n}^{\operatorname{op}}} } \ar[r] &
{ \varinjlim_{\Delta^{\operatorname{op}}} X } \cr
}$$
as desired. This proves the claim.

Finally, since $\mathcal{D}$ is an $n$-category, the localization
$$\xymatrix@C=10mm{
{ \mathcal{P}(\mathcal{C}) } \ar@<.7ex>[r]^-{L} &
{ \operatorname{Ind}(\mathcal{C}) \simeq \mathcal{D} } \ar@<.7ex>[l]^-{\iota} \cr
}$$
factors uniquely through a localization
$$\xymatrix@C=10mm{
{ \tau_{\leq n-1}\mathcal{P}(\mathcal{C}) } \ar@<.7ex>[r]^-{L'} &
{ \mathcal{D} } \ar@<.7ex>[l]^-{\iota'} \cr
}$$
of the $(n-1)$-truncation of $\mathcal{P}(\mathcal{C})$. Hence, it suffices to prove that for every simplicial object $Y \colon \Delta^{\operatorname{op}} \to \mathcal{P}(\mathcal{C})$ in the $\infty$-topos $\mathcal{P}(\mathcal{C})$, the canonical map
$$\xymatrix{
{ \varinjlim_{\Delta_{\leq n}^{\operatorname{op}}} Y|_{\Delta_{\leq n}^{\operatorname{op}}} } \ar[r] &
{ \varinjlim_{\Delta^{\operatorname{op}}} Y } \cr
}$$
is $n$-connective. But we may identify this map with the map
$$\xymatrix{
{ \varinjlim_{\Delta_{\leq n}^{\operatorname{op}}} i^*(Y) \simeq \varinjlim_{\Delta^{\operatorname{op}}} i_!i^*(Y) } \ar[r] &
{ \varinjlim_{\Delta^{\operatorname{op}}} Y } \cr
}$$
induced by the counit of the adjunction given by left Kan extension and restriction along the inclusion $\smash{ i \colon \Delta_{\leq n}^{\operatorname{op}} \to \Delta^{\operatorname{op}} }$, and \cite[Proposition~6.5.3.10]{luriehtt} shows that this map is $n$-connective, because $i$ is fully faithful.
\end{proof}

\section{Conservativity and limits}

The following result is certainly well-known, but we have not been able to find a suitable reference, so we take this opportunity to give the proof here. 

\begin{lemma}
\label{lem:conservative}
If
$\mathcal{C} \colon K^{\operatorname{op}} \to \operatorname{Cat}_{\infty}$ is a diagram of
$\infty$-categories indexed by a small simplicial set $K^{\operatorname{op}}$, then
the restriction
$$\xymatrix{
{ \lim_{K^{\operatorname{op}}}\mathcal{C} } \ar[r]^-{i^*} &
{ \prod_{k \in K_0} \mathcal{C}(k) } \cr
}$$
along the inclusion $i \colon K_0 \to K$ of the set of vertices
is a conservative functor.
\end{lemma}

\begin{proof}We use the description of limits in $\operatorname{Cat}_{\infty}$ in terms of marked simplicial sets developed in~\cite[Chapter~3]{luriehtt}. So let $p \colon X \to K$ be a cartesian  fibration of simplicial sets classified by $\mathcal{C} \colon K^{\operatorname{op}} \to \operatorname{Cat}_{\infty}$, and let $p' \colon X' \to K_0$ be the base-change of $p$ along $i$ in the 1-category of simplicial sets. It follows from~\cite[Corollary~3.3.3.2]{luriehtt} that we may identify the map in the statement with the restriction map
$$\xymatrix{
{ \operatorname{Map}_K^{\flat}(K^{\sharp},X^{\natural}) } \ar[r]^-{i^*} &
{ \operatorname{Map}_{K_0}^{\flat}(K_0^{\sharp},X'{}^{\natural}) } \cr
}$$
from the simplicial set of cartesian sections of $p$ to that of cartesian sections of $p'$. An edge in the domain is a map $\sigma \colon (\Delta^1)^{\flat} \times K^{\sharp} \to X^{\natural}$ of marked simplicial sets over $K^{\sharp}$, and its image in the target is the map $i^*(\sigma) \colon (\Delta^1)^{\flat} \times K_0^{\sharp} \to X'{}^{\natural}$ obtained from $\sigma$ by restriction along $i$. By~\cite[Remark~3.1.3.1]{luriehtt}, the edge $\sigma$ is invertible if and only if it factors through $(\Delta^1)^{\flat} \times K^{\sharp} \to (\Delta^1)^{\sharp} \times K^{\sharp}$. Similarly, the edge $i^*(\sigma)$ is invertible if and only if it factors through $(\Delta^1)^{\flat} \times K_0^{\sharp} \to (\Delta^1)^{\sharp} \times K_0^{\sharp}$.

We now assume that $i^*(\sigma)$ is invertible and prove that $\sigma$ is invertible. An edge in $\Delta^1 \times K$ is a pair $(a,b)$ of an edge $a$ in $\Delta^1$ and an edge $b$ in $K$, and we must show that $\sigma \colon \Delta^1 \times K \to X$ maps all such edges to $p$-cartesian edges in $X$. We already know that $\sigma$ maps all edges $(a,b)$ with $a$ degenerate to $p$-cartesian edges in $X$, and the assumption that $i^*(\sigma)$ be invertible implies that $\sigma$ also maps all edges $(a,b)$ with $b$ degenerate to $p$-cartesian edges in $X$. So $\sigma$ maps all edges $(a,b)$, where either $a$ or $b$ or both are degenerate, to $p$-cartesian edges in $X$. But every edge in $\Delta^1 \times K$ is a composition of edges of this kind, and by~\cite[Proposition~2.4.1.7]{luriehtt}, a composition of two $p$-cartesian edges in $X$ is again a $p$-cartesian edge in $X$. So we conclude that $\sigma$ maps all edges in $\Delta^1 \times K$ to $p$-cartesian edges in $X$. 
\end{proof}

\providecommand{\bysame}{\leavevmode\hbox to3em{\hrulefill}\thinspace}
\providecommand{\MR}{\relax\ifhmode\unskip\space\fi MR }
\providecommand{\MRhref}[2]{%
  \href{http://www.ams.org/mathscinet-getitem?mr=#1}{#2}
}
\providecommand{\href}[2]{#2}

\end{document}